%% file: main.tex
\documentclass{article}

\input{extras/preamble}

\begin{document}

\maketitle

\begin{abstract}
\noindent Using the pro\'etale site, we construct models for the continuous actions of the Morava stabiliser group on Morava E-theory, its \inftycat{} of $K(h)$-local modules, and its Picard spectrum. For the two sheaves of spectra, we evaluate the resulting descent spectral sequences: these can be thought of as homotopy fixed point spectral sequences for the profinite Galois extension $L_{K(h)} \mathbb S \to E_h$. We show that the descent spectral sequence for the Morava E-theory sheaf is the $K(h)$-local $E_h$-Adams spectral sequence. The spectral sequence for the sheaf of Picard spectra is closely related to one recently defined by Heard; our formalism allows us to compare many differentials with the $K(h)$-local $E_h$-Adams spectral sequence, and isolate the exotic Picard elements in the $0$-stem. In particular, we show how this recovers the computation due to Hopkins, Mahowald and Sadofsky of the group $\Pic_1$ at the more difficult prime two. We also show this gives a bound on the Brauer group of $K(h)$-local spectra, and compute this bound at height one.
\end{abstract}

\tableofcontents

\section{Introduction}
\label{sec:intro}

\subfile{Sections/Intro}

\section{The continuous action on Morava E-theory}
\label{sec:Morava E-theory}

\subfile{Sections/Morava_E-theory}

\section{Descent for modules and the Picard spectrum}
\label{sec:descent theory}

In the previous section, we showed that Morava E-theory defines a hypercomplete sheaf of spectra $\mathcal E$ on the pro\'etale classifying site of $\G$. Our next aim is to improve this to a statement about its module $\infty$-category and therefore its Picard spectrum. The main result of this section is the construction of a hypercomplete sheaf of connective spectra $\pic(\mathcal E)$, with global sections $\Gamma \pic(\mathcal E) = \pic_h \coloneqq \pic(\mathcal Sp_{\K})$. The functoriality of the construction via the pro\'etale site will allow us to compare the resulting descent spectral sequence to the $\K$-local $\E$-Adams spectral sequence, including differentials.

\subsection{Descent for \texorpdfstring{$K(n)$}{K}-local module categories} \label{sec:modules}

\subfile{Sections/Modules}

\subsection{The Picard spectrum as a pro\'etale spectrum} \label{sec:picard sheaf}

\subfile{Sections/Picard_sheaf}

\section{Picard group computations}
\label{sec:n=1}

\subfile{Sections/n=1}


\section{Brauer groups} \label{sec:brauer}

\subfile{Sections/Brauer}

\appendix
\section{Appendix: D\'ecalage results}
\label{app:decalage}

\subfile{Sections/zDecalage}

\section{Appendix: The Adams spectral sequence at height one} 
\label{app:ASS n=1}

In this appendix, we compute the $K(1)$-local $E_1$-Adams spectral sequences at all primes. The results are again well-known, and implicit in the computations of \cite{ravenel_localisation, periodic-phenomena}. At the prime two we give a different perspective to \cite{beaudry-goerss-henn_splitting} on the computation, which illustrates how one can make use of the Postnikov tower of a sheaf of spectra on $\proet \G$: this approach may be useful in higher height examples that are more computationally challenging. In particular, we use the finite resolution of the $K(1)$-local sphere, and as such our only input is knowledge of the HFPSS for the conjugation action on $KU$, displayed in \cref{fig:HFPSS for KO to KU}.

\subsection{Odd primes}

\subfile{Sections/zASSodd}

\subsection{\texorpdfstring{$p=2$}{p=2}}

\subfile{Sections/zASSp=2}

\clearpage
\nocite{*}
\addcontentsline{toc}{section}{References}
\printbibliography

\end{document}

%% file: extras/preamble.tex
\usepackage{extras/arxiv}

\pagecolor{white}

\newcommand{\pic}{\mathfrak{pic}}
\newcommand{\gll}{\mathfrak{gl}_1}
\newcommand{\GL}{\mathrm{GL}_1}
\newcommand{\Pic}{\operatorname{Pic}}
\newcommand{\Picc}{\mathfrak{Pic}}
\newcommand{\Br}{\operatorname{Br}}

\newcommand{\Brr}{\mathfrak{Br}}

\newcommand{\Az}{\operatorname{Az}}
\newcommand{\Azz}{\mathfrak{Az}}
\newcommand{\E}{\mathbf{E}}
\newcommand{\K}{\mathbf{K}}
\newcommand{\G}{\mathbb{G}}
\newcommand{\Shh}[2][1]{\widehat{\mathcal{S}h}(#2, #1)}
\newcommand{\Sh}[2][1]{\mathcal{S}h(#2, #1)}

\newcommand{\Modk}[1]{\catname{Mod}_{#1, \K}}

\newcommand{\Prsm}{{\Pr}^{L, \operatorname{smon}}}

\newcommand{\et}[1]{B{#1}_{\operatorname{\acute{e}t}}}
\newcommand{\proet}[1]{B{#1}_{\operatorname{pro\acute{e}t}}}
\newcommand{\Free}[1]{\mathrm{Free}_{#1}}
\newcommand{\Proj}[1]{\mathrm{Proj}_{#1}}
\newcommand{\Cond}{\mathrm{Cond}}
\def\gl1{\mathfrak{gl}_1}
\newcommand{\Algk}[1]{\operatorname{Alg}_{#1, \K}}
\newcommand{\Catk}[1]{\operatorname{Cat}_{#1, \K}}
\def\complete{^{\mathrm{cpl}}}
\def\Lcomplete{^{\heartsuit \mathrm{cpl}}}
\newcommand{\Pro}{\mathrm{Pro}}

\def\Thotimes{\mathrm{Thick}^\otimes}
\newcommand{\Loccon}[2][]{\mathrm{LC}_{#1} \left(#2\right)}
\def\cd{\mathrm{cd}}
\def\vcd{\mathrm{vcd}}
\def\Pyk{\mathrm{Pyk}}

\loadglsentries{zzGlossary}

\usepackage{subfiles}

\title{Picard and Brauer groups of \texorpdfstring{$K(n)$}{K(n)}-local spectra\\via profinite Galois descent}
\author{Itamar Mor}

%% file: zzGlossary.tex
\newglossaryentry{Sphere}{
    name=\ensuremath{\bm 1_A},
    description={$A$-local sphere spectrum},
}

\newglossaryentry{E}{
    name=\ensuremath{\E},
    description={Morava E-theory, $E_n = E(\F[p^n], \Gamma_n)$},
}
        
\newglossaryentry{K}{
    name=\ensuremath{\K},
    description={Morava K-theory, $K(n) = E_n/(p, v_1, \dots, v_{n-1})$},
}


\newglossaryentry{G}{
    name=\ensuremath{\G},
    description={Extended Morava stabiliser group $\G_n = \mathbb S_n \rtimes \mathrm{Gal}(\F[p^n]/\F[p]) = \mathrm{Aut}(\Gamma_n, \F[p^n])$},
}


\newglossaryentry{In}{
    name = \ensuremath{I_n},
    description = {Maximal ideal $(p, v_1, \dots, v_{n-1}) \subset \pi_* E_n$}
}

\newglossaryentry{et}{
    name=\ensuremath{\et G},
    description={\'Etale classifying site of a profinite group $G$},
}

\newglossaryentry{proet}{
    name=\ensuremath{\proet G},
    description={Pro\'etale classifying site of a profinite group $G$},
}

\newglossaryentry{Free}{
    name=\ensuremath{\Free G},
    description={Subsite of free $G$-sets in $\proet G$},
}

\newglossaryentry{Cont}{
    name=\ensuremath{\Cont_G(T, T')},
    description = {Set of continuous $G$-equivariant maps between profinite $G$-sets}
}

\newglossaryentry{mathcal E delta}{
    name=\ensuremath{\mathcal E^\delta},
    description={\'Etale Morava E-theory sheaf}
}

\newglossaryentry{mathcal E}{
    name=\ensuremath{\mathcal E},
    description={Pro\'etale Morava E-theory sheaf}
}

\newglossaryentry{Modk}{
    name=\ensuremath{\Modk{A}},
    description={$K(n)$-local module \inftycat{}, $\Mod_A(\mathcal Sp_{K(n)}) = L_{K(n)} \Mod_A$}
}

\newglossaryentry{Pic}{
    name=\ensuremath{\Pic (-) \,(\Picc(-), \pic(-))},
    description={Picard group (resp. space, spectrum)}
}





\newglossaryentry{Picn}{
    name=\ensuremath{\Pic_n \, (\Picc_n, \,  \pic_n)},
    description={Picard group (resp. space, spectrum) of $K(n)$-local category, $\Pic (\mathcal Sp_{K(n)} )$},
}

\newglossaryentry{Picalg}{
    name=\ensuremath{\Pic_n^{\mathrm{alg}} \, (\Picc_n^{\mathrm{alg}})},
    description={Algebraic Picard group (resp. groupoid) of $K(n)$-local category, $\Pic (\Mod_{\pi_* E_n}^{\G_n} )$},
}

\newglossaryentry{Mod wedge}{
    name=\ensuremath{\Mod_R\Lcomplete},
    description={Category of $L$-complete discrete modules over a (classical) local ring $R$},
}

\newglossaryentry{Mod Lwedge}{
    name=\ensuremath{\Mod_R\complete},
    description={\inftycat{} of complete modules (in spectra) over a local ring $R$},
}
\newglossaryentry{otimes}{
    name=\ensuremath{\otimes},
    description={$K(n)$-localised smash product of spectra, $L_{K(n)} (- \wedge -)$},
}

\newglossaryentry{hat otimes}{
    name=\ensuremath{\widehat \otimes},
    description={Tensor product on $\Mod_R\Lcomplete$, $L_0 (- \otimes_R -)$},
}

\newglossaryentry{hat otimes L}{
    name=\ensuremath{\widehat{\otimes}^{\mathbb L}},
    description={Tensor product on $\Mod_R\complete$, $L (- \otimes_R^{\mathbb L} -)$},
}

\newglossaryentry{E vee}{
    name=\ensuremath{\E_*^\vee (-)},
    description={Completed E-homology functor $\pi_* L_{K(n)} (E_n \wedge - )$},
}

%% file: Sections/Intro.tex
In \cite{hopkins-mahowald-sadofsky}, Hopkins, Mahowald and Sadofsky study the \emph{Picard group} of a symmetric monoidal category: by definition, this is the group of isomorphism classes of invertible objects with respect to the monoidal product. This is a notion that goes back much further, and gives a useful invariant of a ring or scheme. Its particular relevance to homotopy theory comes from the observation that if the category $\mathcal C$ is a Brown category (for example, $\mathcal C$ might be the homotopy category of a compactly-generated stable \inftycat{}), then the representability theorem applies and shows that the Picard group $\Pic(\mathcal C)$ classifies homological automorphisms of $\mathcal C$, each of these being of the form $T \otimes (-)$ for some invertible object. The objective of \opcit\ is to develop techniques for studying Picard groups in some examples coming from chromatic homotopy: the main theorem is the computation of the Picard group of $K(1)$-local spectra at all primes, where $K(1)$ is Morava K-theory at height one.

The aim of this project is to give a new proof of these computations using \emph{Galois descent}, inspired by the formalism developed in \cite{ms}. We will write $\Pic_h$ for the Picard group of $K(h)$-local spectra. There are still many open questions regarding these groups: for example, it is unknown if they are finitely generated as modules over $\Z[p]$. The question of computing $\Pic_2$ has been studied by many authors (for example \cite{kamiya-shimomura, karamanov,ghmr_pic}); using recent work at the prime 2 \cite{bbghps_ass}, our results give a new potential approach to the computation of $\Pic_2$ in \cite{bbghps_pic}. Following \cite{gl}, we also extend these techniques to the \emph{Brauer group} of $\mathcal Sp_{K(h)}$, giving a cohomological approach to these, which allows us to bound the size of the group of $K(1)$-local Azumaya algebras trivialised over $E_1 = KU_p$, at all primes.

The notion of Galois descent in algebra is very classical, and says that if $A \to B$ is a Galois extension of rings, then $\Mod_A$ can be recovered as the category of descent data in $\Mod_B$: in particular, invertible $A$-modules can be recovered as invertible $B$-modules $M$ equipped with isomorphisms $\psi_g: M \cong g^* M$ for each $g \in \operatorname{Gal}(B/A)$, subject to a cocycle condition. This gives an effective way to compute Picard groups and other invariants. One can try to play the same game in higher algebra, and use descent techniques to get a handle on the groups $\Pic_h$. Fundamental to this approach is the notion of a Galois extension of commutative ring spectra, as set down in \cite{rognes}: this is a direct generalisation of the classical axioms in \cite{auslander-goldman}. Given a finite $G$-Galois extension $A \to B$ which is \emph{faithful}, the analogous descent statement (due to \cite{meier_thesis, gl}) is that the canonical functor
\begin{equation} \label{intro:Galois descent}
    \Mod_A \to (\Mod_B)^{hG} \coloneqq \lim \left( \cosimp[\Mod_B]{\prod_G \Mod_B} \right)
\end{equation}
is an equivalence of symmetric monoidal \inftycats{}. Taking the Picard \emph{spectrum} of a symmetric monoidal \inftycat{} preserves (homotopy) limits, and therefore any such Galois extension gives rise to an equivalence $\pic(\Mod_A) \simeq \tau_{\geq 0} (\pic(\Mod_B)^{hG})$. In particular one gets a homotopy fixed points spectral sequence (hereafter HFPSS), whose 0-stem converges to $\Pic(\Mod_A)$. This technique has proved very fruitful in Picard group computations: for example, the Picard groups of $KO$ and $Tmf$ are computed in \cite{ms}, and the Picard group of the higher real K-theories $EO(n)$ in \cite{heard-mathew-stojanoska}. In each case the starting point is a theorem of Baker and Richter \cite{baker-richter}, which says that the Picard group of an even-periodic ring spectrum $E$ with $\pi_0 E$ regular Noetherian is a $\Z/2$-extension of the Picard group of the ring $\pi_0 E$. One can for example study the action of the Morava stabiliser group $\G_h$ on Morava E-theory:

\begin{thm}[\cite{devinatz-hopkins,rognes}]
Write $E_h$ for height $h$ Morava E-theory at the (implicit) prime $p$. The unit
    \[ L_{K(h)} \mathbb S \to E_h \]
is a $K(h)$-local profinite Galois extension for the Goerss-Hopkins-Miller action of $\G_h$. That is, there are $K(h)$-local spectra $E_h^{hU}$ for every open subgroup $U$ of $\G_h$ such that the following hold:
\begin{enumerate}
    \item each $E_h^{hU}$ is an $\mathbb E_\infty$-ring spectrum over which $E_h$ is a commutative algebra,
    \item choosing a cofinal sequence of open subgroups $U$ yields $E_h \simeq L_{K(h)} \varinjlim_U E_h^{hU}$,
    \item for any normal inclusion $V \triangleleft U$ of open subgroups, the map $E_h^{hU} \to E_h^{hV}$ is a faithful $U/V$-Galois extension of $K(h)$-local spectra.
\end{enumerate}
\end{thm}

In order to leverage this to compute the groups $\Pic_h$ it is necessary to understand not just the finite Galois descent technology mentioned above, but how this assembles over the entire system of extensions. Our main theorem is the following descent result for Picard groups:

\begin{introthm} \label{intro:main theorem}
The unit $L_{K(h)} \mathbb S \to E_h$ induces an equivalence of spectra
\begin{equation} \label{intro:pic_n is homotopy fixed points}
    \pic(\mathcal Sp_{K(h)}) \simeq \tau_{\geq 0} \pic(E_h)^{h \G_h},
\end{equation}
where the right-hand side denotes \emph{continuous} homotopy fixed points. The resulting spectral sequence takes the form
\begin{equation} \label{intro:HFPSS}
    E_2^{s,t} = H^{s}_{cont} \left(\G_h, \pi_t \pic(E_h) \right) \implies \pi_{t-s} \pic(\mathcal Sp_{K(h)}).
\end{equation}
In an explicit range, it agrees with the $K(h)$-local $E_h$-Adams spectral sequence, including differentials.
\end{introthm}

Here $\pic(E_h)$ denotes the Picard spectrum of \emph{$K(h)$-local} $E_h$-modules. One of the major tasks is to properly interpret the right-hand side of (\ref{intro:pic_n is homotopy fixed points}), in order to take into account the profinite topology on the Morava stabiliser group; to do so, we make use of the \emph{pro\'etale} (or \emph{condensed/pyknotic}) formalism of \cite{bs,condensed,pyknotic}. We elaborate on our approach later in the introduction, but will first mention some consequences.

\subsubsection*{Picard group computations}
As a first corollary, we show how to recover the computation of $\Pic_1$.

\begin{introthm}[\cite{hopkins-mahowald-sadofsky}; \cref{lem:Pic_1 p odd,lem:Pic_1 n=2}] \label{intro:HMS}
The Picard group of the $K(1)$-local category is as follows:
\begin{enumerate}
    \item At odd primes, $\Pic_1 = \Z[p] \times \Z/{2(p-1)}$,
    \item When $p = 2$, $\Pic_1 = \Z[2] \times \Z/2 \times \Z/4$.
\end{enumerate}
\end{introthm}

This is obtained by computing the spectral sequence \cref{intro:HFPSS} from knowledge of the $K(1)$-local $E_1$-Adams spectral sequence. In fact, we focus on computing the exotic part of $\Pic_1$; in \cref{lem:algebraic Picard group} we show that this is precisely the part detected in \cref{intro:HFPSS} in filtration at least two.

We also consider examples at height bigger than two. When combining the results of \cite{hopkins-mahowald-sadofsky,kamiya-shimomura,karamanov,ghmr_pic,bbghps_pic,pstragowski_picard}, the following gives an algebraic expression for the first undetermined Picard group:

\begin{introthm}[\cite{culver-zhang}; \cref{lem:n^2 = 2p-1,ex:n^2 = 2p-1}] \label{intro:kappa 5}
If $p = 5$, there is an isomorphism
    \[ \kappa_3 \simeq H_0(\mathbb S_3, \pi_{8} E_3)^{\mathrm{Gal}(\F[125]/\F[5])}. \]
\end{introthm}

At height two, spectral sequence \cref{intro:HFPSS} implies the (known) result that exotic elements exist only at the primes $2$ and $3$; the groups $\kappa_2$ at these primes were computed in \cite{bbghps_pic} and \cite{ghmr_pic} respectively. Both cases used the filtation on $\Pic_2$ defined in \cite[Prop.~7.6]{hopkins-mahowald-sadofsky}, and explicit constructions of Picard elements. On the other hand, \Cref{intro:main theorem} gives another possible approach to compute $\kappa_2$, by computing \cref{intro:HFPSS} via comparison with the $K(2)$-local $E_2$-Adams spectral sequence; the latter is studied for example in \cite{henn-karamanov-mahowald,bbghps_ass}. We hope to return to this in future.

\subsubsection*{Brauer group computations}
In \cref{sec:brauer}, we turn our attention to Galois descent computations of the Brauer group of $\mathcal Sp_{K(h)}$. In \cite{gl}, Gepner and Lawson use Galois descent to compute the relative Brauer group $\Br(KO \mid KU)$ from knowledge of the HFPSS for $\pic(KO) \simeq \tau_{\geq 0} \pic(KU)^{hC_2}$, and we show that their method extends to our context. Note that while the Picard group is amenable to computation by a variety of methods, the Brauer group is much more difficult to study in an \emph{ad hoc} manner, and each of \cite{antieau-gepner,gl,antieau-meier-stojanoska} uses some such form of descent to obtain an upper bound; this follows the classical picture, in which descent computations were pioneered by Grothendieck in \cite{dix-exposees_ii}. The results of \cref{sec:brauer descent} strengthen existing descent technology, making the Brauer groups of the $K(h)$-local categories computationally tractable.

In analogy with the Picard case, we will write $\Br_h^0 \coloneqq \Br(\mathcal Sp_{K(h)} \mid E_h)$ for the group of Brauer classes of $\mathcal Sp_{K(h)}$ that become trivial (up to $K(h)$-local Morita equivalence) over Morava E-theory. In \cref{sec:Brauer n=1}, we consider the bound on this group imposed by spectral sequence \cref{intro:HFPSS}.

\begin{introthm}[\cref{lem:Brauer n = 1 p odd} and \cref{lem:Brauer n = 1 p = 2}] \label{intro:brauer}
\begin{enumerate}
    \item If $p$ is odd, then $\Br_1^0 \subset \Z/p-1$.
    \item At the prime $2$, we have $\lvert \Br_1^0 \rvert \leq 32$.
\end{enumerate}
\end{introthm}

We also give candidates for populating these groups; in \cite{relativebrauer}, we complete this computation. Combining this with the main theorem of \cite{hl} would give a description of the full Brauer group $\Br_h$, but we do not pursue this here.

\setcounter{introthm}{1}

\subsubsection*{Methods for profinite descent}
We now summarise our approach to continuity of the $\G_h$-action on $\pic(E_h)$. In the case of Morava E-theory itself, this was explored in the work of Davis \cite{davis_thesis, davis} and Quick \cite{quick_actions}, both of whom (using different methods) gave model-theoretic interpretations for continuous actions of profinite groups. In particular, both approaches recover the $K(h)$-local $E_h$-based Adams spectral sequence as a type of HFPSS for the action on $E_h$.

We will make crucial use of the \emph{pro\'etale classifying site} of the profinite group $\G_h$. This is closer in flavour to the approach of Davis, who uses the \emph{\'etale} classifying site; his strategy is briefly recounted in Section \ref{sec:Morava E-theory}, where we make the connection explicit. Namely, Davis makes use of the fact that viewed as an object of the $K(h)$-local category, Morava E-theory is a \textit{discrete} $\G_h$-object, meaning that it is the filtered colimit of the objects $E_h^{hU}$. Such $\G_h$-actions can be effectively modelled using the site of \emph{finite} $\G_h$-sets, and the requisite model category of simplicial sheaves was developed by Jardine \cite{jardine} as a way of formalising Thomason's work on descent for $K(1)$-localised algebraic K-theory. Of course, the action on $E_h$ is \emph{not} discrete when we view it as a plain spectrum, as can be seen on homotopy groups. Likewise, the induced action on its Picard spectrum is not discrete. As a model for more general continuous actions of a profinite group $G$, we therefore use the \inftycat{} of sheaves on the \emph{pro\'etale} site $\proet{G}$, whose objects are the profinite $G$-sets. This was studied in \cite{bs}; as shown there, in many cases it gives a site-theoretic interpretation of continuous group cohomology. Even when the comparison fails, pro\'etale sheaf cohomology exhibits many desirable properties absent in other definitions.

The equivalence in \cref{intro:main theorem} is therefore interpreted as the existence of a sheaf of connective spectra $\pic (\mathcal E_h)$ on the pro\'etale site, having
    \[ \Gamma(\G_h/*, \pic(\mathcal E_h)) \simeq \pic(\Mod_{E_h}(\mathcal Sp_{K(h)})) \qquad \text{and} \qquad \Gamma(\G_h/\G_h, \pic(\mathcal E_h)) \simeq \pic(\mathcal Sp_{K(h)}). \]

To this end, we begin by proving a descent result for Morava E-theory itself.

\begin{introsubthm}[\cref{lem:nu^p mathcal E is a sheaf} and \cref{lem:decalage for sheaves of spectra}]
There is a hypercomplete sheaf of spectra $\mathcal E_h$ on $\proet{(\G_h)}$ with
\begin{equation}
    \Gamma(\G_h/*, \mathcal E_h) \simeq E_h \qquad \text{and} \qquad \Gamma(\G_h/\G_h, \mathcal E_h) \simeq L_{K(h)} \mathbb S.    
\end{equation}
Its homotopy sheaves are given by $\pi_t \mathcal E_h \simeq \Cont_\G(-, \pi_t E_h)$, and its descent spectral sequence agrees with the $K(h)$-local $E_h$-Adams spectral sequence (including differentials).
\end{introsubthm}

 This may be of independent interest as it gives a novel construction of the $K(h)$-local $E_h$-Adams spectral sequence, which may be extended to an arbitrary spectrum $X$; its $E_2$-page (given \emph{a priori} in terms of sheaf cohomology on the pro\'etale site) is continuous group cohomology for suitable $X$ (\cref{rem:E2 page of ASS for arbitrary X}). Note moreover that pro\'etale cohomology enjoys excellent functoriality properties, and that the category of $L$-complete abelian sheaves on $\proet \G$ is abelian, as opposed to $L$-complete $(E_h)_*^\vee E_h$-comodules.

Next, we deduce a descent result for \inftycats{} of $K(h)$-local modules, which is really an extension to the condensed world of the following significant theorem:

\begin{thm}[\cite{mathew_galois}]
The diagram of symmetric monoidal \inftycats{}
    \[ \augcosimp[\mathcal Sp_{K(h)}]{\Mod_{E_h}(\mathcal Sp_{K(h)})}{\Mod_{L_{K(h)} E_h \wedge E_h}(\mathcal Sp_{K(h)})} \]
is a limit cone.
\end{thm}

Namely, in \cref{sec:descent theory} we prove the following profinite Galois descent result, which can be seen as the identification $\mathcal Sp_{K(h)} \simeq \left(\Mod_{E_h}(\mathcal Sp_{K(h)}) \right)^{h\G_h}$ analogous to \cref{intro:Galois descent}.

\begin{introsubthm}[\cref{lem:hyperdescent for presheaf extension of Modk mathcal E}]
There is a hypercomplete sheaf $\Mod_{\mathcal E_h}(\mathcal Sp_{K(h)})$ of symmetric monoidal \inftycats{} on $\proet{(\G_h)}$ with
    \[ \Gamma(\G_h/*, \Mod_{\mathcal E_h}(\mathcal Sp_{K(h)})) \simeq \Mod_{E_h}(\mathcal Sp_{K(h)}) \qquad \text{and} \qquad \Gamma(\G_h/\G_h, \Mod_{\mathcal E_h}(\mathcal Sp_{K(h)})) \simeq \mathcal Sp_{K(h)}. \]
\end{introsubthm}

One recovers the first part of \Cref{intro:main theorem} by taking Picard spectra pointwise. For the second part of that theorem, we must identify the $E_2$-page of the descent spectral sequence, which \emph{a priori} begins which sheaf cohomology on the pro\'etale site. The results of \cite{bs} allow us to deduce this, as a consquence of the following identification:

\begin{introsubthm}[\cref{lem:pi_* pic mathcal E}]
There is a hypercomplete sheaf of connective spectra $\pic(\mathcal E_h)$ on $\proet{(\G_h)}$ with
\begin{equation}
    \Gamma(\G_h/*, \pic(\mathcal E_h)) \simeq \pic(E_h) \qquad \text{and} \qquad \Gamma(\G_h/\G_h, \pic(\mathcal E_h)) \simeq \pic(\mathcal Sp_{K(h)}).
\end{equation}
The homotopy sheaves of $\pic(\mathcal E_h)$ are
    \[ \pi_t \pic(\mathcal E_h) = \operatorname{Cont}_{\G}(-, \pi_t \pic(E_h)), \]
i.e. represented by the homotopy groups of $\pic(E_h)$ (with their natural profinite topology).
\end{introsubthm}

The identification of homotopy groups is immediate for $t \geq 1$ (by comparing with $\pi_t \mathcal E_h$), but there is some work to do for $t = 0$; this is the same issue that accounts for the uncertainty in degree zero in the descent spectral sequence of \cite{heard}.

In fact we go a bit further, relating the spectral sequence of Theorem A to the $K(h)$-local $E_h$-Adams spectral sequence. In degrees $t \geq 2$, the homotopy groups of the Picard spectrum of an \Einfty-ring are related by a shift to those of the ring itself. It is a result of \cite{ms} that this identification lifts to one between truncations of the two spectra, in a range that grows with $t$: that is, for every $t \geq 2$ there is an equivalence
    \[ \tau_{[t, 2t-2]} \pic(A) \simeq \tau_{[t,2t-2]} \Sigma A, \]
functorial in the ring spectrum $A$. Using the pro\'etale model, it is quite straightforward to deduce the following comparison result, as proven in \opcit\ for finite Galois extensions.

\begin{introsubthm}[\cref{lem:additive differentials} and \cref{lem:nonlinear differential}; c.f. \cite{ms}] \label{intro:differentials}
\begin{enumerate}[leftmargin=\parindent, itemindent=*]
    \item Suppose $2 \leq r \leq t-1$. Under the identification
        \[ E_2^{s,t} = H^{s}(\G_h, \pi_t \pic(E_h)) \cong H^{s}(\G_h, \pi_{t-1} E_h) = E_2^{s,t-1}(ASS), \]
    the $d_r$-differential on the group $E_r^{s,t}$ in \cref{intro:HFPSS} agrees with the differential on $E_r^{s,t-1}(ASS)$ on classes that survive to $E_r$ in both.
    
    \item If $x \in H^t(\G_h, \pi_t \pic(E_h)) \cong H^t(\G_h, \pi_{t-1} E_h)$ (and $x$ survives to the $t$-th page in both spectral sequences), then the differential $d_t(x)$ is given by the following formula in the $K(h)$-local $E_h$-Adams spectral sequence:
        \[ d_t(x) = d^{ASS}_t(x) + x^2. \]
\end{enumerate}
\end{introsubthm}

\subsubsection*{Comparison with Morava modules}
Finally, we will says some words on how to derive \Cref{intro:kappa 5,intro:brauer} from the main result. Recall that a useful technique for computing Picard groups, originating already in \cite{hopkins-mahowald-sadofsky}, is to use completed E-theory to compare the category of $K(h)$-local spectra to the category $\Mod_{\pi_* E_h}^{\G_h}$ of \emph{Morava modules}, i.e. $L$-complete $\pi_* E_h$-modules equipped with a continuous action of the Morava stabiliser group $\G_h$:
\begin{equation} \label{eqn:Morava modules}
    (E_h)_*^{\vee}(-) \coloneqq \pi_* L_{K(h)}(E_h \wedge (-) ): \mathcal Sp_{K(h)} \to \Mod_{\pi_* E_h}^{\G_h}.
\end{equation}

This carries invertible $K(h)$-local spectra to invertible Morava modules, and hence induces a homomorphism on Picard groups. The category on the right-hand side is completely algebraic in nature, and its Picard group $\Pic_h^{\mathrm{alg}}$ can (at least in theory) be computed as a $\Z/2$-extension of $\Pic_h^{\mathrm{alg},0} = H^1(\G_h, (\pi_0 E_h )^\times)$; the strategy is therefore to understand the comparison map and how much of $\Pic_h$ it can see. A remarkable theorem of Pstrągowski says that the map $(E_h)_*^{\vee}: \Pic_h \to \Pic_h^{\mathrm{alg}}$ is an isomorphism if $p \gg n$ (more precisely, if $2p -2 > n^2 + n$). This reflects the more general phenomenon that chromatic homotopy theory at large primes is well-approximated by algebra, as is made precise in \cite{chromatic-asymptotically, monochromatic-asymptotically}.

As noted in \cite{pstragowski_picard}, the existence of a spectral sequence of the form (\ref{intro:HFPSS}) immediately yields an alternative proof, by sparseness of the $K(h)$-local $E_h$-Adams spectral sequence at large primes (in fact, this improves slightly the bound on $p$). Heard gave such a spectral sequence in \cite{heard}, and our results can be seen as a conceptual interpretation of that spectral sequence, analogous to the relation of \cite{davis, quick_actions} to \cite{devinatz-hopkins}. Beyond the conceptual attractiveness, our derivation of the spectral sequence also clarifies certain phenomena: for example, we give a proof of the claim made in \cite{heard} that the exotic part of the Picard group is given precisely by those elements in filtration at least 2:

\setcounter{introthm}{4}

\begin{introthm}[\cref{lem:algebraic Picard group}]
For any pair $(n,p)$, the algebraic Picard group is computed by the truncation to filtration $\leq 1$ of \cref{intro:HFPSS}, and the exotic Picard group $\kappa_n$ agrees with the subgroup of $\Pic_h$ in filtration at least 2 for \cref{intro:HFPSS}.
\end{introthm}

For example, when $n^2 = 2p-1$ this leads to the description of the exotic Picard group given in \Cref{intro:kappa 5}.

Note also that the group in bidegree $(s,t) = (1,0)$ of Heard's Picard spectral sequence is undetermined, which is an obstruction to computing Brauer groups; as discussed in \cref{sec:picard sheaf}, the relevant group in \cref{intro:HFPSS} really is $H^1(\G_h, \Pic(E_h))$. As expected, the computation simplifies at sufficiently large primes, and this should give rise to an algebraicity statement for the group $\Br_h^0$. We intend to explore the algebraic analogue $\Br_h^{0, alg}$ in future work.

\subsection{Outline}
In \cref{sec:Morava E-theory}, we collect the results we need on the Devinatz-Hopkins action and the pro\'etale site, showing how to define the sheaf of spectra $\mathcal E_h$. One can in fact deduce it is a sheaf from \cref{lem:hyperdescent for presheaf extension of Modk mathcal E}. Nevertheless, we wanted to give a self-contained proof of the spectrum-level hyperdescent; we also explain how this compares with Davis' construction of the continuous action on $E_h$. In the second half of \cref{sec:Morava E-theory} we compute the homotopy sheaves of $\mathcal E_h$, and explain how this leads to the identification with the $K(h)$-local $E_h$-Adams spectral sequence; the requisite d\'ecalage results are collected in \cref{app:decalage}.

In \cref{sec:descent theory} we categorify, obtaining descent results for the presheaf of $K(h)$-local module \inftycats{} over $\mathcal E_h$. We discuss how the Picard spectrum functor yields a sheaf of connective spectra exhibiting the identification of fixed points $\pic(\mathcal Sp_{K(h)}) \simeq \tau_{\geq 0} \pic(E_h)^{h\G_h}$, and investigate the resulting descent spectral sequence.

In \cref{sec:n=1} we use the previous results to compute Picard groups. We first identify the algebraic and exotic Picard groups in the descent spectral sequence. Combining this with the well-known form of the $K(1)$-local $E_1$-Adams spectral sequence allows us to reprove the results of \cite{hopkins-mahowald-sadofsky}: we are particularly interested in computing the exotic Picard group at the prime 2. We also consider Picard groups in the boundary case $n^2 = 2p-1$. In \cref{app:ASS n=1} we give a method to compute the height one Adams spectral sequence at $p=2$ using the Postnikov tower for the sheaf $\mathcal E_h$.

Finally, in \cref{sec:brauer} we show how to use our results in Brauer group computations. We show that the $(-1)$-stem in the Picard spectral sequence gives an upper bound for the relative Brauer group, and compute this bound at height one.

\subsection{Relation to other work}
As already mentioned, Heard has obtained a spectral sequence similar to that of \Cref{intro:main theorem}, and one of our objectives in this work was to understand how to view that spectral sequence as a HFPSS for the Goerss-Hopkins-Miller action.  Recently, a result close to \Cref{intro:main theorem} was proven by Guchuan Li and Ningchuan Zhang \cite{li-zhang}. Their approach differs somewhat from ours, using Burklund's result on multiplicative towers of generalised Moore spectra to produce pro-object presentations of $\Mod_{E_h}(\mathcal Sp_{K(h)})$ and $\pic_n$; a detailed comparison between the two would certainly be of interest.

\subsection{Acknowledgements}
I am grateful to Dustin Clausen for an outline of the argument of \cref{lem:hyperdescent for presheaf extension of Modk mathcal E}; in fact, part of my motivation for this project was his talk on Morava E-theory at the University of Copenhagen masterclass \emph{Condensed Mathematics}. I'd also like to thank Guchuan Li and Ningchuan Zhang for sharing with me their alternative approach. In addition, I have benefited greatly from conversations with Jack Davies, Ivo Dell'Ambrogio, David Gepner, Zachary Halladay, Luka Ilic, Shai Keidar, Thomas Nikolaus, Arthur Pander-Maat, Lucas Piessevaux, Maxime Ramzi, Vignesh Subramanian, Ivan Tomašić and Paul VanKoughnett, and would like to extend my thanks to them all. Above all, I am indebted to my supervisors Lennart Meier and Behrang Noohi, for their guidance, support and careful reading of previous drafts. Finally, I'd like to thank Utrecht University for their hospitality. This work forms a part of my thesis, supported by the Engineering and Physical Sciences Research Council [grant number EP/R513106/1].

\subsection{Notation and conventions}
\begin{itemize}
\item Throughout, we will work at a fixed prime $p$ and height $h$, mostly kept implicit. Also implicit is the choice of a height $h$ formal group law $\Gamma_h$ defined over $\mathbb F_{p^h}$; for concreteness we fix the Honda formal group, but this choice will not be important. For brevity, we will therefore write $\E, \K$ and $\G$ for Morava E-theory, Morava K-theory and the extended Morava stabiliser group, respectively. These will be our principal objects of study.

\item We will freely use the language of \inftycats{} (modeled as quasi-categories) as pioneered by Joyal and Lurie \cite{htt, ha, sag}. In particular, all (co)limits are $\infty$-categorical. We will mostly be working internally to the $K(h)$-local category, and as such we stress that \textbf{the symbol $\otimes$ will denote the $K(h)$-local smash product throughout}; where we have a need for it, we will use the notation $\wedge$ for the smash product of spectra (despite this having become archaic in some circles). On the other hand, we will distinguish $K(h)$-local colimits by writing for example $L_{\K} \colim X$ or $L_{\K} \varinjlim X_i$, as we feel that not to do so would be unnecessarily confusing. We use the notation $\varinjlim$ to denote a \emph{filtered} colimit, and similarly for cofiltered limits. In particular, if $T$ is a profinite set, we will use the expression `$T = \varprojlim T_i$' to refer to a presentation of $T$ as a pro-object, leaving implicit that each $T_i$ is finite. We will also assume throughout we have made a fixed choice of decreasing open subgroups $U_i \subset \G$ with trivial intersection; the symbols $\varprojlim_i$ and $\varinjlim_i$ will always refer to the (co)limit over such a family. Likewise, we will assume we have chosen a sequence of ideals $I \subset \pi_0 E_h$ generating the $\mathfrak m$-adic topology; without loss of generality, the ideals $I$ will be chosen so that there exists a tower $M_I$ of generalised Moore spectra, with $\pi_* M_I = (\pi_* E_h)/I$.

\item We only consider spectra with group actions, and not any more sophisticated equivariant notion. When $G$ is a profinite group, we will write $H^*(G,M)$ for \emph{continuous} group cohomology with pro-$p$ (or more generally profinite) coefficients, as defined for example in \cite{symonds-weigel} (resp. \cite{jannsen}).

\item A few words about spectral sequences. When talking about the `Adams spectral sequence', we always have in mind the $K(h)$-local $E_h$-based Adams spectral sequence; the classical Adams spectral sequence (based on $H\mathbb F_p$) makes no appearance in this document. We will freely use abbreviations such as `ASS', `HFPSS', `BKSS'. We will also use the name `descent spectral sequence' for either the t-structure or \v Cech complex definition, since in the cases of interest we show they agree up to reindexing; when we need to be more explicit, we refer to the latter as the `Bousfield-Kan' or `\v Cech' spectral sequence. The name `Picard spectral sequence' will refer to the descent spectral sequence for the sheaf $\pic(\mathcal E_h)$.

\item In \cref{sec:Morava E-theory,sec:descent theory} we form spectral sequences using the usual t-structure on spectral sheaves; this is useful for interpreting differentials and filtrations, for example in \cref{lem:algebraic Picard group}. To obtain familiar charts, we will declare that the spectral sequence associated to a filtered object starts at the $E_2$-page; in other words, this is the page given by homotopy groups of the associated graded object. Thus our spectral sequences run
    \[ E_2^{s,t} = H^s(G, \pi_t E) \implies \pi_{t-s} E^{hG}, \]
with differentials $d_r$ of $(s,t-s)$-bidegree $(r,-1)$, and this is what we display in all figures. However, we also make use of the Bousfield-Kan definition of the descent spectral sequence using the \v Cech complex of a covering, and we relate the two formulations by d\'ecalage (see \cref{app:decalage}): there is an isomorphism between the two spectral sequences that reads
    \[ E_2^{s,t} \cong \check E_3^{2s-t, s} \]
if we use the same grading conventions for each of the underlying towers of spectra. We will always use $s$ for filtration, $t$ for internal degree, and $t-s$ for stem.

Our Lyndon-Hochschild-Serre spectral sequences are in cohomological Serre grading.

\item Finally, we largely ignore issues of set-theory, since these are discussed at length in \cite{condensed} and \cite{pyknotic}. For concreteness, one could fix a hierarchy of strongly inaccessible cardinals $\kappa < \delta_0 < \delta_1$ such that $|\G_h| < \kappa$ and the unit in $\mathcal Sp_{K(h)}$ is $\kappa$-compact, and work throughout over the `$\delta_1$-topos' of sheaves of $\delta_1$-small spaces on profinite $\G_h$-sets of cardinality less than $\delta_0$.
\end{itemize}

\printunsrtglossary[
    title = {\normalsize List of symbols}
]

%% file: Sections/Morava_E-theory.tex
\subsection{Discrete Morava E-theory} \label{sec:discrete E-theory}
Let $\E \coloneqq E(\F[p^h], \Gamma_h)$ be Morava $E$-theory based on the Honda formal group at height $h$ and prime $p$; $h$ and $p$ will henceforth be fixed, and kept implicit to ease notation. Let $\K \coloneqq K(\F[p^h], \Gamma_h)$ be Morava K-theory, its residue field. Recall that $\E$ is the $\K$-local Landweber exact spectrum whose formal group is the universal deformation of $\Gamma_h$ to the Lubin-Tate ring $\pi_* \E = \mathbb W(\F[p^h])[[u_1, \dots, u_{n-1}]][u_h^{\pm 1}]$. Functoriality yields an action of the extended Morava stabiliser group $\G = \G_h \coloneqq \operatorname{Aut}(\F[p^h], \Gamma_h)$ on the homotopy ring spectrum $\E$, and celebrated work of Goerss, Hopkins, Miller and Lurie \cite{goerss-hopkins,ellipticii} promotes $\E$ to an \Einfty-ring and the action to one by \Einfty\ maps. This action controls much of the structure of the $\K$-local category, and is the central object of study in this document. In this section, we formulate the action of $\G$ on $\E$ in a sufficiently robust way for our applications; to do so, we will present $\E$ as a sheaf of spectra on the \emph{pro\'etale classifying site} of $\G$. Descriptions of the $\K$-local $\E$-Adams spectral sequence have been previously given, notably in work of Davis and of Quick \cite{davis_thesis, davis, quick_actions, behrens-davis, davis-quick}, who described a number of formulations of this action as the \emph{continuous} action of the \emph{profinite} group $\G$.

Recall that continuous actions and continuous cohomology of a topological group $G$ are generally much more straightforward when we assume our modules to have \emph{discrete} topology. There are notable categorical benefits in this case: for example, it is classical that the category of discrete $G$-modules is abelian with enough injectives, which is not true of the full category of topological modules. Further, in the discrete context we can understand actions completely by looking at the induced actions of all finite quotients of $G$. This was pioneered by Thomason in his study of $K(1)$-local descent for algebraic $K$-theory, and formalised in a model-theoretic sense by Jardine. 

Any profinite group $G$ has an \emph{\'etale classifying site}, denoted $\et G$, whose objects are the (discrete) finite $G$-sets and whose coverings are surjections; as shown in \cite[\S6]{jardine}, the category of sheaves of abelian groups on $\et G$ gives a category equivalent to the category $\Ab_G^\delta$ of discrete $G$-modules in the sense of \cite{serre_galois}. As noted below, sheaf cohomology on $\et G$ corresponds to continuous group cohomology with discrete coefficients, again in the sense of Serre. Motivated by this, Jardine defines a model structure on presheaves of \emph{spectra} on $\et G$, which models the category of `discrete' continuous $G$-spectra, i.e. those that can be obtained as the filtered colimit of their fixed points at open subgroups. We will more generally refer to objects of $\Shh[\mathcal C]{\et G}$ 
as \emph{discrete $G$-objects} of an arbitrary (cocomplete) \inftycat{} $\mathcal C$. Davis uses this as his starting point, and we observe below that in this formalism it is easy to pass to the $\infty$-categorical setting. Namely, we begin by collecting the following facts:

\begin{thm}[Devinatz-Hopkins, Davis, Rognes, Dugger-Hollander-Isaksen] \label{lem:definition of mathcal E}
    There is a hypercomplete sheaf of spectra $\mathcal E^\delta$ on $\et \G$, such that
    \begin{enumerate}
        \item Any ordered, cofinal sequence $(U_i)$ of open subgroups of $\G$ induces $L_{\K} \varinjlim_i \mathcal E^\delta (G/U_i) \simeq \E$,
        \item On global sections, $\Gamma \mathcal E^\delta \coloneqq \Gamma(\G/\G, \mathcal E^\delta) \simeq \bm 1_{\K}$ is the $\K$-local sphere spectrum,
        \item $\mathcal E^\delta$ lifts to $\operatorname{CAlg}(\mathcal Sp_{\K})$,
        \item For any normal inclusion of open subgroups $V \subset U \subset \G$ , the map $\mathcal E^\delta(\G/U) \to \mathcal E^\delta(\G/V)$ is a faithful $U/V$-Galois extension.
    \end{enumerate}
\end{thm}

\begin{proof}
The presheaf of spectra $\mathcal E^\delta$ is constructed in \cite[\S4]{devinatz-hopkins}, with $(ii)$ being part of Theorem 1 therein and the identification $(i)$ being the trivial case of Theorem 3; see also \cite[\S2]{barthel-beaudry-goerss-stojanoska} for a nice summary. Devinatz and Hopkins construct $\mathcal E^\delta$ by hand (by taking the limit of the \emph{a priori} form of its Amitsur resolution in $K(h)$-local $E_h$-modules); they denote $\Gamma(\G/U, \mathcal E^\delta)$ by $E^{hU}_h$, but we copy \cite{behrens-davis} and write $E^{dhU}_h$. By \cite[Theorem~4]{devinatz-hopkins}, $\mathcal E^\delta$ is a sheaf on $\et \G$; in fact, Devinatz and Hopkins already define $\mathcal E^\delta$ to land in $\K$-local $\mathbb E_\infty$-rings, and so we may consider it as an object of $\Sh[\operatorname{CAlg}(\mathcal{S}p_{\K})]{\et \G}$ since limits in $\operatorname{CAlg}(\mathcal Sp_{\K})$ are computed at the level of underlying spectra \cite[Corollary~3.2.2.5]{ha}. Item $(iv)$ is \cite[Theorem~5.4.4]{rognes}, so it remains to show that $\mathcal E$ is hypercomplete; this we will deduce from Davis' work.

More specifically, Davis utilises the Jardine model structure on the category of presheaves of spectra on $\et \G$, denoted $\operatorname{Spt}_{\G}$\footnote{We have taken slight notational liberties: in \cite{davis_thesis}, $\operatorname{Spt}_\G$ denotes the category of spectra based on discrete $\G$-sets, which is equipped with a model structure lifted from the Jardine model structure on the (equivalent) category $\mathbf{ShvSpt}$ of sheaves of spectra.}; this is defined in such a way that there is a Quillen adjunction
    \[ \operatorname{Spt} \xrightleftarrows{\text{const}}{(-)^{h\G}} \operatorname{Spt}_{\G}. \]
Recall that the main result of \cite{dugger-hollander-isaksen} says that the fibrant objects of $\operatorname{Spt}_{\G}$ are precisely those projectively fibrant presheaves that satisfy $(i)$ the (1-categorical) sheaf condition for coverings in $\et \G$, and $(ii)$ descent for all hypercovers, and so the \inftycat{} associated to $\operatorname{Spt}_{\G}$ is a full subcategory of $\Shh[\mathcal Sp]{\et \G}$.

In this setting, Davis shows that the spectrum $F_h \coloneqq \varinjlim E^{dhU}_h$ (the colimit taken in plain spectra) defines a fibrant object of $\operatorname{Spt}_{\G}$, $\mathcal F \colon \G/U \mapsto F_h^{hU}$ \cite[Corollary~3.14]{davis}; Behrens and Davis show that $E_h^{hU} \simeq L_{\K} F_h^{hU} \simeq E_h^{dhU}$ for any open subgroup $U \subset \G$ \cite[Discussion above Prop.~8.1.2 and Lemma~6.3.6 respectively]{behrens-davis}. In fact, \cite[Theorem~8.2.1]{behrens-davis} proves the same equivalence for any \emph{closed} subgroup $H$, but restricting our attention to open subgroups cuts out some of the complexity and makes clear the equivalences happen naturally in $U$ (we remark that they also appeared in Chapter 7 of Davis' thesis \cite{davis_thesis}). This provides a projective equivalence between $\mathcal E^\delta$ and $L_{\K} \mathcal F$, and hence an equivalence of \emph{presheaves} between $\mathcal E^\delta$ and a hypercomplete sheaf of spectra.
\end{proof}

If we pick a cofinal sequence of open normal subgroups $(U_i)$ we can identify the starting page of the descent spectral sequence for $\mathcal E^\delta$:

\begin{lem} \label{lem:constructing the descent spectral sequence}
    Let $G$ be a profinite group and $\mathcal{F} \in \Shh[\mathcal C]{\et G}$, where $\mathcal C = \mathcal Sp$ or $\mathcal Sp_{\geq 0}$. There is a spectral sequence with starting page
        \[E_2^{s,t} = H^s(G, \pi_t \varinjlim \mathcal{F}(G/U_i)),\]
    and converging conditionally to $\pi_{t-s} \Gamma \lim_t \tau_{\leq t} \mathcal{F}$.
\end{lem}

\begin{proof}
This is the spectral sequence for the Postnikov tower of $\mathcal{F}$, formed as in \cite[\S1.2.2]{ha}. Its starting page is given by sheaf cohomology of the graded abelian sheaf $\pi_* \mathcal F$ on $\et G$. Explicitly, form the Postnikov tower in sheaves of spectra
\[\begin{tikzcd}[column sep=small]
	{\mathcal F} & \cdots & {\tau _{\leq 1}\mathcal F} & {\tau _{\leq 0}\mathcal F} &  \cdots \\
	&& {\Sigma \pi_1 \mathcal F} & {\pi_0 \mathcal F} &
	\arrow[from=1-1, to=1-2]
	\arrow[from=1-2, to=1-3]
	\arrow[from=1-3, to=1-4]
	\arrow[from=1-4, to=2-4]
	\arrow[from=1-3, to=2-3]
	\arrow[from=2-4, to=1-3]
	\arrow[from=2-3, to=1-2]
	\arrow[from=1-4, to=1-5]
\end{tikzcd}\]
Applying global sections and taking homotopy groups gives an exact couple, and we obtain a spectral sequence with
    \[ E_2^{s,t} = \pi_{t-s} \Gamma \Sigma^t \pi_t \mathcal F = R^{s} \Gamma \pi_t \mathcal F = H^{s} (\et \G, \pi_t \mathcal F), \]
with abutment $\pi_{t-s} \Gamma \lim \tau_{\leq t} \mathcal F$. To identify this with continuous group cohomology, we make use of the equivalence
\begin{align*}
    \Ab^\delta_G & \to \Sh[\Ab]{\et G} \\
    M & \mapsto \left(\coprod G/U \mapsto \prod M^U\right)
\end{align*}
of \cite{jardine}, whose inverse sends $\mathcal F \mapsto \varinjlim_i \mathcal F(G/U_i)$. Under this equivalence the fixed points functor on $\Ab^\delta_G$ corresponds to global sections, and so taking derived functors identifies sheaf cohomology on the right-hand side with derived fixed points on the left; as in \cite[\S2.2]{serre_galois}, when we take discrete coefficients this agrees with the definition in terms of continuous cochains.
\end{proof}

Writing $\Modk{(-)} \coloneqq \operatorname{Mod}_{L_{\K} (-)}(\mathcal Sp_{\K}) = L_{\K} \operatorname{Mod}_{(-)}$, our strategy is is to apply the functor $\pic \circ \Modk{(-)}: \operatorname{CAlg} \to \mathcal Sp_{\geq 0}$ pointwise to the sheaf $\mathcal E^\delta$, in order to try to obtain a sheaf $\pic_{\K}(\mathcal E^\delta) \in \Shh[\mathcal{S}p_{\geq 0}]{\et \G}$; we'd then like to apply the above lemma to deduce the existence of the descent spectral sequence. This does not quite work for the same reason that the lemma applied to $\mathcal E^\delta$ does not recover the $\K$-local $\E$-Adams spectral sequence: while $\E$ is $\K$-locally discrete, it is certainly not discrete as a $\G$-spectrum (for example, even the action on $\pi_2 \E$ is not discrete). Nevertheless, it is worth remarking that the first step of this approach does work: since $\E$ is a discrete $\G$-object of \emph{$\K$-local spectra}, $\Modk{\E}$ is discrete as a presentable \inftycat{} with $\G$-action. This is a consequence of the following two results.


\begin{lem}
    The composition $\Modk{\mathcal E^\delta}: {\et \G}^{op} \to \operatorname{CAlg}(\mathcal{S}p_{\K}) \to \Prsm$ is a sheaf.
\end{lem}

\begin{proof}
To check the sheaf condition for $\mathcal{F}: \et \G^{op} \to \mathcal{C}$ we need to show that finite coproducts are sent to coproducts, and that for any inclusion $U \subset U'$ of open subgroups the canonical map $\mathcal{F}(\G/U') \to \Tot \mathcal{F}(\G/U^{\times_{\G/U'}\bullet + 1})$ is an equivalence. For the presheaf $\Modk{\mathcal E^\delta}$, the first is obvious (using the usual idempotent splitting), while the second is \textit{finite} Galois descent \cite[Proposition~6.2.6]{meier_thesis} or \cite[Theorem~6.10]{gl}, at least after refining $U$ to a normal open subgroup of $U'$.
\end{proof}

Fixing again a cofinal sequence $(U_i)$, we write $F_{ij} \dashv R_{ij}: \Modk{\mathcal E^\delta(\G/U_i)} \rightleftarrows \Modk{\mathcal E^\delta (\G/U_j)}$ and $F_j \dashv R_j$ for the composite adjunction $\Modk{\mathcal E^\delta(\G/U_i)} \rightleftarrows \Modk{\E}$, and
    \[ \varinjlim \Modk{\mathcal E^\delta(\G/U_i)} \stackrel[R]{F}{\rightleftarrows} \Modk \E \]
for the colimit (along the functors $F_{ij}$) in $\Prsm$.

\begin{prop} \label{lem: Mod E is colimit of Mod mathcal E(S_i)}
    The functors $F$ and $R$ define an adjoint equivalence $\varinjlim \Modk{\mathcal E^\delta(\G/U_i)} \simeq \Modk \E$.
\end{prop}

More generally:

\begin{prop} \label{lem: colim of modules is modules over colim}
    Let $\mathcal C$ be a presentably symmetric monoidal stable $\infty$-category. Suppose $A_{(-)}: I \to \operatorname{CAlg}(\mathcal C)$ is a filtered diagram, and write $A$ for a colimit (formed equivalently in $\mathcal C$ or $\operatorname{CAlg}(\mathcal C)$). Then the induced adjunction
        \[\varinjlim \operatorname{Mod}_{A_i}(\mathcal{C}) \stackrel[R]{F}{\rightleftarrows} \operatorname{Mod}_A(\mathcal C) \]
    is an equivalence of \emph{presentable} symmetric-monoidal $\infty$-categories.
\end{prop}


\begin{proof}
This is implied by \cite[Corollary~4.8.5.13]{ha}, since filtered categories are weakly contractible and filtered colimits in $\CAlg_{\mathcal C}(\Pr^L)$ computed in $\Pr^L$. Alternatively, one can give an explicit description of the unit and counit and verify that both are natural equivalences.
\end{proof}

\begin{warning} \label{warning:discreteness}
The result above \textit{fails} if we consider the same diagram in $\operatorname{Cat}_\infty^{\mathrm{smon}}$ (after restricting to $\kappa$-compact objects for $\kappa$ a regular cardinal chosen such that $\bm 1_{\K} \in \mathcal Sp_{\K}$ is $\kappa$-compact, say). Indeed, one can see that the homotopy type of the mapping spectra $\operatorname{map}(\bm 1, X)$ out of the unit in this colimit would be different from that in $\Modk{\E}$, an artefact of the failure of $\K$-localisation to be smashing.

Applying $\mathfrak{pic}: \Prsm \to \mathcal Sp_{\geq 0}$ to the coefficients, one obtains a sheaf of connective spectra on $\et \G$ given by
    \[G/U_i \mapsto \mathfrak{pic}(\Modk{\mathcal E^\delta(G/U_i)}). \]
In particular, $\Gamma \mathfrak{pic}( \Modk{\mathcal E^\delta}) = \pic_h \coloneqq \mathfrak{pic}(\mathcal Sp_{\K})$. Unfortunately, this sheaf is unsuitable for the spectral sequence we would like to construct: the $E_1$-page will be group cohomology with coefficients in the homotopy of $\varinjlim \mathfrak{pic}(\Modk{\mathcal E^\delta(G/U_i)})$, which by \cite[Proposition~2.2.3]{ms} is the Picard spectrum of the colimit of module categories computed in $\operatorname{Cat}^{\mathrm{smon}}_\infty$; as noted above, this need \emph{not} agree with $\pic(\Modk{\E})$.
\end{warning}

\subsection{Pro\'etale homotopy theory} \label{sec:condensed}
In \cref{sec:discrete E-theory}, we showed how the work of Devinatz-Hopkins and Davis defines Morava E-theory as a discrete $\G$-object in the category of $\K$-local \Einfty-rings, and that using this it is straightforward to construct a discrete $\G$-action on the presentable \inftycat{} $\Modk{\E}$, having $\mathcal Sp_{\K}$ as fixed points. As noted in Warning \ref{warning:discreteness}, this does not suffice to prove our desired descent result for Picard spectra. Likewise, the descent spectral sequence for the sheaf $\mathcal E^\delta$ on $\et \G$ is \emph{not} the $\K$-local $\E$-Adams spectral sequence: the issue is that the action on the (unlocalised) spectrum $\E$ is not discrete. We are led to the following solution: rather than work with the site of \emph{finite} $\G$-sets, we think of Morava E-theory as a sheaf $\mathcal E$ on \emph{profinite} $\G$-sets.

\begin{defn}
    Let $G$ be a $\delta_0$-small profinite group. We denote by $\proet G$ the \emph{pro\'etale site}, whose underlying category consists of all $\delta_0$-small continuous profinite $G$-sets and continuous equivariant maps. $\proet G$ is equipped with the topology whose coverings are collections $\{S_\alpha \to S\}$ for which there is a \emph{finite} subset $A$ with $\coprod_{\alpha \in A} S_\alpha \twoheadrightarrow A$. 
\end{defn}

The pro\'etale site and resulting $1$-topos were extensively studied by Bhatt and Scholze in \cite{bs}; when $G$ is trivial, one recovers the \emph{condensed}/\emph{pyknotic} formalism of \cite{condensed, pyknotic}, up to a choice of set-theoretic foundations. One key feature for our purposes is that sheaf cohomology on $\proet G$ recovers continuous group cohomology for a wide range of coefficient modules; this will allow us to recover the desired homotopy sheaves, which we saw we could not do using $\et \G$---in particular, we will recover the $\K$-local $\E$-Adams spectral sequence as a descent spectral sequence. Moreover, $\proet \G$ gives a site-theoretic definition of continuous group cohomology, which makes the passage to actions on module categories (and hence Picard spectra) transparent.

 To begin, we will discuss some generalities of pro\'etale homotopy theory. Most important for the rest of the document will be \cref{sec:extending discrete objects}, in which we give criteria for a sheaf on $\et G$ to extend without sheafification to one on $\proet G$; the other sections will be used only cursorily, but are included for completeness.

\subsubsection{Free \texorpdfstring{$G$}{G}-sets} \label{sec:Free G}
A significant difference between the \'etale and pro\'etale sites is that the latter contains $G$-torsors. In fact, for any profinite set $T$ we have an object $G \times T \in \proet G$. It will often be useful to restrict to such free objects.

\begin{lem}
If $S$ is a continous $G$-set, the projection $S \to S/G$ is split: that is, there is a continuous section
    \[ S/G \to S. \]
In particular, if the $G$-action on $S$ is free then $S \cong S/G \times G$.
\end{lem}

\begin{proof}
This is proven as part of \cite[Prop.~3.7]{scholze}.
\end{proof}

\begin{defn}
The subsite $\Free G \subset \proet G$ is given by the full subcategory of $G$-sets with \emph{free} $G$-action; equivalently, these are $G$-sets isomorphic to ones of the form $T \times G$ for $T$ a profinite set with trivial $G$-action.
\end{defn}

\begin{lem} \label{lem:Free generates proet}
The subsite $\Free G$ generates $\proet G$: any $T \in \proet G$ admits a covering by a free $G$-set. Consequently, restriction induces an equivalence
    \[ \mathcal Sh(\proet G) \xrightarrow{\sim} \mathcal Sh(\Free G). \]
\end{lem}

\begin{proof}
The map $S \times G \to S$ is surjective, and the (diagonal) $G$-action on the domain is clearly free.
\end{proof}

\begin{remark}
In fact, one can restrict further to the category $\Proj G = \proet{G}^{\mathrm{wc}}$ of \emph{weakly contractible} $G$-sets: these are the free $G$-sets of the form $G \times T$, where $T$ is extremally disconnected. While this is not a site (it does not have pullbacks), $\Proj G$ does generate the \emph{hypercomplete} pro\'etale topos, in the sense that restriction induces an equivalence
    \[ \widehat{\mathcal Sh}(\proet G) \xrightarrow{\sim} \mathcal P_{\Sigma}(\Proj G), \]
where the codomain is the full subcategory of $\mathcal P(\Proj G)$ spanned by multiplicative presheaves---i.e., those that send binary coproducts to products.
\end{remark}

\subsubsection{Descent spectral sequence}
Suppose that $G$ is a profinite group, and $\mathcal C$ a prestable \inftycat{}. Then $\Sh[\mathcal Sp(\mathcal C)]{\proet G}$ obtains a t-structure, whose coconnective part is 
    \[ \Sh[\mathcal Sp(\mathcal C)]{\proet G}_{\leq 0} \coloneqq \Sh[\mathcal Sp(\mathcal C)_{\leq 0}]{\proet G}. \]
In particular, if $\mathcal F \in \Sh[\mathcal Sp(\mathcal C)]{\proet G}$ then $\mathcal F$ has a Postnikov tower
    \[ \mathcal F \to \cdots \to \mathcal F_{\leq n} \to \cdots \]
whose sections $\pi_t \mathcal F$ live in $\Sh[\mathcal C^\heartsuit]{\proet G}$. If $\mathcal F$ is \emph{hypercomplete}, the Postnikov tower converges by \cite{mondal-reinecke}. If $H$ is any subgroup, we may evaluate the Postnikov tower at the $G$-set $G/H$, obtaining a tower in $\mathcal Sp(\mathcal C)$ and hence a spectral sequence
    \[ E_2^{s,t} = \pi_{t-s} \Gamma(G/H, \pi_t \mathcal F) \implies \pi_{t-s} \Gamma(G/H, \mathcal F) \]
valued in the abelian category $\mathcal C^\heartsuit$.

\begin{cor} \label{lem:proetale descent spectral sequence}
If $\mathcal F \in \Shh[\mathcal Sp]{\proet G}$ or $\Shh[\mathcal Sp_{\geq 0}]{\proet G}$, there is a spectral sequence
    \[ E_2^{s,t} = H^s(\proet G, \pi_t \mathcal F) \implies \pi_{t-s} \Gamma \mathcal F. \]
\end{cor}

When $\mathcal A = \Cont_G(-, M) \in \Sh[\Ab]{\proet G}$ is represented by some topological $G$-module, there is a comparison map
    \[ \Phi: H^s_{\mathrm{cont}}(G, M) \to H^s(\proet G, A) \]
from continuous group cohomology, which arises as the edge map in a \v Cech-to-derived functor spectral sequence. In this context, \cite[Lemma~4.3.9]{bs} gives conditions for $\Phi$ to be an isomorphism.

\subsubsection{Extending discrete objects} \label{sec:extending discrete objects}
For a profinite group $G$, the pro\'etale site $\proet G$ is related to the site $\et G$ of finite $G$-sets by a map of sites $\nu: \et G \hookrightarrow \proet G$, which induces a geometric morphism at the level of topoi. More generally, if $\mathcal C$ is any complete and cocomplete $\infty$-category we write
    \[\nu^*: \Sh[\mathcal C]{\et G} \rightleftarrows \Sh[\mathcal C]{\proet G}: \nu_*\]
for the resulting adjunction; the left adjoint is the sheafification of the presheaf extension, given in turn by
    \[ \nu^p \mathcal F: S = \varprojlim_i S_i \mapsto \varinjlim_i \mathcal F(S_i) \numberthis \label{eqn:expression for nu^p for proetale site}\]
when each $S_i$ is a finite $G$-set. When $\mathcal C = \operatorname{Set}$, it is a basic result of \cite{bs} that $\nu^* = \nu^p$: essentially, this is because \begin{enumerate}
    \item the sheaf condition is a finite limit,
    \item it therefore commutes with the filtered colimit in (\ref{eqn:expression for nu^p for proetale site}).
\end{enumerate} 
This fails for sheaves valued in an arbitrary $\infty$-category $\mathcal C$ with limits and filtered colimits, where both properties might fail; as we now show, one can nevertheless make certain statements about sheaves of spaces or spectra, for example under cohomological dimension assumptions on $G$. 


\begin{defn}
Let $G$ be a profinite group. The \emph{mod $p$ cohomological dimension} of $G$ is
    \[ \cd_p(G) \coloneqq \sup \left\{ d: \text{there exists a $p$-power torsion $G$-module $M$ with } H^d(G, M) \neq 0 \right\} \in [0, \infty]. \]
The \emph{mod $p$ virtual cohomological dimension} $\vcd_p(G)$ is the smallest mod $p$ cohomological dimension of an open subgroup $U \subset G$.
\end{defn}

\begin{defn}
\begin{enumerate}[itemindent = *, leftmargin = \parindent]
    \item Let $Y: \Z[\geq 0]\op \to \mathcal C$ be a pretower in a stable \inftycat{} $\mathcal C$. If $\mathcal C$ has sequential limits, we can form the map
        \[ f: \{\lim Y\} \to \{Y_n\} \]
    in the \inftycat{} $\operatorname{Fun} (\Z[\geq 0]\op, \mathcal C)$ (in which the source is a constant tower). Recall that $Y$ is \emph{$d$-rapidly convergent to its limit} if each $d$-fold composite in $\fib(f) \in \operatorname{Fun} (\Z[\geq 0]\op, \mathcal C)$ is null (c.f. \cite[Definition~4.8]{cm}; see also \cite{hopkins-palmieri-smith}). More generally, if $Y: (\Z[\geq 0] \cup \{\infty\})\op \to \mathcal C$ is a tower, we say $Y$ is \emph{$d$-rapidly convergent} if the same condition holds for the map
        \[ f: \{Y_{\infty}\} \to \{Y_n\}. \]
    In particular, this implies that $Y_{\infty} \simeq \lim Y$; in fact, $F(Y_{\infty}) \simeq \lim F(Y)$ for any exact functor $F: \mathcal C \to \mathcal D$.
    
    \item Now suppose $\mathcal C$ also has finite limits, and $X^\bullet: \Delta_+ \to \mathcal C$ is an augmented cosimplicial object. Recall that $X$ is \emph{$d$-rapidly convergent} if the $\Tot$-tower $\{\Tot_{n} X^\bullet\}$ is $d$-rapidly convergent. In particular $X_{\infty} \simeq \Tot X^\bullet$, and the same is true after applying any exact functor $F$.
\end{enumerate}
\end{defn}

\begin{ex} \label{lem: A otimes cobar complex is 1-rapidly convergent}
Given an $\mathbb E_1$-algebra $A$ in any presentably symmetric monoidal stable \inftycat{} $\mathcal C$, the Adams tower $T(A,M)$ for a module $M$ over $A$ is defined by the property that $A \otimes T(A,M)$ is 1-rapidly convergent. It is well-known to agree with the $\Tot$-tower for the Amitsur complex for $M$ over $A$ (for example, this is worked out in detail in \cite[\S2.1]{mathew-naumann-noel}), and in particular the cosimplicial object $A^{\otimes \bullet + 2} \otimes M$, given by smashing the Amitsur complex with a further copy of $A$, is always $1$-rapidly convergent.
\end{ex}

\begin{remark}
Let $G$ be a profinite group and $\mathcal F$ a presheaf on $\proet G$, valued in a stable \inftycat{}. We will make use of the following assumption:
\begin{itemize}
    \item[$(\star)$] There exists $d \geq 0$ such that for any inclusion $V \subset U$ of normal subgroups, the \v Cech complex
        \[ \augcosimp[\mathcal F(G/U)]{\mathcal F(G/V)}{\mathcal F(G/V \times_{G/U} G/V)} \]
    is $d$-rapidly convergent.
\end{itemize}
The key observation is that if $p \colon I \to \Fun(\Delta^1, \proet{G})$ is a diagram of coverings and $\mathcal F$ satisfies $(\star)$, then $\mathcal F$ also satisfies descent for the covering $\lim(p)$; more generally, if $q \colon I \to \Fun(\Delta, \mathcal C)$ is a diagram of $d$-rapidly convergent cosimplicial objects in any stable $\mathcal C$, then $\varinjlim (q) \colon \Delta \to \mathcal C$ is also $d$-rapidly convergent.
\end{remark}

When $G$ has finite cohomological dimension mod $p$, one has the following corollary of \cite[Theorem~4.16]{cm} (which itself goes back to Thomason's work on $L_{K(1)} K(-)$).

\begin{lem} \label{lem:proetale descent for hypercomplete etale sheaves}
    Suppose $G$ is a profinite group with $\cd_p(G) = d$, and $k$ any $p$-local spectrum. Suppose $\mathcal F \in \Sh[\mathcal Sp_k]{\et G}$ is a \emph{hypercomplete} sheaf on $\et G$. Then $\nu^p \mathcal F$ is a sheaf of $k$-local spectra on $\proet G$.
\end{lem}

\begin{proof}
As above, $\nu^p \mathcal F$ is given by the formula
    \[ S = \varprojlim S_i \mapsto L_k \varinjlim \mathcal F(S_i). \]
Suppose we are given a surjection $S' = \lim_{j' \in J'} S'_{j'} \xrightarrow{\alpha} S = \lim_{j \in J} S_j$ of profinite $G$-sets. By \cite[Lemma~3.6]{scholze}, one can present this as the limit of finite covers $S'_i \to S_i$. We therefore want to show that the following is a limit diagram:
\begin{align} \label{eqn:Amistur complex for sheaf on et G}
    \augcosimp[L_k \varinjlim \mathcal F(S_i)]{L_k \varinjlim \mathcal F(S_i')}{L_k \varinjlim \mathcal F(S_i' \times_{S_i} S_i')}.
\end{align}

Localisation at $k$ commutes with totalisations, since these preserve the category of $k$-acyclic spectra: indeed, finite limits of acyclics are clearly acyclic, while $\Z$-shaped limits of acyclics are acyclic by virtue of the Milnor sequence. It therefore suffices to check (\ref{eqn:Amistur complex for sheaf on et G}) is a limit diagam before localising at $k$, instead working in $\mathcal Sp_{(p)}$. 

We claim that each of the cosimplicial objects
    \[ \augcosimp[\mathcal F(S_i)]{\mathcal F(S_i')}{\mathcal F(S_i' \times_{S_i} S_i')}. \]
is $d$-rapidly rapidly convergent. This is a property preserved by colimits, and will therefore yield the equivalence in the second line below:
\begin{align*} \label{eqn:commuting filtered colimits with Tot for proetale descent}
    \varinjlim_i \mathcal F(S_i) &\simeq \varinjlim_i \lim_n \Tot_n \mathcal F({S'_i}^{\times_{S_i} \bullet + 1}) \\
    &\simeq \lim_n \varinjlim_i \Tot_n \mathcal F({S'_i}^{\times_{S_i} \bullet + 1}) \\
    &\simeq \lim_n \Tot_n \varinjlim_i \mathcal F({S'_i}^{\times_{S_i} \bullet + 1}) \simeq \lim_{\bm \Delta} \nu^p \mathcal F({S'_i}^{\times_{S_i} \bullet + 1}). \numberthis
\end{align*}

Suppose therefore that $S' \to S$ is a \emph{finite} covering. Decomposing $S$ into transitive $G$-sets splits (\ref{eqn:Amistur complex for sheaf on et G}) into a product of the \v Cech complexes for $S'_H = S' \times_S G/H \to G/H$, and further writing $S'_H \coloneqq \coprod_i G/K_i$ we will reduce to the case $G/K \to G/H$, where $K \subset H$ are open subgroups.

For this last point, note that if $S'_1, S_2' \to S$ are coverings then $S'_1 \sqcup S'_2 \to S$ factors as $S'_1 \sqcup S'_2 \to S \sqcup S'_2 \to S$. For the first map we have $(S'_1 \sqcup S'_2)^{\times_{S \sqcup S'_2} \bullet} \simeq {S'_1}^{\times_S \bullet} \sqcup {S'_2}$, so $d$-rapid convergence for $S \to S'_1$ implies the same for $S'_1 \sqcup S'_2 \to S \sqcup S'_2$. On the other hand, the second is split, and so the \v Cech complex for $S'_1 \sqcup S'_2 \to S$ is a retract of that for $S'_1 \sqcup S'_2 \to S \sqcup S'_2$. In particular, $d$-rapid convergence for the covering $S'_1 \to S$ implies the same for $S'_1 \sqcup S'_2 \to S$. This allows us to reduce to covers by a single finite $G$-orbit, as desired.

Writing $\operatorname{cd}_p(G) = d < \infty$, we are in the context of \cite[Lemma~4.16]{cm}, which says that hypercompleteness is \emph{equivalent} to condition $(\star)$ when the mod $p$ cohomological dimension is finite. Thus each of the diagrams
    \[ \augcosimp[\mathcal F(G/H)]{\mathcal F(G/K)}{\mathcal F(G/K \times_{G/H} G/K)}  \]
is $d$-rapidly convergent.
\end{proof}




When $p-1 \mid h$, the Morava stabiliser group has $p$-torsion elements, and therefore infinite cohomological dimension mod $p$. To obtain results uniform in $h$ it will therefore be important to know that the same results hold under the more general assumption of \emph{descendability} in the sense of \cite{mathew_galois}.

\begin{defn}
\begin{enumerate}[itemindent = *, leftmargin = \parindent]
    \item Suppose $(\mathcal C, \otimes, \bm 1)$ is symmetric monoidal, and $Z \in \mathcal C$. Then $\Thotimes (Z)$ is the smallest full subcategory of $\mathcal C$ containing $Z$ and closed under extensions, retracts, and $(-) \otimes X$ for every $X \in \mathcal C$. Note that $\Thotimes(Z)$ is the union of full subcategories $\Thotimes_r(Z)$, for $r \geq 1$, spanned by retracts of those objects that can be obtained by at most $r$-many extensions of objects $Z \otimes X$; each of these is a $\otimes$-ideal closed under retracts, but not thick.
    
    \item Suppose now that $A \in \Alg(\mathcal C)$. Recall that $M \in \mathcal C$ is \emph{$A$-nilpotent} (\cite[Definition~7.1.6]{orangebook}) if $M \in \Thotimes(A)$, and that $A$ is \emph{descendable} (\cite[Definition~3.18]{mathew_galois}) if $\Thotimes(A) = \mathcal C$. Thus $A$ is descendable if and only if $\bm 1$ is $A$-nilpotent.
\end{enumerate}
\end{defn}

In \cite[Proposition~3.20]{mathew_galois} it is shown that $A \in \Alg(\mathcal C)$ is descendable if and only if the $\Tot$-tower of the Amitsur complex
    \[ \augcosimp[\bm 1]{A}{A \otimes A} \]
defines a constant pro-object converging to $\bm 1$. It will be useful to have the following quantitative refinement of this result, which explicates some of the relations between various results in \opcit\ with \cite{cm, mathew-naumann-noel, mathew_thick_2015}:

\begin{lem} \label{lem:different types of nilpotence}
    Let $\mathcal C$ be stable and symmetric monoidal. Consider the following conditions:
    \renewcommand{\labelenumi}{$(\arabic{enumi})$}
    \begin{enumerate}
        \item[$(1)_d$] \label{item:nilpotence d-rapidly convergent} The Amitsur complex for $A$ is $d$-rapidly convergent to $\bm 1$.

        \item[$(2)_d$] \label{item:nilpotence retract} The canonical map $\bm 1 \to \Tot_d A^{\otimes \bullet +1}$ admits a retraction.

        \item [$(3)_d$] \label{item:nilpotence exponent} The canonical map $\Tot^d A^{\otimes \bullet +1} \coloneqq \fib( \bm 1 \to \Tot_d A^{\otimes \bullet +1}) \to \bm 1$ is null. In the notation of \cite[\S4]{mathew-naumann-noel}, this says that $\exp_A(\bm 1) = d$.

        \item[$(4)_d$] \label{item:nilpotence vanishing lines} For any $X \in \mathcal C$, the spectral sequence for the tower of mapping spectra
        \begin{equation} \label{eqn:tower of mapping spectra for nilpotence vanishing lines}
            \cdots \to \map_\mathcal C (X, \Tot_n A^{\otimes \bullet}) \to \cdots \to \map_{\mathcal C} (X, \Tot_0 A^{\otimes \bullet})
        \end{equation}
        collapses at a finite page, with a horizontal vanishing line at height $d$.
    \end{enumerate}

    Then we have implications $(1)_d \Leftrightarrow (4)_d$, $(2)_d \Leftrightarrow (3)_d$, and $(1)_d \Rightarrow (2)_d \Rightarrow (1)_{d+1}$.
    
\end{lem}

\begin{proof}
    We begin with the first three conditions. The implication $(1)_d \Rightarrow (2)_d$ is immediate from the diagram
\[\begin{tikzcd}[ampersand replacement=\&,column sep=small]
	{\Tot^d A^{\otimes \bullet + 1}} \& {\bm 1} \& {\Tot_dA^{\otimes \bullet + 1}} \\
	{\Tot^0 A^{\otimes \bullet + 1}} \& {\bm 1} \& {\Tot_ 0A^{\otimes \bullet + 1}}
	\arrow[no head, from=1-2, to=2-2, equals]
	\arrow[from=2-2, to=2-3]
	\arrow[from=1-3, to=2-3]
	\arrow[from=1-2, to=1-3]
	\arrow[from=1-1, to=1-2]
	\arrow[from=1-1, to=2-1, "0"']
	\arrow[from=2-1, to=2-2]
    \arrow[from=1-3, to=2-2, dashed, "\exists"]
\end{tikzcd}\]
    in which the rows are cofibre sequences. The equivalence $(2)_d \Leftrightarrow (3)_d$ is clear, so we now prove $(2)_d \Rightarrow (1)_{d+1}$.
    
    Recall from the proof of \cite[Proposition~3.20]{mathew_galois} that the full subcategory
        \[ \{ X : X \otimes A^{\otimes \bullet +1} \text{ is a constant pro-object} \} \subset \mathcal C \]
    is a thick $\otimes$-ideal containing $A$ (since $A \otimes A^{\otimes \bullet + 1}$ is split), and therefore contains $\bm 1$. We can in fact define full subcategories
        \[ \mathcal C_r \coloneqq \{ X : X \otimes A^{\otimes \bullet +1} \text{ is $r$-rapidly convergent} \} \subset \mathcal C, \]
    which are closed under retracts and $(-) \otimes X$ for any $X \in \mathcal C$. Each of these is \emph{not} thick, but if $X \to Y \to Z$ is a cofibre sequence with $X \in \mathcal C_r$ and $Z \in \mathcal C_{r'}$, then $Y \in \mathcal C_{r + r'}$ (this follows by contemplating the diagram
\[ \numberthis \label{eqn:nilpotence is additive}
\begin{tikzcd}[ampersand replacement=\&,column sep=small]
	{X \otimes \Tot_n A^{\otimes \bullet +1}} \& {Y \otimes \Tot_n A^{\otimes \bullet +1}} \& {Z \otimes \Tot_n A^{\otimes \bullet +1}} \\
	{X \otimes \Tot_{n+r'} A^{\otimes \bullet +1}} \& {Y \otimes \Tot_{n+r'} A^{\otimes \bullet +1}} \& {Z \otimes \Tot_{n+r'} A^{\otimes \bullet +1}} \\
	{X \otimes \Tot_{n+r+r'} A^{\otimes \bullet +1}} \& {Y \otimes \Tot_{n+r+r'} A^{\otimes \bullet +1}} \& {Z \otimes \Tot_{n+r+r'} A^{\otimes \bullet +1}}
	\arrow[from=1-1, to=2-1]
	\arrow["0"', from=2-1, to=3-1]
	\arrow[from=3-1, to=3-2]
	\arrow[from=3-2, to=3-3]
	\arrow[from=2-2, to=3-2]
	\arrow[from=2-1, to=2-2]
	\arrow[from=1-1, to=1-2]
	\arrow[from=1-2, to=1-3]
	\arrow["0", from=1-3, to=2-3]
	\arrow[from=2-3, to=3-3]
	\arrow[from=2-2, to=2-3]
	\arrow[from=1-2, to=2-2]
\end{tikzcd} \]
    as in \cite[Lemma~2.1]{hopkins-palmieri-smith} or \cite[Proposition~3.5]{mathew_thick_2015}). In particular, $\Tot_r A^{\otimes \bullet + 1} \in \mathcal C_{r + 1}$ since this can be constructed iteratively by taking $r + 1$-many extensions by free $A$-modules (while $A \in \mathcal C_1$ by \cref{lem: A otimes cobar complex is 1-rapidly convergent}). Now assumption $(2)_d$ implies that $\bm 1 \in \mathcal C_{d+1}$, that is, $(1)_{d+1}$ holds.

    The implication $(1)_d \Leftrightarrow (4)_d$ is proven in \cite[Proposition~3.12]{mathew_thick_2015} (in the case $\mathfrak U = \mathcal C$), by keeping track of the index $d$ (called $N$ therein).
\end{proof}
    

\begin{remark}
The implications in this lemma are special to the Amitsur/cobar complex. For an arbitrary cosimplicial object (even in spectra), condition $(2)_d$ does not in general imply $(1)_{d'}$ for any $d'$ (although $(1)_d \Rightarrow (2)_d$ still holds).
\end{remark}

\begin{lem} \label{lem: mathcal E of a pullback}
Let $G$ be a profinite group, $\mathcal C$ a stable homotopy theory, and $\bm 1 \to A$ a faithful $G$-Galois extension in $\mathcal C$ corresponding to $\mathcal A^\delta \in \Shh[\CAlg(\mathcal C)]{\et G}$. Given any cospan of $G$-sets
    \[ S_1 \to S_0 \leftarrow S_2 ,\]
we have 
\begin{equation} \label{eqn:mathcal E of pullback is pushout}
    \nu^p \mathcal A^\delta (S_1 \times_{S_0} S_2) = \nu^p \mathcal A^\delta (S_1) \otimes_{\nu^p \mathcal A^\delta(S_0)} \nu^p \mathcal A^\delta(S_2).
\end{equation}
\end{lem}

\begin{proof}
Assume first that each $S_i$ is of the form $G/U_i$ where $U_i$ is an open subgroup. 

We can write $S_1 \times_{S_0} S_2 = \coprod_{g \in U_1 \backslash U_0 / U_2} G / U_1 \cap g^{-1} U_2 g$, so that
    \[ \nu^p \mathcal A^\delta (S_1 \times_{S_0} S_2) = \bigoplus_{U_1 \backslash U_0 / U_2} \mathcal A^\delta(U_1 \cap g^{-1} U_2 g). \]
This admits an algebra map from the pushout $\mathcal A^\delta (S_1) \otimes_{\mathcal A^\delta(S_0)} \mathcal A^\delta(S_2)$, and we can check this map is an equivalence after base-change to $A$, a faithful algebra over $\mathcal A^\delta(S_0)$. But 
\begin{align*}
    A \otimes_{\mathcal A^\delta(S_0)} \mathcal A^\delta(S_1) \otimes_{\mathcal A^\delta(S_0)} \mathcal A^\delta(S_2)  &\simeq \bigoplus_{U_0/U_1} A \otimes_{\mathcal A^\delta(S_0)} \mathcal A^\delta(S_2) \\
    &\simeq \bigoplus_{U_0/U_1} \bigoplus_{U_0/U_2} A \\
    & \simeq \bigoplus_{U_1 \backslash U_0 /U_2} \bigoplus_{U_0/U_1 \cap g^{-1} U_2 g} A \\
    & \simeq A \otimes_{\mathcal A^\delta(S_0)} \bigoplus_{U_1 \backslash U_0 /U_2} \mathcal A^\delta(U_0/U_1 \cap g^{-1} U_2 g),
\end{align*}
using the isomorphism of $U_0$-sets
    \[ U_0/U_1 \times U_0/U_2 \simeq \coprod_{g \in U_1 \backslash U_0 /U_2} U_0/U_1 \cap g^{-1} U_2 g. \]

The result for finite $G$-sets follows from the above by taking coproducts, and for arbitrary $G$-sets by taking cofiltered limits of finite $G$-sets.
\end{proof}





Finally, we prove the desired analogue of \cref{lem:proetale descent for hypercomplete etale sheaves} for descendable Galois extensions.

\begin{prop} \label{lem:nu^p is a sheaf for descendable extensions}
Let $G$ be a profinite group, $\mathcal C$ a stable homotopy theory, and $\bm 1 \to A$ a $G$-Galois extension in $\mathcal C$ corresponding to $\mathcal A^\delta \in \Shh[\CAlg(\mathcal C)]{\et G}$. Suppose moreover that $\bm 1 \to A$ is descendable. Then $\nu^p \mathcal A^\delta \in \mathcal P(\proet G, \mathcal C)$ is a hypercomplete sheaf.
\end{prop}

\begin{proof}
Following the method of \cref{lem:proetale descent for hypercomplete etale sheaves}, to show that $\mathcal A = \nu^p \mathcal A^{\delta}$ is a sheaf it will suffice to show that there is some $d \geq 1$ such that 
    \[ \augcosimp[\mathcal A(G/U)]{\mathcal A(G/V)}{\mathcal A(G/V \times_{G/U} G/V)} \]
is $d$-rapidly convergent for \emph{all} finite coverings $G/V \to G/U$, even when $\cd_p(U) = \infty$. This claim implies the equivalence in the third line of \cref{eqn:commuting filtered colimits with Tot for proetale descent}, and the other identifications are formal. In fact, this will give us hypercompleteness too: this follows from \cite[Proposition~2.25]{cm}, noting that the implication $(2) \Rightarrow (1)$ therein uses only that the site is finitary and no assumption of finite cohomological dimension.

By \cref{lem: mathcal E of a pullback}, the \v Cech complex in question is the cobar complex
\begin{equation} \label{eqn:Amitsur complex for U subset V}
	\augcosimp[A^{hU}]{A^{hV}}{A^{hV} \otimes_{A^{hU}} A^{hV}}.
\end{equation}

For every $r \geq 1$, we shall consider the following full subcategory of $\mathcal D \coloneqq \Mod_{A^{hU}}(\mathcal C)$:
    \[ \mathcal D_r = \mathcal D_r(U, V) \coloneqq \{X: X \otimes_{A^{hU}} (A^{hV})^{\otimes_{A^{hU}} \bullet + 1} \text{ is $r$-rapidly convergent} \}. \]
Our aim is to show that $A^{hU} \in \mathcal D_{d}$. To this end, we observe the following properties of $\mathcal D_r$:
\begin{enumerate}
    \item $\mathcal D_r$ is closed under retracts.
    \item $\mathcal D_r$ is a $\otimes$-ideal: if $X \in \mathcal D_r$ then $X \otimes Y \in \mathcal D_r$ for any $Y$.
    \item $\mathcal D_r$ is not thick. However, if $X \to Y \to Z$ is a cofibre sequence such that $X \in \mathcal D_r$ and $Z \in \mathcal D_{r'}$, then $Y \in \mathcal D_{r + r'}$.
\end{enumerate}

But $A \in \mathcal D_1$, since by \cite[Example~4.7.2.7]{ha} the cobar complex \cref{eqn:Amitsur complex for U subset V} splits after applying
    \[ A \otimes_{A^{hU}} (-) \simeq A \otimes_{A^{hV}} A^{hV} \otimes_{A^{hU}} (-). \]
By \cite[Corollary~4.13]{mathew-naumann-noel} we have
    \[ \exp_{A \otimes A^{hU}}(A^{hU}) \leq \exp_A(\bm 1) = d. \]
But $A^{hU}$ splits over $A$, with $A \otimes A^{hU} \simeq \prod_{G/U} A \simeq \bigoplus_{G/U} A \in \Mod_A(\mathcal C)$. Thus any free $A \otimes A^{hU}$-module is also free over $A$, and $\exp_{A \otimes A^{hU}}(A^{hU}) \geq \exp_A(A^{hU})$. Thus $A^{hU} \in \mathcal D_d$ as desired.
\end{proof}

\subsubsection{Functoriality in \texorpdfstring{$G$}{G}}

Let $f \colon H \to G$ be a continuous homomorphism between profinite groups. In this section, we will consider the functoriality of the construction
    \[ G \mapsto \Shh[\mathcal C]{\proet G}. \]
For example, taking $p \colon G \to *$ we obtain (homotopy) invariants and coinvariants functors; functoriality will imply the iterated fixed points formula
    \[ (X^{hH})^{hG/H} \simeq X^{hG} \]
for $H \subset G$ a normal subgroup. 

\begin{remark}
For any presentable $\mathcal C$ we have $\Sh[\mathcal C]{\proet G} \simeq \mathcal Sh(\proet G) \otimes \mathcal C$ by \cite[Remark~1.3.1.6]{sag}, and so we restrict ourselves to sheaves of spaces.    
\end{remark}

\begin{prop} \label{lem:functoriality of Sh(proet)}
\begin{enumerate}
    \item Any map of profinite groups $f \colon H \to G$ gives rise to a geometric morphism
        \[ f^* : \mathcal Sh(\proet G) \rightleftarrows \mathcal Sh(\proet H) : f_* \]
    This defines a functor $\mathrm{Grp}(\Profin) \to \mathcal{LT}\mathrm{op}_\infty$ to the \inftycat{} of $\infty$-topoi and left adjoints.

    \item The left adjoint $f^* \colon \mathcal Sh(\proet G) \to \mathcal Sh(\proet H)$ admits a further left adjoint $f_!$.
\end{enumerate}
\end{prop}

\begin{proof}
Given $f \colon H \to G$, restriction determines a morphism of sites
    \[ \operatorname{res}_f : \proet G \to \proet H, \]
hence a geometric morphism
    \[ f^* : \mathcal Sh(\proet G) \rightleftarrows \mathcal Sh(\proet H) : f_* \]
whose right adjoint is given by precomposition with $\operatorname{res}_f$.

For $(ii)$, we apply \cref{lem:morphism of sites with left adjoint induces f_!} to the adjunction
    \[ (-) \times_H G : \proet H \rightleftarrows \proet G : \mathrm{res}_f, \]
noting that if $S \to S'$ is a surjection of $H$-sets then $S \times_H G \to S' \times_H G$ is also surjective.
\end{proof}

\begin{lem} \label{lem:morphism of sites with left adjoint induces f_!}
Given a morphism of sites $f \colon \mathcal C \to \mathcal C'$, if $f$ admits a left adjoint $g$ then $f^*$ is given by the sheafification of $g_*$, 
    \[ f^* \mathcal F = (g_* \mathcal F)^{+}. \]
If $g$ is itself a morphism of sites, then $f^* = g_*$ admits a further left exact left adjoint $f_! = g^*$.
\end{lem}

\begin{proof}
In general, $f^*$ is the sheafification of the presheaf extension $f^p = \mathrm{Lan}_f(-)$. Since $f^p \colon \mathcal P(\mathcal C) \to \mathcal P(\mathcal C')$ is left adjoint to $f_*$ \cite[Proposition~4.3.3.7]{htt}, precomposition with the counit $\eta$ of the $g \dashv f$ adjunction is a map
    \[ g_* f_* \simeq (gf)_* \xrightarrow{\eta_*} id_{\mathcal P(\mathcal C')}, \]
and hence induces a natural transformation $\alpha \colon g_* \to f^p$ \cite[Lemma~2.3.7]{riehl-verity}. On vertices,
    \[ f^p \mathcal F : C' \mapsto \varinjlim_{\mathcal C \downarrow C'} \mathcal F(C), \]
and the above natural transformation is the map $\mathcal F(g(C')) \to f^p \mathcal F(C)$ induced by the counit $gf(C') \to C'$, viewed as an object of $\mathcal C \downarrow C'$. But the counit is terminal in $\mathcal C \downarrow C'$, so
    \[ \alpha_{C'} : g_* \mathcal F(C') \xrightarrow{\sim} f^p \mathcal F (C'), \]
which gives the first claim after sheafifying. If moreover $g$ preserves coverings then $g_* \mathcal F$ is already a sheaf, and so $f^* \simeq g_*$ admits $g^*$ as a left adjoint.
\end{proof}

\begin{lem}
Suppose that $i \colon H \subset G$ is a normal subgroup with quotient $p \colon G \to G/H$, and $X \in \Sh[\mathcal C]{\proet G}$. Then
    \[ (X^{hH})^{hG/H} \simeq X^{hG}. \]
\end{lem}

\begin{proof}
Writing $\Gamma_G$ for evaluation on $G/G \in \proet G$ and likewise for $G/H$, the claim is that 
    \[ \Gamma_{G/H} p_* X \simeq \Gamma_G X. \]
This is immediate: if $\pi^G \colon G \to *$ denotes the projection, then \cref{lem:functoriality of Sh(proet)} implies the middle equivalence in
    \[ \Gamma_G \simeq \Gamma_* \pi^G_* \simeq \Gamma_* 
    \pi^{G/H}_* p_* \simeq \Gamma_{G/H} p_* . \qedhere \]
\end{proof}

\subsubsection{Comparison with pyknotic \texorpdfstring{$G$}{G}-objects}

Given a profinite group $G$ and a base \inftycat{} $\mathcal C$, we have described the \inftycat{} $\widehat{\mathcal Sh}(\proet G)$ as a good model for the $\infty$-topos of continuous $G$-spaces. Using the theory of pyknotic objects, there are two other natural candidates:
\begin{enumerate}
    \item $G$ has image $\underline{G} \in \Pyk(\Set)$ under the embedding of compactly generated spaces in pyknotic sets, and $\underline G$ is a pyknotic group. Since $\mathcal S$ is tensored over sets (and hence $\Pyk(\mathcal S)$ over $\Pyk(\Set)$), we can form a category of $\underline G$-modules in $\Pyk(\mathcal S)$.

    \item $G$ has a \emph{pyknotic classifying space} $BG \in \Pyk(\mathcal S)$, given by the limit in $\Pyk(\mathcal S)$ of the classifying spaces $BG_i \in \mathcal S \hookrightarrow \Pyk(\mathcal S)$, where $G = \lim G_i$ is a presentation of $G$ as a profinite group. We can therefore consider the $\infty$-topos $\Pyk(\mathcal S)_{/BG}$.
\end{enumerate}

In this section we show that all these notions agree. See \cite[Theorem~A.5.6.1]{sag} for a related result in the context of \emph{pro-$\pi$-finite spaces}.

\begin{remark}
We have chosen to use the $\infty$-topos $\Pyk(\mathcal S)$ in this section, as some of our results depend on \cite{wolf}. Note however that a translation of these results to $\Cond(\mathcal S)$ is given in \cite[Remark~4.19]{wolf}.
\end{remark}

\begin{prop} \label{lem:proet is G-objects}
For any profinite group $G$, the adjunction
    \[ i^* : \widehat{\mathcal Sh}(\proet G) \rightleftarrows \widehat{\mathcal Sh}(\proet *) = \Pyk(\mathcal S) : i_* \]
induced by $i \colon * \to G$ is comonadic, and induces an equivalence
    \[ \widehat{\mathcal Sh}(\proet G) \simeq \Mod_{\underline{G}}(\Pyk(\mathcal S)). \]
In particular, for any presentable $\mathcal C$ we have 
    \[ \Shh[\mathcal C]{\proet G} \simeq \Mod_{\underline{G}}(\Pyk(\mathcal C)). \]
\end{prop}

\begin{proof}
To see that the adjunction is monadic, we must verify the conditions of the Barr-Beck-Lurie theorem \cite[Theorem~4.7.3.5]{ha}; in fact, we will verify that this holds before hypercompletion. Since both sides are presentable, this amounts to checking:
\begin{enumerate}
    \item $i^*$ is conservative,
    \item $i^*$ preserves totalisations of $i^*$-split simplicial objects.
\end{enumerate}
But $i^*$ admits a further left adjoint, which implies $(ii)$. On the other hand, $i^* \mathcal F = i^p \mathcal F \colon T \mapsto \mathcal F(T \times G)$ by \cref{lem:morphism of sites with left adjoint induces f_!}, so that $i^*$ is conservative by \cref{lem:Free generates proet}.

As a result, we see that $\mathcal Sh(\proet G) \simeq \mathrm{Mod}_{\mathbb T}(\Pyk(\mathcal S))$, where $\mathbb T$ is the comonad $i^* i_*$. To identify this with the coaction comonad for $\underline G$ it suffices to do so for compact generators of $\mathcal Sh(\proet *)$, since both comonads preserve colimits. But if $T \in \proet *$ then
    \[ \mathbb T(\underline T)(T') = \Cont(\operatorname{res}_i (T' \times G), T) \cong \Cont(\operatorname{res}_i (T'), \Cont(\operatorname{res}_i(G), T)) = \Map(\underline G, \underline T)(T')  \]
since $\proet *$ is Cartesian closed; that is, $\mathbb T = \Map(\underline G, -)$ so
    \[ \mathcal Sh(\proet G) \simeq \Mod_{\underline G}(\mathcal Sh(\proet *)). \]
For general $\mathcal C$, this implies the equivalences
\begin{align*}
    \Sh[\mathcal C]{\proet G} &\simeq \mathcal Sh(\proet G) \otimes \mathcal C \\
    &\simeq \Mod_{\underline G}(\Pyk(\mathcal S)) \otimes \mathcal C \\
    &\simeq \Mod_{\underline G}(\Pyk(\mathcal S)) \otimes_{\Pyk(\mathcal S)} \Pyk(\mathcal S) \otimes \mathcal C \\
    &\simeq \Mod_{\underline G}(\Pyk(\mathcal S)) \otimes_{\Pyk(\mathcal S)} \Pyk(\mathcal C) \\
    &\simeq \Mod_{\underline G}(\Pyk(\mathcal C)).
\end{align*}
Here we have used \cite[Remark~2.3.5]{pyknotic} for the fourth equivalence, and \cite[Theorem~4.8.4.6]{ha} for the final equivalence.
\end{proof}

\begin{prop} \label{lem:proet is Pyk over BG}
For any profinite group $G$, there are natural equivalences
    \[ \widehat{\mathcal Sh}(\proet G) \simeq \mathrm{Fun}^{\mathrm{cts}}(BG, \Pyk(\mathcal S)) \simeq \Pyk(\mathcal S)_{/BG}. \]
\end{prop}

\begin{proof}
This is a very mild modification of \cite[Corollary~1.2]{wolf}. In the trivially stratified case, the main theorem of \opcit\ provides a natural \emph{exodromy} equivalence
    \[ \mathcal X^{\Pyk} \coloneqq \widehat{\mathcal Sh}_{\mathrm{eff}}(\Pro(\mathcal X^{\mathrm{coh}}_{< \infty})) \xrightarrow[\mathrm{ex}]{\sim} \mathrm{Fun}^{\mathrm{cts}}(\widehat{\Pi}_\infty \mathcal X, \Pyk(\mathcal S)) \numberthis \label{eqn:pyknotic exodromy} \]
between the \emph{pyknotification} of an $\infty$-topos $\mathcal X$ with the \inftycat{} of \emph{pyknotic presheaves} on its profinite shape, viewed as a pyknotic space\footnote{Recall \cite[Remark~13.4.4]{exodromy} that the inclusion $y \colon \mathcal S_\pi^\wedge \hookrightarrow \widehat{\mathcal Sh}_{\mathrm{eff}}(\mathcal S_\pi^\wedge) \simeq \Pyk(\mathcal S)$ preserves all small limits. Combining this with \cite[13.4.11]{exodromy} we see that $y \{ X_i \} = \lim y X_i = \lim \Gamma^* X_i$ is the limit in $\Pyk(\mathcal S)$ of the $\pi$-finite spaces $X_i$, viewed as discrete pyknotic spaces.}. Any $X \in \Pro(\mathcal Sh(\et G)^{\mathrm{coh}}_{< \infty})$ may be covered by a profinite $G$-set: this follows just as in the case $G = *$, which is \cite[Proposition~13.4.9]{exodromy}. Thus the subcategory $j \colon \proet G = \Pro(\et G) \hookrightarrow \Pro(\mathcal Sh(\et G)^{\mathrm{coh}}_{< \infty})$ generates the same topos, and we obtain natural equivalences
    \[ \widehat{\mathcal Sh}(\proet G) \xleftarrow[j_*]{\sim} \widehat{\mathcal Sh}(\Pro(\mathcal Sh(\et G)^{\mathrm{coh}}_{< \infty})) = \mathcal Sh(\et G)^{\Pyk} \xrightarrow[\mathrm{ex}]{\sim} \mathrm{Fun}^{\mathrm{cts}}(\widehat{\Pi}_\infty \mathcal Sh(\et G), \Pyk(\mathcal S)), \]
where $j_*$ is restriction and $\mathrm{ex}$ denotes the pyknotic exodromy equivalence \cref{eqn:pyknotic exodromy}. We now identify $\widehat{\Pi}_\infty \mathcal Sh(\et G) \in \mathcal S_\pi^\wedge = \Pro(\mathcal S_\pi)$: if we write $G = \lim G_i$, then by \cite[Construction~4.5]{cm} there is an equivalence
    \[ \mathcal Sh(\et G) \simeq \varprojlim \mathcal S_{/BG_i}. \]
By \cite[Theorem~E.2.4.1]{sag}, the pro-extension $\Psi \colon \mathcal S_\pi^\wedge \to \mathcal{LT}\mathrm{op}_\infty$ of $X \mapsto \mathcal S_{/X}$ is a fully faithful right adjoint to $\widehat{\Pi}_\infty$, so
    \[ \widehat{\Pi}_\infty \mathcal Sh(\et G) \simeq \widehat{\Pi}_\infty \varprojlim \Psi(BG_i) \simeq \widehat{\Pi}_\infty \Psi(\{BG_i\}) \simeq \{BG_i\} \in \mathcal S_\pi^\wedge. \]

By definition, $BG$ is image of $\{BG_i\}$ in $\Pyk(\mathcal S)$, so continuous unstraightening \cite[Corollary~3.20]{wolf} yields the desired equivalence
    \[ \widehat{\mathcal Sh}(\proet G) \simeq \mathrm{Fun}^{\mathrm{cts}}(BG, \Pyk(\mathcal S)) \simeq \Pyk(\mathcal S)_{/BG}. \qedhere \]
\end{proof}

\begin{remark}
The results of \cite[Appendix~E.5]{sag} imply that the profinite group $G$, as a $0$-truncated object of $\mathrm{Grp}(\mathcal S_\pi^\wedge)$, has a profinite classifying space; that is, a delooping in $\mathcal S_\pi^\wedge$. In fact, in the proof of \cite[Theorem~E.5.0.4]{sag} (\cite[E.5.6]{sag}) Lurie observes that this is $\{ BG_i \}$; passing to $\Pyk(\mathcal S)$ we see that $BG$ is a delooping of $\underline G$. In the hypercomplete case, one can use this to deduce \cref{lem:proet is G-objects} from \cref{lem:proet is Pyk over BG}.
\end{remark}

\begin{cor} \label{lem:6 functors}
For any presentable $\mathcal C$, the assignment $G \mapsto \Shh[\mathcal C]{\proet G}$ extends to a Wirthm\"uller context: there is a $6$-functor formalism on the geometric context of profinite groups and continuous homomorphisms, with $f^! \simeq f^*$ for any $f$.
\end{cor}

\begin{proof}
By \cite[Prop.~A.5.10]{mann}, it remains to check:
\begin{enumerate}
    \item For every map $f \colon H \to G$ of profinite groups, the functor $f^*$ admits a left adjoint $f_! \colon \widehat{\mathcal Sh}(\proet H) \to \widehat{\mathcal Sh}(\proet G)$.
    
    \item If moreover $g \colon G' \to G$ and $H' \coloneqq H \times_G G'$, then the following square is left adjointable:
        \[\begin{tikzcd}
	   {\widehat{\mathcal Sh}(\proet G)} & {\widehat{\mathcal Sh}(\proet H)} \\
	   {\widehat{\mathcal Sh}(\proet G')} & {\widehat{\mathcal Sh}(\proet H')}
	   \arrow["{g^*}"', from=1-1, to=2-1]
	   \arrow["{f^*}", from=1-1, to=1-2]
	   \arrow["{g'^*}", from=1-2, to=2-2]
	   \arrow["{f'^*}"', from=2-1, to=2-2]
        \end{tikzcd}\]
    That is, the push-pull transformation $f'_! g'^* \to g^* f_!$ is an equivalence.

    \item If $M \in \widehat{\mathcal Sh}(\proet G)$ and $N \in \widehat{\mathcal Sh}(\proet H)$, the natural map
        \[ f_!(N \times f^*M) \to f_! N \times M \]
    is an equivalence.
\end{enumerate}

In fact, these properties will all follow formally from \cref{lem:proet is Pyk over BG}. The functor $G \mapsto \widehat{\mathcal Sh}(\proet G)$ is the restriction to classifying spaces of profinite groups of the functor $X \mapsto \Pyk(\mathcal S)_{/X}$, which in the language of \cite{martini_yoneda} is the \emph{universe} $\Omega_{\Pyk(\mathcal S)}$. But the universe $\Omega$ of any $\infty$-topos $\mathcal B$ satisfies properties $(i)$ to $(iii)$. Explicitly:

\begin{enumerate}
    \item For any $f \colon A \to B$ in $\mathcal B$, the \'etale geometric morphism $f^*: \mathcal B_{/B} \to \mathcal B_{/A}$ is given by basechange along $f$, and admits a left adjoint $f_! = f \circ - $.

    \item For a pullback square
\[\begin{tikzcd}[ampersand replacement=\&]
	A \& B \\
	C \& D
	\arrow["{g'}", from=1-1, to=1-2]
	\arrow["f", from=1-2, to=2-2]
	\arrow["g"', from=2-1, to=2-2]
	\arrow["{f'}"', from=1-1, to=2-1]
\end{tikzcd}\]
in $\mathcal B$, the evaluation of the push-pull transformation $f_!' {g'}^* \to g^* f_!$ on $X \to B$ may be identified with the left-hand vertical map in the extended diagram
\[\begin{tikzcd}[ampersand replacement=\&]
	{X \times_B A} \& A \& B \\
	{X \times_D C} \& C \& D
	\arrow["{g'}", from=1-2, to=1-3]
	\arrow["f", from=1-3, to=2-3]
	\arrow["g"', from=2-2, to=2-3]
	\arrow["{f'}"', from=1-2, to=2-2]
	\arrow[from=1-1, to=1-2]
	\arrow[from=2-1, to=2-2]
	\arrow[from=1-1, to=2-1]
\end{tikzcd}\]
    and is therefore an equivalence.

    \item For $f \colon A \to B$ in $\mathcal B$, the projection formula is given in \cite[Remark~6.3.5.12]{htt}. \qedhere
\end{enumerate}

\end{proof}

\subsection{Morava E-theory as a pro\'etale spectrum}
The aim of this section is to show that the pro\'etale site allows us to capture the continuous action on $\E$ as a sheaf of spectra on $\proet \G$. Our first task is to exhibit Morava E-theory itself as a pro\'etale sheaf of spectra, so as to recover the $\K$-local $\E$-Adams spectral sequence as a descent spectral sequence. This will allow us to compare it to the descent spectral sequence for the Picard spectrum.

\begin{prop} \label{lem:nu^p mathcal E is a sheaf}
    The presheaf of $\K$-local spectra
        \[ \mathcal E \coloneqq \nu^p \mathcal E^\delta: S = \varprojlim_i S_i \mapsto L_{\K} \varinjlim_i \mathcal E^\delta (S_i) \]
    is a hypercomplete sheaf on $\proet \G$.
\end{prop}

\begin{remark}
Since $\mathcal Sp_{\K} \subset \mathcal Sp$ is a right adjoint, the same formula defines a sheaf of spectra (or even of $\mathbb E_\infty$-rings). We stress however that \textbf{$\nu^p$ will always refer to the left Kan extension internal to $\mathcal Sp_{\K}$, or equivalently the $\K$-localisation of the Kan extension in spectra}.
\end{remark}

\begin{proof}
This is immediate from \cref{lem:nu^p is a sheaf for descendable extensions}: by \cite[Proposition~10.10]{mathew_galois}, Morava E-theory $\E \in \CAlg(\mathcal Sp_{\K})$ is descendable. This is a consequence of descendability of $L_h \mathbb S \to \E$, which is proven in \cite{orangebook} (the proof crucially uses the fact that $\G$ has finite \emph{virtual} cohomological dimension). In fact, when $p$ is sufficiently large (or more generally, when $p-1$ does not divide $h$), \cref{lem:nu^p mathcal E is a sheaf} is a consequence of the descent results of \cite{cm}.
\end{proof}

As an aside, note that the proof of \cref{lem:nu^p is a sheaf for descendable extensions} allows us likewise to consider $\E$-homology of any spectrum $X$:

\begin{cor} \label{lem: mathcal E otimes X is a hypercomplete sheaf}
If $X$ is any $\K$-local spectrum, then the presheaf $\mathcal E \otimes X$ is a hypercomplete sheaf.
\end{cor}

\begin{proof}
Given a covering $S' \to S$ in $\proet \G$, the proof of \cref{lem:nu^p is a sheaf for descendable extensions} showed that the Amitsur complex for $S' \to S$ is $d$-rapidly convergent. Since the functor $(-) \otimes X: \mathcal Sp_{\K} \to \mathcal Sp_{\K}$ preserves finite limits, the same is true for the augmented cosimplicial object given by tensoring everywhere by $X$.
\end{proof}

We can now define a descent spectral sequence for the sheaf $\mathcal E$. 
\begin{remark}
\cref{lem:nu^p mathcal E is a sheaf} gives us a sheaf $\mathcal E \in \Shh[\mathcal Sp_{\K}]{\proet \G}$. On applying the forgetful functor $\mathcal Sp_{\K} \hookrightarrow \mathcal Sp$, we obtain $\mathcal E \in \Shh[\mathcal Sp]{\proet \G}$; note that this sheaf is \emph{not} left Kan extended from $\et \G$. Nevertheless, applying \cref{lem:proetale descent spectral sequence}, we obtain the following result.
\end{remark}

\begin{cor}
    There is a conditionally convergent spectral sequence of the form
    \begin{equation} \label{eqn:descent spectral sequence for nu E}
        E_{2,+}^{s,t} = H^{s}(\proet \G, \pi_{t} \mathcal E) \implies \pi_{t-s} \Gamma \lim \tau_{\leq j} \mathcal E.
    \end{equation}
\end{cor}

To identify the $E_2$-page and the abutment, we need to identify the homotopy sheaves of $\mathcal E$.

\begin{lem} \label{lem:pi* nu* mathcal E of free G-sets}
    The homotopy sheaves of $\mathcal E$ are given 
    by
        \begin{equation}
        \pi_t \mathcal E: S \mapsto 
        \operatorname{Cont}_\G(S, \pi_t \E),
        \end{equation}
    continuous equivariant maps from $S$ to the homotopy groups of Morava E-theory (equipped with their profinite topology).
\end{lem}

\begin{proof}
To prove the lemma, it is enough to prove that the homotopy \textit{presheaves} $\pi_t^p \mathcal E$ take the form $\operatorname{Cont}_\G(-, \pi_t \E)$, since the topology is subcanonical. Moreover, since $\Free \G$ generates the pro\'etale topos, it is enough to show that $(\pi_t^p \mathcal E)|_{\Free \G}: S \mapsto \operatorname{Cont}_\G(S, \pi_t \E) = \operatorname{Cont}(S/\G, \pi_t \E)$.

Using Lemma \ref{lem: mathcal E of a pullback}, we have for any free $\G$-set of the form $S = T \times \G$ (with trivial action on $T$),
\begin{equation}
    \mathcal E (S) \simeq \mathcal E (T) \otimes \mathcal E( \G) \simeq L_{\K} \varinjlim_j \prod_{T_j} \E. \label{eqn:mathcal E of free G-sets}
\end{equation}

Since $\E$-localisation is smashing, the spectrum $\varinjlim_j \prod_{T_j} \E$ is $\E$-local, and so its $\K$-localisation can be computed by smashing with a tower of generalised Moore spectra $M_I$; see for example \cite{hovey-strickland}. Thus
    \[ \pi_t \left( \mathcal E (S) \right) = \pi_t \lim_I (\varinjlim_j \prod_{T_j} \E \otimes M_I ), \]
and we obtain a Milnor sequence
    \[ 0 \to \lim_I  \pi_t \varinjlim_j \prod_{T_j} \E \otimes M_I \to \pi_t \left( \mathcal E (S) \right) \to {\lim_I}^1 \pi_t \varinjlim_j \prod_{T_j} \E \otimes M_I \to 0. \]
Now observe that
\begin{align*}
    \pi_t \varinjlim_j \prod_{T_j} \E \otimes M_I &= \varinjlim_j \prod_{T_j} \pi_t \E \otimes M_I \\
    &= \varinjlim_j \prod_{T_j} (\pi_t \E) /I \\
    &= \varinjlim_j \operatorname{Cont}(T_j, (\pi_t \E) /I ) \\
    &= \operatorname{Cont}(T, (\pi_t \E) /I ) ,
\end{align*}
using for the last equality that the target is finite. In particular, each $\operatorname{Cont}(T, (\pi_t \E) /I') \to \operatorname{Cont}(T, (\pi_t \E) /I )$ is surjective, since the inclusion $I' \subset I$ induces a surjection of finite sets $(\pi_t \E)/I' \twoheadrightarrow (\pi_t \E)/I$, and so admits a (set-theoretic) splitting. Thus $\lim^1$ vanishes and
    \[ \pi_t \left( \mathcal E (S) \right) = \lim_I \operatorname{Cont}(T, (\pi_t \E) /I ) = \operatorname{Cont}(T, \pi_t \E ). \qedhere \] 
\end{proof}

\begin{cor} \label{lem:hyperdescent for presheaf extension of mathcal E}
The Postnikov tower of the sheaf $\mathcal E$ converges. Thus \cref{eqn:descent spectral sequence for nu E} converges conditionally to $\pi_* \Gamma \mathcal E = \pi_* \bm 1_{\K}$.
\end{cor}

\begin{proof}
Both properties are local, so we can restrict to the subsite $\Free \G$. There, the proof of \cref{lem:pi* nu* mathcal E of free G-sets} showed that the homotopy presheaves of $\mathcal E$ are its homotopy sheaves. Taking cofibres and limits of presheaves (both of which preserve sheaves), we see that the truncation $\tau_{\leq t} \mathcal E$ is the presheaf $U \mapsto \tau_{\leq t} \mathcal E(U)$ (that is, no sheafification is necessary). But Postnikov towers in \emph{presheaves} of spectra converge since the same is true in spectra; thus 
	\[ \mathcal E \simeq \lim_t \tau_{\leq t} \mathcal E. \qedhere \]
\end{proof}

\begin{remark}
In fact, we will compare this spectral sequence to the $\K$-local $\E$-Adams spectral sequence (\cref{lem:decalage of descent spectral sequence}) to show existence of a horizontal vanishing line (at least for $t \geq 2$); thus the spectral sequence converges completely in this region.
\end{remark}

\begin{remark}
Similar arguments are used in \cite[Theorem~3.2.3]{bs} and \cite[Proposition~A.10]{mathew_tr} to prove Postnikov completeness results. Specifically, Bhatt and Scholze prove that Postnikov towers of hypercomplete objects in the $\infty$-topos of $\proet \G$ converge. Using this fact, one can also deduce \cref{lem:hyperdescent for presheaf extension of mathcal E} directly from \cref{lem:nu^p mathcal E is a sheaf}; in fact, this proves Postnikov completeness of $\mathcal E \otimes X$ for any spectrum $X$. We remark that it was proven in \cite{mondal-reinecke} that Postnikov towers converge in the hypercomplete topos of an arbitrary \emph{replete} site, answering \cite[Question~3.1.12]{bs}.
\end{remark}

\begin{cor} \label{lem:E1 page of descent spectral sequence}
The starting page of the spectral sequence (\ref{eqn:descent spectral sequence for nu E}) is given by continuous group cohomology:
\begin{equation}
    E_{2,+}^{s,t} = H^{s}(\G, \pi_t \E).
\end{equation}
\end{cor}

\begin{proof}
Using the identification in \cref{lem:pi* nu* mathcal E of free G-sets} of the homotopy sheaves, this follows from \cite[Lemma~4.3.9(4)]{bs}, which implies that the canonical map
    \[ \Phi_M: H^*(\G, M) \to H^*(\proet \G, \operatorname{Cont}_G(-, M))\]
for a topological $\G$-module $M$ is an isomorphism whenever $M$ can be presented as the limit of a countable tower of finite $\G$-modules.
\end{proof}

\begin{remark}[c.f. \cite{barthel-heard}] \label{rem:E2 page of ASS for arbitrary X}
The proof of \cref{lem:pi* nu* mathcal E of free G-sets} also goes through for the hypercomplete sheaf $\mathcal E \otimes X$, as long as $\pi_* \E \wedge X \wedge M_I \simeq (\E_* X)/I$ for each of the ideals $I$. Thus one obtains a conditionally convergent descent spectral sequence
    \[ E_{2}^{s,t}(X) = H^{s} \left( \proet \G, \Cont_\G(-,\E_t^\vee X) \right) \implies \pi_{t-s} X. \]
for any such $\K$-local spectrum $X$, and if each of the $\G$-modules $\E_t X$ satisfies one of the assumptions of \cite[Lemma~4.3.9]{bs}, then the $E_2$-page is given by the continuous group cohomology $H^{s}(\G, \E_t^\vee X)$. For example, this happens when $X$ is finite. On the other hand, writing $\mathcal E_t^\vee X \coloneqq \pi_t (\mathcal E \otimes X)$, one obtains a conditionally convergent spectral sequence
    \[ E_2^{s,t} = H^s(\proet \G, \mathcal E_t^\vee X) \implies \pi_{t-s} X \]
for \emph{any} $\K$-local spectrum, which gives one indication that working with pro\'etale cohomology may be preferable.
\end{remark}

We have thus identified the $E_2$-page of the descent spectral sequence for $\mathcal E$ with the $E_2$-page of the $\K$-local $\E$-Adams spectral sequence. The next step is to show this extends to an identification of the spectral sequences. We do this by using the d\'ecalage technique originally due to Deligne \cite{hodgeii}; the following theorem is standard, but for the sake of completeness (and to fix indexing conventions, one of the great difficulties in the subject) we include the argument in \cref{app:decalage}.

\begin{prop}[\cref{lem:decalage for sheaves of spectra appendix}] \label{lem:decalage for sheaves of spectra}
    Let $\mathcal F$ be a sheaf of spectra on a site $\mathcal C$, and let $X \twoheadrightarrow *$ be a covering of the terminal object. Suppose that for every $t$ and $q > 0$ we have $\Gamma(X^q, \tau_t \mathcal F) = \tau_t \Gamma(X^q, \mathcal F).$ Then there is an isomorphism between the descent and Bousfield-Kan spectral sequences, up to reindexing: for all $r$,
        \[ E_r^{s,t} \cong \check E_{r+1}^{2s-t,s}. \]
\end{prop}

\begin{prop} \label{lem:decalage of descent spectral sequence}
D\'ecalage of the Postnikov filtration induces an isomorphism between the following spectral sequences:
\begin{align*}
    E_{2,+}^{s,t} = \pi_s \Gamma \tau_t \mathcal E = H^{s}(\G, \pi_t \E) &\implies \pi_{t-s} \bm 1_{\K}, \phantom{E_{2,+}^{s,t} = \pi_s \Gamma \tau_t \mathcal E =} \\
    \check E_{3,+}^{2s-t,s} = \pi^{s}(\pi_t \E^{\otimes \bullet + 1}) = H^{s}(\G, \pi_t \E) & \implies \pi_{t-s} \bm 1_{\K}.
\end{align*}
The first is the descent spectral sequence for the sheaf $\mathcal E$, and the second is the $\K$-local $\E$-Adams spectral sequence.
\end{prop}

\begin{proof}
By Lemma \ref{lem: mathcal E of a pullback}, the $\Tot$-filtration associated to the cosimplicial spectrum $\Gamma (\G^\bullet, \mathcal E)$ is precisely the Adams tower for the Amitsur complex of $\bm 1_{\K} \to \E$. The resulting spectral sequence is the $\K$-local $\E$-Adams spectral sequence by definition.

According to Lemma \ref{lem:decalage for sheaves of spectra}, all that remains to check is that each spectrum
    \[ \Gamma (\G^q, \tau_t \mathcal E), \qquad q>0 \]
is Eilenberg-Mac Lane. Note that when $q=1$ this is immediate since any cover of $\G$ in $\proet \G$ is split; when $q>1$ the profinite set $\G^{q-1}$ is not extremally disconnected, and we will deduce this from Lemma \ref{lem:pi* nu* mathcal E of free G-sets}. Indeed, we know that $\Gamma (\G^q, \tau_t \mathcal E)$ is $t$-truncated, while for $s \geq t$ we have
\begin{align*}
    \pi_s \Gamma(\G^q, \tau_t \mathcal E) &= H^{t-s}({\proet \G}_{/\G^q}, \pi_t \mathcal E) \\
    &= H^{t-s}(\mathrm{Profin}_{/\G^{q-1}}, \operatorname{Cont}(-, \pi_t \E)) \\
    &= H^{t-s}_{\mathrm{cond}}(\G^{q-1},  \operatorname{Cont}(-, \pi_t \E)),
\end{align*}
where $H^*_{\mathrm{cond}}$ denotes condensed cohomology. Here we have used the equivalence ${\proet \G}_{/\G^q} \simeq \Profin_{/\G^{q-1}}$ sending $S \mapsto S/\G$. We now argue that the higher cohomology groups vanish, essentially as in the first part of \cite[Theorem~3.2]{condensed}. Namely, the sheaves $\operatorname{Cont}(-, \pi_t \E)$ on ${\proet \G}_{/\G^q}$ satisfy the conditions of \cite[Lemma~4.3.9(4)]{bs}, and so the cohomology groups in question can be computed by \v Cech cohomology: the \v Cech-to-derived spectral sequence collapses, since the higher direct images of $\operatorname{Cont}(-, \pi_t \E)$ vanish. As a result, to check they vanish it will be enough to check that the \v Cech complex 
    \[ \operatorname{Cont}_\G(\G^{q}, \pi_t \E) \to \operatorname{Cont}_\G(S, \pi_t \E) \to \operatorname{Cont}_\G(S \times_{\G^q} S, \pi_t \E) \to \cdots \numberthis \label{eqn:complex for showing Gamma Gn is E-M}\]
is exact for any surjection $S \twoheadrightarrow \G^{q}$. Writing this as a limit of surjections of finite $\G$-sets $S_i \twoheadrightarrow S'_i$ (with $\lim S'_i = \G^q$ and $\lim S_i = S$), and writing $A_{i,I}^j \coloneqq \operatorname{Cont}_\G({S_i}^{\times_{S'_i} j}, \pi_t \E/I)$ for brevity, (\ref{eqn:complex for showing Gamma Gn is E-M}) is the complex
    \[ \lim_I \colim_i A_{i,I}^0 \to \lim_I \colim_i A_{i,I}^1 \to \lim_I \colim_i A_{i,I}^2 \to \cdots. \]
Its cohomology therefore fits in a Milnor sequence\footnote{According to \cite[Prop.~3.5.8]{weibel}, this is true whenever the system $A_{i, *}^j \to A_{i, *}^j$ is Mittag-Leffler for each fixed $j$. But we have already noted that any inclusion $J \subset I$ induces a surjection $A_{i, J}^j \twoheadrightarrow A_{i, I}^j$.}
    \[ 0 \to {\lim_I}^1 H^{j-1} (\colim_i A_{i,I}^*) \to H^j (\lim_I \colim_i A_{i,I}^*) 
    \to \lim_I H^j ( \colim_i A_{i,I}^*) \to 0. \numberthis \label{eqn:Milnor sequence for showing Gamma Gn is E-M} \]
But for each fixed pair $(i,I)$, the complex $A_{i,I}^*$ is split; the same is therefore also true of the colimit, which thus has zero cohomology. Now (\ref{eqn:Milnor sequence for showing Gamma Gn is E-M}) shows that (\ref{eqn:complex for showing Gamma Gn is E-M}) is exact, and so
\begin{align*}
    \pi_s \Gamma (\G^q, \tau_t \mathcal E) = \left \{ 
    \begin{array}{ll}
        \pi_t \Gamma (\G^q, \mathcal E) & s = t \\
        0 & s < t
    \end{array} \right.
\end{align*}
That is, $\Gamma(\G^q, \tau_t \mathcal E) = \Sigma^t \pi_t \Gamma(\G^q, \mathcal E)$ as required.
\end{proof}

\begin{remark} \label{rem:condensed cohomology of profinite abelian groups}
The same proof shows that $H^*_{\mathrm{cond}}(T, \operatorname{Cont}(-, M)) = 0$ in positive degrees, for any profinite set $S$ and profinite abelian group $M$ with a presentation as a directed colimit $M = \lim_{\N \op} M_n$ satisfying the Mittag-Leffler condition. 
\end{remark}

\begin{remark} \label{rem:Morava E theory is Kan extended}
There are alternative ways to obtain a pro\'etale or condensed object from $\E$, and it is sometimes useful to make use of these. Since we have proven that $\mathcal E = \nu^p \mathcal E^\delta$ is Kan extended from the \'etale classifying site, it is uniquely determined as an object of $\Shh[\mathcal Sp_{\K}]{\proet G}$ by the following properties:
\begin{enumerate}
    \item it has underlying $\K$-local spectrum $\Gamma(\G/*, \mathcal E) = \E$.
    \item its homotopy sheaves are $\Cont_\G(-, \pi_t \E)$, where $\pi_t \E$ is viewed as a topological $\G$-module with its adic topology.
\end{enumerate}
\end{remark}


%% file: Sections/Modules.tex
A similar strategy to that of Proposition \ref{lem:nu^p mathcal E is a sheaf} fails for the sheaf of module categories, since it is not clear that filtered colimits in $\Prsm$ are exact; moreover, we would need an analogue of \cite[Theorem~4.16]{cm}. This section is devoted to the following result, which builds on the descent results of \cite{ms}:

\begin{thm} \label{lem:hyperdescent for presheaf extension of Modk mathcal E}
    The presheaf $\nu^p \Modk{\mathcal E^\delta}: \proet{\G}^{op} \to \Prsm$ satisfies hyperdescent.
\end{thm}

\begin{remark}
By \cref{lem: colim of modules is modules over colim}, $\nu^p \Modk{\mathcal E^{\delta}} = \Modk{\nu^p \mathcal E^{\delta}} = \Modk{\mathcal E}$. Thus taking endomorphisms of the unit gives an alternative proof of \cref{lem:nu^p mathcal E is a sheaf} as an immediate corollary.
\end{remark}

We will deduce \cref{lem:hyperdescent for presheaf extension of Modk mathcal E} from a more general result, true for any descendable Galois extension in a sufficiently nice stable homotopy theory:

\begin{thm} \label{lem:hyperdescent for presheaf extension of Modk mathcal A}
    Let $\mathcal C$ be a compactly assembled stable homotopy theory, and $G$ a profinite group. Suppose that $\bm 1 \to A$ is a descendable $G$-Galois extension in $\mathcal C$, corresponding to $\mathcal A^\delta \in \Sh[\mathcal C]{\et G}$. Then the presheaf
        \[ \nu^p \Mod_{\mathcal A^\delta}(\mathcal C) \in \mathcal P(\proet G, \Prsm) \]
    satisfies hyperdescent.
\end{thm}

Before giving the proof, let us make a few remarks.

\begin{remark}
\begin{enumerate}
    \item By \cite[Proposition~6.15]{mathew_galois}, descendable Galois extensions are faithful.
    
    \item Recall \cite[§21.1.2]{sag} that a presentable \inftycat{} is said to be \emph{compactly assembled} if it is a retract in $\Pr^L$ of a compactly generated presentable \inftycat{}; if $\mathcal C$ is stable, this is equivalent to being dualisable in the symmetric monoidal structure on $\Pr^L_{\mathrm{st}}$.

    \item For the proof, we do \emph{not} require that $\mathcal C \in \CAlg(\Pr^L_{\omega})$. In particular, $\mathcal Sp_{\K}$ is compactly generated (though not by the unit) and so is an example of such an \inftycat{}, and \cref{lem:hyperdescent for presheaf extension of Modk mathcal E} follows since $\bm 1_{\K} \to \E$ is descendable.
\end{enumerate}
\end{remark}

The functor $\pic$ preserves limits of symmetric monoidal \inftycats{} by \cite[Proposition~2.2.3]{ms}, and so we obtain the first part of our main result:

\begin{cor}
There is a hypercomplete sheaf $\pic(\mathcal E) \coloneqq \pic \circ \Modk{\mathcal E} \in \Shh[\mathcal Sp_{\geq 0}]{\proet \G}$ having
    \[ \Gamma(\G, \pic(\mathcal E)) \simeq \pic(\Modk \E) \qquad \text{ and } \qquad \Gamma(*, \pic (\mathcal E)) \simeq \pic (\mathcal Sp_{\K}) = \pic_h.\]
In particular, we get a conditionally convergent spectral sequence
\begin{equation}
\label{eqn:picard spectral sequence} E_{2}^{s,t} = H^{s}(\proet \G, \pi_t \pic(\mathcal E)) \implies \pi_{t-s} \pic_h.
\end{equation}
\end{cor}

\begin{remark}
    In \cref{sec:picard sheaf}, we will evaluate the homotopy sheaves $\pi_t \pic(\mathcal E)$ and hence identify the $E_2$-page with continuous cohomology $H^s(\G, \pi_t \pic(\E))$.
\end{remark}

We begin the proof of \cref{lem:hyperdescent for presheaf extension of Modk mathcal E} with the following observation, stated for later use in slightly greater generality than is needed for this section.

\begin{lem} \label{lem: basechange of descent}
Let $\mathcal C^\otimes$ be a symmetric monoidal $\infty$-category with geometric realisations, and $1 \leq k \leq \infty$. Suppose that $A \in \mathcal C$ is an $\mathbb E_k$-algebra, and that $B \in \operatorname{CAlg}(\mathcal C)$ is such that
\begin{equation}    \label{eqn:descendability}
    \mathcal C \simeq \lim \left( \cosimp[\Mod_B(\mathcal C)]{\Mod_{B \otimes B}(\mathcal C)} \right).
\end{equation}
Then
    \[ \RMod_A(\mathcal C) \simeq \lim \left( \cosimp[\RMod_{A \otimes B}(\mathcal C)]{\RMod_{A \otimes B \otimes B}(\mathcal C)} \right). \]
\end{lem}

\begin{remark}
We call $B \in \CAlg(\mathcal C)$ a \emph{descent algebra} if \cref{eqn:descendability} holds.
\end{remark}

\begin{proof}
Write $B' \coloneqq A \otimes B$; this is itself an $\mathbb E_k$-algebra. We will verify the hypotheses of the Barr-Beck-Lurie theorem. Namely,
\begin{enumerate}
    \item $(-) \otimes_A B' \simeq (-) \otimes B$ is conservative, by virtue of the equivalence (\ref{eqn:descendability}).
    
    \item $\RMod_A(\mathcal C)$ has limits of $B'$-split cosimplicial objects: given $M^\bullet: \bm \Delta \to \RMod_A(\mathcal C)$ with $M^\bullet \otimes_A B'$ split, we can form the limit $M$ in $\mathcal C$, by the descendability assumption. Since $\RMod_A(\mathcal C) \subset \mathcal C$ is closed under cosifted (in fact, all) limits, $M$ is also a limit in $\RMod_A(\mathcal C)$. This limit is clearly preserved by $(-) \otimes_A B'$. \qedhere
\end{enumerate}
\end{proof}


Equipped with this and Lemma \ref{lem: mathcal E of a pullback}, we can now prove the main theorem of this section. I am grateful to Dustin Clausen for the strategy of the following result, which forms the key complement to the results of \cite{mathew_galois}.


\begin{prop} \label{lem:nu* Mod of free G-sets} 
    Let $\mathcal C$ be a compactly assembled stable homotopy theory, and $G$ a profinite group. Suppose that $\bm 1 \to A$ is a descendable $G$-Galois extension in $\mathcal C$, corresponding to $\mathcal A^\delta \in \Sh[\mathcal C]{\et G}$. Then the restriction
        \[ \nu^p \Mod_{\mathcal A^\delta}(\mathcal C)|_{\Free G} \in \mathcal P(\Free G, {\Pr}^L) \]
    is a hypercomplete sheaf.
\end{prop}

\begin{proof}
We do this in a number of steps.

\renewcommand{\labelenumi}{(\arabic{enumi})}
\renewcommand{\labelenumii}{(\textit{\roman{enumii}})}
\begin{enumerate}[leftmargin=\parindent, itemindent=*]
\item As noted in \cref{sec:Free G}, any $G$-set is covered by a free one and every free $G$-set is split. A consequence is that for any free $G$-set $S$, the functor $S' \mapsto S'/G$ is an equivalence $({\proet G})_{/ S} \simeq \{G\} \times \Profin_{/ (S/G)}$, since any $G$-set over $S$ is itself free. Thus suppose $T = \lim_i T_i$ is a profinite set over $S/G$, and choose a convergent sequence of neighbourhoods $U_j \subset G$ of the identity; applying Lemmas \ref{lem: colim of modules is modules over colim} and \ref{lem: mathcal E of a pullback} we deduce the equivalences
\begin{align*}
    \nu^p \Mod_{\mathcal A^\delta}(T \times G) &= \varinjlim_{i,j} \Mod_{\mathcal A^\delta(T_i \times G/U_j)}\\
    & \simeq \varinjlim_{i,j} \Mod_{\mathcal A^\delta(T_i) \otimes \mathcal A^\delta(G/U_j)} \\
    & \simeq \varinjlim_{i} \Mod_{ \varinjlim_j \mathcal A^\delta(T_i) \otimes \mathcal A^\delta(G/U_j)} \\
    & \simeq \varinjlim_{i} \Mod_{\mathcal A^\delta(T_i) \otimes A} \\
    & \simeq \varinjlim_{i} \Mod_{\prod_{T_i} A} \\
    & \simeq \varinjlim_{i} \prod_{T_i} \Mod_{A}.
\end{align*}

Under the aforementioned equivalence we think of this as a presheaf on $\Profin_{/(S/G)}$: that is, if $T$ is a profinite set over $S/G$, then
    \[ T = \lim_i T_i  \mapsto \varinjlim_i \prod_{T_i} \Mod_ A. \]

\item Since $T_i$ is a finite set, $\nu^p \Mod_{\mathcal A}(T_i \times G) \simeq \prod_{T_i} \Mod_ A \simeq \Sh[\Mod_ A]{T_i}$; we claim that the same formula holds for arbitrary $T$. Writing $T = \lim_i T_i$, it will suffice to prove that the adjunction
    \[q^*: \varinjlim \Sh[\Mod_{A}]{T_i} \leftrightarrows \Sh[\Mod_{A}]{T} : q_*, \]
induced by the adjunctions $(q_i)^* \dashv (q_i)_* $ for each projection $q_i: T \to T_i$, is an equivalence; then the claim follows by passing to a limit of finite $T_i$. In fact, since the adjunction is obtained by tensoring the adjunction
    \[ q^*: \varinjlim \mathcal Sh(T_i) \leftrightarrows \mathcal Sh(T) : q_*, \label{eqn:colimit adjunction of sheaves on a profinite space} \numberthis \]
with $\Mod_{A}$ (combine \cite[Remark~1.3.1.6 and Prop.~1.3.1.7]{sag}), it will suffice to prove that \cref{eqn:colimit adjunction of sheaves on a profinite space} is an equivalence.

Let us assume for notational convenience that the diagram $T_i$ is indexed over a filtered poset $J$, which we may do without loss of generality. Then the claim is a consequence of the fact that the topology on $T$ is generated by subsets $q_i^{-1}(x_i)$, where $q_i \colon T \to T_i$ is a finite quotient.  Indeed, let $\mathcal O \coloneqq \mathrm{Open}(T)$; subsets of the form $q_i^{-1}(x_i)$ form a clopen basis, and we will write $\mathcal B \subset \mathcal O$ for the full subcategory spanned by such. Note that
    \[ \mathcal Sh(T) \simeq \mathcal P_{\Sigma}(\mathcal B), \]
where the right-hand side denotes the full subcategory of presheaves that send binary coproducts to products. On the other hand, an object of $\varinjlim \mathcal Sh(T_i)$ is a Cartesian section of the fibration determined by $i \mapsto \mathcal Sh(T_i)$; abusively, we will denote such an object by $(\mathcal F_i)$, where $\mathcal F_i \in \mathcal Sh(T_i)$, leaving the coherence data implicit. Write also 
    \[ (q_{j,\infty})^* : \mathcal Sh(T_j) \rightleftarrows \varinjlim \mathcal Sh(T_i) : (q_{j, \infty})_* \] 
for the colimit adjunction. By applying Yoneda, one verifies that the map
\begin{align*}
    \varinjlim (q_i)^* \mathcal F_i \to q^* ((\mathcal F_i)),
\end{align*}
obtained by adjunction from the maps $(q_{j,\infty})_* (\eta) \colon \mathcal F_j = (q_{j,\infty})_* ((\mathcal F_i)) \to (q_{j, \infty})_* q_* q^* ((\mathcal F_j)) = (q_j)_* q^* ((\mathcal F_i))$ as $j$ varies, is an equivalence. Likewise, adjunct to $q^* ((q_i)_* \mathcal F) \simeq \varinjlim (q_i)^* (q_i)_* \mathcal F \to \mathcal F$ is an equivalence
    \[ ((q_i)_* \mathcal F) \xrightarrow{\sim} q_* \mathcal F. \]

Restricting to $\mathcal B$ (where no sheafification is required for forming the left adjoint), we will show that the unit and counit of the adjunction are equivalences. For the counit $q^* q_* \mathcal F \to \mathcal F$, this is clear:
\begin{align*}
    \left [q^* q_* \mathcal F \right](q_i^{-1}(x_i)) &\simeq \left[ \varinjlim_{j \geq i} (q_j)^* (q_j)_* \mathcal F \right] (q_{i}^{-1}(x_i)) \\
    &\simeq \varinjlim_{j \geq i} \left[ (q_j)^* (q_j)_* \mathcal F (q_{i}^{-1}(x_i)) \right] \\
    &\simeq \varinjlim_{j \geq i} \left[ (q_j)_* \mathcal F (q_{ij}^{-1}(x_i)) \right] \\
    &\simeq \varinjlim_{j \geq i} \left[ \mathcal F (q_i^{-1}(x_i) \right] \simeq \mathcal F (q_i^{-1}(x_i)).
\end{align*}

For the unit $(\mathcal F_j) \to q_* q^* (\mathcal F_j)$, it will suffice to prove for each $i$ that the canonical map
\begin{equation} \label{eqn:unit for equivalence between sheaves on T and limit of sheaves on Ti}
    \mathcal F_i \to (q_i)_* \varinjlim_{j \geq i} (q_j)^* \mathcal F_j 
\end{equation}
is an equivalence. But
\begin{align*}
    \left[ (q_i)_* \varinjlim_{j \geq i} (q_j)^* \mathcal F_j \right](x_i) &\simeq \left[\varinjlim_{j \geq i} (q_j)^* \mathcal F_j \right](q_i^{-1} (x_i)) \\
    &\simeq \varinjlim_{j \geq i} \left[ (q_j)^* \mathcal F_j (q_i^{-1} x_i) \right] \\
    &\simeq \varinjlim_{j \geq i} \left[ \mathcal F_j (q_{ij}^{-1} x_i) \right] \\
    &\simeq \varinjlim_{j \geq i} \left[ (q_{ij})_* \mathcal F_j (x_i) \right] ,
\end{align*}
with respect to which \cref{eqn:unit for equivalence between sheaves on T and limit of sheaves on Ti} is the structure map for $j = i$. This is an equivalence: since each of the coherence maps
    \[ \mathcal F_i \to (q_{ij})_* \mathcal F_j \]
defining the colimit is an equivalence by definition of $\varinjlim \mathcal Sh(T_i)$, the diagram $j \mapsto (q_{ij})_* \mathcal F_j(x_i)$ factors through the groupoid completion $J_{i/}^{\mathrm{gpd}}$. But since $J_{i/}$ is filtered, both inclusions 
    \[ J_{i/} \hookrightarrow J_{i/}^{\mathrm{gpd}} \hookleftarrow \{i\} \]
are cofinal by \cite[Corollary~4.1.2.6]{htt}.

\item We are left to prove that $T \mapsto \Sh[\Mod_{A}]{T} \in \Prsm$ is a hypercomplete sheaf on $\operatorname{Profin}_{/(S/G)}$. This is precisely the content of \cite[Theorem~0.5]{haine_descent_2022}, noting that
\begin{enumerate}
    \item limits in $\Prsm$ are computed in $\Cat_\infty$; 
    \item $\mathcal C$, and so $\Mod_ A(\mathcal C)$, is compactly assembled;
    \item any profinite set $T$ has homotopy dimension zero by \cite[Theorem~7.2.3.6 and Remark~7.2.3.3]{htt}, so that Postnikov towers in $\mathcal Sh(T)$ converge: $\mathcal Sh(T) \simeq \lim_n \Sh[\mathcal S_{\leq n}]{T}$ by \cite[Theorem~7.2.1.10]{htt}. \qedhere
\end{enumerate}
\end{enumerate}
\end{proof}

The proof of \cref{lem:hyperdescent for presheaf extension of Modk mathcal A} follows by combining the previous proposition with \cite[Proposition~3.22]{mathew_galois}:

\begin{proof}[Proof (Prop. \ref{lem:hyperdescent for presheaf extension of Modk mathcal A})]
Let $T_\bullet \twoheadrightarrow T_{-1} = S$ be a hypercovering in $\proet G$, and form the diagram

\[\begin{tikzcd}[column sep=small]
	\vdots & \vdots & \vdots \\
	{S \times G \times G} & {T_0 \times G \times G} & {T_1 \times G \times G} & \cdots \\
	{S \times G} & {T_0 \times G} & {T_1 \times G} & \cdots \\
	{S } & {T_0} & {T_1} & \cdots
	\arrow[from=4-2, to=4-1]
	\arrow[shift left=1, from=4-3, to=4-2]
	\arrow[shift left=2, from=4-4, to=4-3]
	\arrow[from=3-1, to=4-1]
	\arrow[shift right=1, from=2-1, to=3-1]
	\arrow[shift left=1, from=2-2, to=3-2]
	\arrow[from=3-2, to=4-2]
	\arrow[from=3-3, to=4-3]
	\arrow[shift right=1, from=2-3, to=3-3]
	\arrow[shift right=1, from=2-3, to=2-2]
	\arrow[from=3-2, to=3-1]
	\arrow[shift right=1, from=3-3, to=3-2]
	\arrow[shift left=2, from=3-4, to=3-3]
	\arrow[shift left=2, from=1-1, to=2-1]
	\arrow[shift right=2, from=1-2, to=2-2]
	\arrow[shift left=2, from=1-3, to=2-3]
	\arrow[shift right=2, from=2-4, to=2-3]
	\arrow[shift left=1, from=3-3, to=3-2]
	\arrow[shift left=1, from=2-3, to=2-2]
	\arrow[from=2-2, to=2-1]
	\arrow[shift right=1, from=2-2, to=3-2]
	\arrow[shift left=1, from=2-3, to=3-3]
	\arrow[shift left=2, from=2-4, to=2-3]
	\arrow[from=2-4, to=2-3]
	\arrow[shift right=2, from=3-4, to=3-3]
	\arrow[from=3-4, to=3-3]
	\arrow[shift right=2, from=4-4, to=4-3]
	\arrow[from=4-4, to=4-3]
	\arrow[shift right=2, from=1-3, to=2-3]
	\arrow[from=1-3, to=2-3]
	\arrow[shift left=2, from=1-2, to=2-2]
	\arrow[from=1-2, to=2-2]
	\arrow[shift right=2, from=1-1, to=2-1]
	\arrow[from=1-1, to=2-1]
	\arrow[shift left=1, from=2-1, to=3-1]
	\arrow[shift right=1, from=4-3, to=4-2]
\end{tikzcd}\]

To prove descent for the hypercover, it is enough to show descent for each column and for each non-negative row. For each such row we are in the context of \cref{lem:nu* Mod of free G-sets}, and so obtain a limit diagram after applying $\nu^p \Mod_{\mathcal A^\delta}$. On the other hand, writing $T_j = \lim T_{ij}$ and applying $\nu^p \Mod_{\mathcal A^\delta}$ to a column we obtain the complex
    \[ \augcosimp[\varinjlim_i \Mod_{\mathcal A^\delta(T_{ij})}]{\varinjlim_{i, k} \Mod_{\mathcal A^\delta(T_{ij} \times G/U_k)}}{\varinjlim_{i,k } \Mod_{\mathcal A^\delta(T_{ij} \times G/U_k \times G/U_k)}} \]

Using Lemma \ref{lem: colim of modules is modules over colim} and the equivalences
    \[{\varinjlim}_{i,k} \mathcal A^\delta(T_{ij} \times (G/U_k)^n) \simeq \mathcal A(T_j) \otimes A^{\otimes n} \]
one identifies this with the complex
    \[ \augcosimp[\Mod_{\mathcal A(T_j)}]{\Mod_{\mathcal A(T_j) \otimes A}}{\Mod_{\mathcal A(T_j) \otimes A \otimes A}}. \numberthis \label{eqn:complex of module categories for nu^p Modk}\]
    
When $T_j = *$, this is the complex
    \[ \augcosimp[\mathcal C]{\Mod_{A}}{ \Mod_{A \otimes A}} \]
which is a limit diagram according to \cite{mathew_galois}.
For general $T$, (\ref{eqn:complex of module categories for nu^p Modk}) is a limit diagram by combining the $T = *$ case with Lemma \ref{lem: basechange of descent}.
\end{proof}

%% file: Sections/Picard_sheaf.tex
In \cref{sec:modules} we showed that $\K$-local modules determine a sheaf of symmetric-monoidal \inftycats{} on the site $\proet \G$. Since limits in $\Prsm$ are computed in $\Cat_\infty^{\mathrm{smon}}$, the functor
    \[ \pic: \Prsm \to \mathcal Sp_{\geq 0}, \]
preserves them \cite[Prop.~2.2.3]{ms}, and so the composite $\pic(\mathcal E) = \pic \circ \Modk{\mathcal E}$ is immediately seen to be a sheaf of connective spectra. As a result, the $0$-stem in its descent spectral sequence \cref{eqn:picard spectral sequence} converges conditionally to $\Pic_h$. Its $E_2$-page consists of cohomology of the homotopy sheaves $\pi_* \pic(\mathcal E)$, and as in \cref{lem:E1 page of descent spectral sequence} we'd like to identify this with group cohomology with coefficients in the continuous $\G$-module $\pi_* \pic(\E)$. In order to deduce this once again from \cite[Lemma~4.3.9]{bs}, we need to show that
\begin{align} \label{eqn:Picard sheaf homotopy WTS}
    \pi_t \pic (\mathcal E) = \Cont_\G(-, \pi_t \pic(\E)).
\end{align}
for $t = 0, 1$; for $t \geq 2$ the result follows from \cref{lem:pi* nu* mathcal E of free G-sets} and the isomorphism $\pi_t \pic(A) \simeq \pi_{t-1} A$, natural in the ring spectrum $A$. The first aim of this section is to prove \cref{eqn:Picard sheaf homotopy WTS}. Having done this, we will evaluate the resulting spectral sequence.

\begin{thm} \label{lem:pi_* pic mathcal E}
The homotopy sheaves of the Picard sheaf are given by
\begin{align} \label{eqn:pi_* pic mathcal E}
    \pi_t \pic(\mathcal E) = \mathrm{Cont}_\G \left(-, \pi_t \pic(\E) \right) = \left\{ \begin{array}{ll}
         \mathrm{Cont}_{\G} \left(-, \Pic(\E) \right) & t = 0 \\
         \mathrm{Cont}_{\G} \left(-, (\pi_0 \E)^\times \right) & t = 1 \\
         \mathrm{Cont}_\G \left(-, \pi_{t-1} \E \right) & t \geq 2 
    \end{array} \right.
\end{align}
\end{thm}

Before proving this, we give the desired identification of the $E_2$-page in \cref{eqn:picard spectral sequence}:

\begin{cor} \label{lem:E1 page of Picard spectral sequence}
The starting page of the descent spectral sequence for the Picard spectrum $\pic(\mathcal E)$ is given by continuous group cohomology:
    \begin{equation} \label{eqn:Picard spectral sequence corrected}
    E_{2}^{s,t} = H^s(\G, \pi_t \pic (\E)).
\end{equation}
\end{cor}

We now focus on the proof \cref{lem:pi_* pic mathcal E}.

\begin{proof}
As in \cref{lem:pi* nu* mathcal E of free G-sets} we appeal to \cite[Lemma~3.9.4(4)]{bs}. We have already noted that the requisite conditions hold for the sheaves $\pi_t \pic(\mathcal E) \simeq \pi_{t-1} \mathcal E$ when $t \geq 2$; all that remains to justify is what happens at $t = 0$ and $1$. But $\pi_0 \pic (\mathcal E) \simeq \Cont_\G(-,\Z/2)$ certainly satisfies the hypotheses of \cite[Lemma~4.3.9]{bs}, while $\pi_1 \pic (\mathcal E) = \Cont_\G(-, (\pi_0 \E)^\times)$, and $(\pi_0 \E)^\times$ is the limit of finite $\G$-modules $\left(\pi_0 \E / I \right)^\times$.
\end{proof}



\begin{remark} \label{remark:K locally invertible mathcal E modules for free G sets}
As in \cref{sec:Morava E-theory}, it will suffice to prove that the homotopy sheaves take this form on $\Free \G$. If $T$ is a finite set, then $\mathcal E(T \times \G) \simeq \bigoplus_T \E$, and so
    \[ \Cont(T, \Pic(\E)) \cong \bigoplus_T \Pic(\E) \cong \Pic(\Modk{\mathcal E(T \times \G)}). \]
If $T$ is an arbitrary profinite set, the isomorphisms above induce a canonical map
\begin{align*} \label{eqn:comparison map between Cont Pic and Pic Cont}
    \chi: \Cont(T, \Pic(\E)) \to \Pic(\Modk{\mathcal E(T \times \G)}),
\end{align*}
which may be described explicitly as follows: any continuous map $f \colon T \to \Pic(\E) = \Z/2 \{\Sigma \E\}$ defines a clopen decomposition $T = T^0 \sqcup T^1$, with $T^i = f^{-1}(\Sigma^i \E)$. Projecting to finite quotients gives $T_i = T^0_i \sqcup T^1_i$ for $i$ sufficiently large, and since $\mathcal E(T \times \G) = L_{\K} \varinjlim_i \bigoplus_{T_i} \E$ we set
    \[ \chi(f) = L_{\K} \varinjlim_i \left(\bigoplus_{T^0_i} \E \oplus \bigoplus_{T^1_i} \Sigma \E \right) \in \Pic(\Modk{\mathcal E(T \times \G)}). \numberthis \label{eqn:formula for chi} \]
\end{remark}

We will deduce \cref{lem:pi_* pic mathcal E} by showing that the map $\chi$ is an isomorphism. The key point will be the following:

\begin{prop} \label{lem:invertible modules are locally free}
Let $T$ be a profinite set and $X \in \Pic(\Modk{\mathcal E(T \times \G)})$. Then $X$ is in the image of $\chi$.
\end{prop}

Given this, it is straightforward to prove the main result of the section:

\begin{proof}[Proof (\cref{lem:pi_* pic mathcal E}).]
For $t \geq 1$, this follows from the equivalence
    \[ \Omega \Picc \left(\Modk{\mathcal E} \right) \simeq \GL(\mathcal E). \]
and the identification of the homotopy sheaves $\pi_t \mathcal E$ in \cref{lem:pi* nu* mathcal E of free G-sets}.

For $t = 0$, \cref{lem:invertible modules are locally free} implies that $\chi$ is surjective, while injectivity is clear from \cref{eqn:formula for chi}. This yields isomorphisms of sheaves on $\Free \G$,
    \[ \Cont_\G(-, \Pic(\E)) \simeq \Cont((-)/\G, \Pic(\E)) \simeq \pi_0 \pic(\mathcal E).  \qedhere \]
\end{proof}

To prove \cref{lem:invertible modules are locally free} it will be convenient to work in the context of sheaves of $\E$-modules, using the equivalence
    \[ \Modk{\mathcal E(T \times \G)} \simeq \Sh[\Modk \E]{T} \]
from the proof of \cref{lem:hyperdescent for presheaf extension of Modk mathcal A}. We begin by recording some basic lemmas.

\begin{lem}[\cite{stacks-project}, \href{https://stacks.math.columbia.edu/tag/0081}{Tag 0081}] \label{lem:constant sheaves on a profinite set}
Let $T$ be a topological space, and $A$ a set. The constant sheaf $A_T$ on the $T$ takes the form
\begin{equation} \label{eqn:constant sheaves on a profinite set}
    U \mapsto \Loccon{U, A},
\end{equation}
that is, locally constant functions $U \to A$. \qedhere
\end{lem}

\begin{lem} \label{lem:homotopy classes of maps E_T to mathcal F}
For any $X \in \Sh[\Modk{\E}]{T}$,
    \[ \left[ \E_T, X \right] \simeq \Hom(\pi_* \E_T, \pi_* X) \]
\end{lem}

\begin{proof}
We have isomorphisms
\begin{align*}
    \left[ \E_T, X \right] &\simeq \pi_0 \Map(\E_T, X) \\
    &\simeq \pi_0 \Gamma X \\
    &\simeq \Gamma \pi_0 X \\
    &\simeq \Hom(\pi_* \E_T, \pi_* X),
\end{align*}
because the descent spectral sequence for $\Gamma X$ collapses immediately to the $0$-line. Indeed, profinite sets have homotopy dimension zero, and therefore cohomological dimension zero \cite[Corollary~7.2.2.30]{htt}.
\end{proof}

\begin{remark}
After sheafification, postcomposition with the functor $\K^{\E}_*(-) \coloneqq \pi_* (\K \otimes_{\E} -) \colon \Modk \E \to \Mod_{\pi_* \K}$ defines a functor
    \[ \K^\E_* : \Sh[\Modk \E]{T} \to \Sh[\Mod_{\pi_* \K}]{T}. \]
Similarly to \cite{hopkins-mahowald-sadofsky,baker-richter}, our strategy will be to use monoidality of this functor to deduce results about invertible objects in $\Sh[\Modk \E]{T}$ by first proving them at the level of $\K_*^{\E}(-)$.
\end{remark}

\begin{lem} \label{lem: K* of invertible sheaf of E modules is invertible}
Let $T$ be a profinite set, and $X \in \Pic(\Sh[\Modk \E]{T})$. Then $\K^{\E}_* X \in \Pic(\Sh[\Mod_{\pi_* \K}]{T})$.
\end{lem}

\begin{proof}
The functor $\K^{\E}_*(-)$ on spectra is (strict) monoidal: this is \cite[Corollary~33]{baker-richter}. It follows that the induced functor
    \[ \K^\E_* : \Sh[\Modk \E]{T} \to \Sh[\Mod_{\pi_*\K}]{T} \]
is also strictly monoidal, and hence preserves invertible objects.
\end{proof}

\begin{lem} \label{lem:invertible sheaves of k modules are locally free}
Let $T$ be a profinite set and $k$ a field graded by an abelian group $A$. Then any $X \in \Pic(\Sh[\Mod_k]{T})$ is locally free of rank one.
\end{lem}

\begin{proof}
Since $\Sh[\Mod_k]{T} = \Mod_{k_T}$ for $k_T$ the constant sheaf, the result follows from \cite[\href{https://stacks.math.columbia.edu/tag/0B8M}{Tag 0B8M}]{stacks-project} in the ungraded case. In the graded case, \cite[\href{https://stacks.math.columbia.edu/tag/0B8K}{Tag 0B8K}]{stacks-project} still shows that $X$ is locally a summand of a finite free module (though not necessarily one in degree zero); since $k$ is a field it follows that $X$ is locally a shift of $k_T$ by some $a \in A$.
\end{proof}

\begin{proof}[Proof (\cref{lem:invertible modules are locally free}).]
Note that under the equivalence
    \[ \Modk{\mathcal E(T \times \G)} \simeq \Sh[\Modk \E]{T}, \]
the image of $\chi$ in $\Sh[\Modk \E]{T}$ consists of those invertible sheaves that are locally free of rank one: indeed, $T$ has a basis of finite clopen covers $T = \bigsqcup_{x \in T_i} U_x$, and $\chi(\Cont(T_i, \Pic(\E)))$ is the subgroup of invertible sheaves that are constant along $\bigsqcup_{x \in T_i} U_x$.

We will deduce that any $X \in \Pic(\Sh[\Modk \E]{T})$ is locally free of rank one by proving in turn each of the following statements:
\begin{enumerate}
    \item $\K_*^{\E} X \in \Sh[\Mod_{\pi_* \K}]{T}$ is locally free of rank one,
    \item for every $i_0, \dots, i_{n-1} \geq 1$, the sheaf $\pi_*(X/(p^{i_0}, \dots, u_{n-1}^{i_{n-1}})) \in \Sh[\Mod_{\pi_* \E}]{T}$ is locally constant with value $\Sigma^\varepsilon \pi_* \E / (p^{i_0}, \dots, u_{n-1}^{i_{n-1}})$ \footnote{Here $\varepsilon \in \{0, 1\}$ may vary on $T$.},
    \item $X \in \Sh[\Modk \E]{T}$ is locally free of rank one.
\end{enumerate}

Essentially, this follows the proof of \cite[Theorem~3.7]{baker-richter}.

\begin{enumerate}
    \item This follows immediately by combining \cref{lem:invertible sheaves of k modules are locally free,lem: K* of invertible sheaf of E modules is invertible}. Since all our claims are local, we will assume for simplicity that $\K^\E_* X \simeq \pi_* \K_T$.

    \item For $i_0, \dots, i_{n-1} \geq 1$, we first show that  that $\pi_* (X/(p^{i_0}, \dots, u_{n-1}^{i_{n-1}}))$ admits a surjection by $\pi_* \E_T$. This follows by induction on $\sum_j i_j$, with the base case being $(i)$. For the induction step, note that $\Sh[\Modk \E]{T}$ is tensored over $\Modk \E$, and so the cofibre sequences in \cite[Lemma~34]{baker-richter} give rise to cofibre sequences
    \begin{align} \label{eqn:fibre sequence for multiplication by u_j nakayama}
        X/(p^{i_0}, \dots, u_{n-1}^{i_{n-1}}) \xrightarrow{u_j} & X/(p^{i_0}, \dots, u_{n-1}^{i_{n-1}}) \to X/(p^{i_0}, \dots, u_j, \dots u_{n-1}^{i_{n-1}}) \oplus \Sigma X/(p^{i_0}, \dots, u_j, \dots u_{n-1}^{i_{n-1}})
    \end{align}
    and
    \begin{align} \label{eqn:fibre sequence for multiplication by u_j evenness}
        X/(p^{i_0}, \dots, u_j^{i_j - 1}, \dots, u_{n-1}^{i_{n-1}}) \xrightarrow{u_j} & X/(p^{i_0}, \dots, u_j^{i_j}, \dots, u_{n-1}^{i_{n-1}}) \to X/(p^{i_0}, \dots, u_j, \dots u_{n-1}^{i_{n-1}})
    \end{align}
    for $i_j  > 1$; in particular, using \cref{eqn:fibre sequence for multiplication by u_j evenness} and the inductive hypothesis, we can assume that $\pi_* (X/(p^{i_0}, \dots, u_{n-1}^{i_{n-1}}))$ is concentrated in even degrees. Thus \cref{eqn:fibre sequence for multiplication by u_j nakayama} implies we have an exact sequence
        \[ 0 \to u_j\pi_{2t}(X/(p^{i_0}, \dots, u_j^{i_j}, \dots, u_{n-1}^{i_{n-1}})) \xrightarrow{u_j} \pi_{2t}(X/(p^{i_0}, \dots, u_j^{i_j}, \dots, u_{n-1}^{i_{n-1}})) \to \pi_{2t}(X/(p^{i_0}, \dots, u_j, \dots, u_{n-1}^{i_{n-1}})) \to 0 \numberthis \label{eqn:exact sequence for multiplication by u_j nakayama} \]
    for any $t$. Since $\pi_* \E_T$ is free, we can choose a lift $\alpha$ below:
    \[\begin{tikzcd}
	& {\pi_*(X/(p^{i_0}, \dots, u_j^{i_j}, \dots, u_{n-1}^{i_{n-1}}))} \\
	{\pi_* \E_T} & {\pi_*(X/(p^{i_0}, \dots, u_j, \dots, u_{n-1}^{i_{n-1}}))}
	\arrow[from=1-2, to=2-2]
	\arrow[two heads, from=2-1, to=2-2]
	\arrow["\alpha", dashed, from=2-1, to=1-2]
    \end{tikzcd}\]
    Now \cref{eqn:exact sequence for multiplication by u_j nakayama} implies that $\alpha$ is a surjection, since Nakayama's lemma implies that it is so on stalks. As in \cite[Lemma~36]{baker-richter}, this implies $(ii)$ by induction on $\sum i_j$.

    \item Since $X$ is invertible (and in particular dualisable),
        \[ {\lim}_{i_0, \dots, i_{n-1}} (\E/(p^{i_0}, \dots, u_{n-1}^{i_{n-1}}) \otimes_{\E} X) \simeq {\lim}_{i_0 \dots, i_{n-1}} (\E/(p^{i_0}, \dots, u_{n-1}^{i_{n-1}})) \otimes_{\E} X \simeq X. \]
    Now $(ii)$ implies that all $\lim^1$-terms vanish and that
    \begin{align*}
        \pi_* X &\cong \pi_* {\lim}_{i_0, \dots, i_{n-1}} (\E/(p^{i_0}, \dots, u_{n-1}^{i_{n-1}}) \otimes_{\E} X) \\
        &\cong {\lim}_{i_0, \dots, i_{n-1}} \pi_* \E_T/(p^{i_0}, \dots, u_{n-1}^{i_{n-1}})
        \cong \pi_* \E_T.
    \end{align*}
    By \cref{lem:homotopy classes of maps E_T to mathcal F} we can lift the inverse to a map $\E_T \to X$, which is necessarily an equivalence. \qedhere
\end{enumerate}
\end{proof}

This completes the construction of the Picard sheaf $\pic(\mathcal E)$, and the proof that its homotopy sheaves take the desired form; \cref{lem:E1 page of Picard spectral sequence} follows. We will now study the resulting descent spectral sequence, as we did for the descent spectral sequence of $\mathcal E$ in \cref{sec:Morava E-theory}. The construction of the descent spectral sequence using a Postnikov-style filtration is useful in making the comparison with the Picard spectral sequence (\ref{eqn:Picard spectral sequence corrected}). Namely, we make the following observation:

\begin{lem}
    Suppose that $\mathcal C$ is a site, and $\mathcal D$ a presentable $\infty$-category. Let $X \in \Sh[\mathcal D]{\mathcal C}$, and $p,q \colon \mathcal D \to \mathcal Sp_{\geq 0}$ two limit-preserving functors related by a natural equivalence $\tau_{[a,b]} p \simeq \tau_{[a,b]} q$. Then the descent spectral sequences for $p X$ and $q X'$ satisfy the following:
    \begin{enumerate}
        \item The $E_2$-pages agree in a range: $E_{2,p X}^{s,t} \simeq E_{2, q X}^{s,t}$ if $t \in [a,b]$.
        \item Under the isomorphism induced by $(i)$, we have $d_{r, p X}^{s,t} = d_{r, q X}^{s,t}$ if $2 \leq r \leq b-t + 1$.
    \end{enumerate}
\end{lem}

\begin{proof}
This can be seen directly by comparing the Postnikov towers, since both claims depend only on their $[a,b]$-truncations.
\end{proof}

\begin{remark}
We are not asserting that $x \in E_{2, p X}$ survives to $E_r$ if and only if the corresponding class in $E_{2, q X}$ does; this should really be taken as another assumption on the class $x$.
\end{remark}

In our setting, the equivalence
\begin{align} \label{eqn:equivalence of truncations}
    \tau_{[t, 2t-2]} \pic(A) \simeq \tau_{[t,2t-2]} \Sigma A
\end{align}
of \cite{heard-mathew-stojanoska}, valid for $t \geq 2$ and \textit{functorial} in the ring spectrum $A$, implies immediately the following corollary:

\begin{prop} \label{lem:additive differentials}
Let $A \in \Sh[\operatorname{CAlg_{(p)}}]{\proet \G}$, and consider the two spectral sequences \cref{eqn:descent spectral sequence for nu E} and \cref{eqn:Picard spectral sequence corrected}. Then
\begin{enumerate}
    \item $E_{2}^{s,t} \simeq E_{2, +}^{s,t-1}$ if $t \geq 2$.
    \item The differentials $d_{r}$ and $d_{r, +}$ on $E_{r}^{s,t} \simeq E_{r, +}^{s,t-1}$ agree as long as $r \leq t-1$ (whenever both are defined).
\end{enumerate}
\end{prop}


Finally, we want to prove d\'ecalage for the sheaf $\pic(\mathcal E)$, as we did for the sheaf $\mathcal E$ itself. This will allow us to determine the differential $d_{t}$ on classes in $E_{2}^{s,t}$ too.

\begin{prop} \label{lem:decalage of Picard spectral sequence}
D\'ecalage of the Postnikov filtration induces an isomorphism between the following spectral sequences:
\begin{align*}
    E_2^{s,t} = \pi_{t-s} \Gamma \tau_t \mathcal E = H^{s}(\G, \pi_t \pic(\E)) &\implies \pi_{t-s} \pic_h, \\
    \check E_3^{2s-t,s} = \pi^{s} \pi_t \pic( \E^{\otimes \bullet + 1}) & \implies \pi_s \pic_h.
\end{align*}
The first is the Picard spectral sequence \cref{eqn:Picard spectral sequence corrected}, and the second is the Bousfield-Kan spectral sequence for the cosimplicial spectrum $\Gamma (\G^{\bullet + 1}, \pic(\mathcal E)) = \pic(\Modk{\E^{\otimes \bullet + 1}})$.
\end{prop}

\begin{remark}
    The Bousfield-Kan spectral sequence of \cref{lem:decalage of Picard spectral sequence} is by definition Heard's spectral sequence \cite[Theorem~6.13]{heard}.
\end{remark}

\begin{proof}
As in \cref{lem:decalage of descent spectral sequence}, this follows from \cref{lem:decalage for sheaves of spectra} once we prove the following equivalences:
\begin{align*}
    \Gamma(\G^q, \tau_0 \pic(\mathcal E)) &\simeq \pi_0 \Gamma(\G^q, \pic(\mathcal E)), \\
    \Gamma(\G^q, \Sigma \pi_1 \pic(\mathcal E)) &\simeq \tau_1 \Gamma(\G^q, \pic(\mathcal E)).
\end{align*}
Using \cref{lem:pi_* pic mathcal E}, we compute that
    \[ H^{t-s}(\G^q, \Pic(\E)) \simeq H^{t-s}(\G^q, (\pi_0 \E)^\times) = 0 \]
for $t-s > 0$, since condensed cohomology with profinite coefficients vanishes (\cref{rem:condensed cohomology of profinite abelian groups}).
\end{proof}

Following the method of \cite{ms}, we can now identify the first new differential in the Picard spectral sequence:

\begin{cor} \label{lem:nonlinear differential}
Suppose that $t \geq 2$ and $x \in E_{2}^{t,t}$. We abuse notation to identify $x$ with its image in the starting page of the descent spectral sequence for $\mathcal E$, and assume that both classes survive to the respective $E_{t}$-pages. The first nonadditive differential on $x$ in the Picard spectral sequence is
\begin{equation} \label{eqn:formula for first unstable differential}
    d_{t} x = d_{t}^{ASS} x + x^2
\end{equation}
The ring structure in the right hand side is that of the $\K$-local $\E$-Adams spectral sequence \cref{eqn:descent spectral sequence for nu E}.
\end{cor}

\begin{proof}
The same proof that appears in \cite{ms} goes through: namely, this formula holds for the universal cosimplicial \Einfty-ring having a class in $E_2^{t,t}$ of its Bousfield-Kan spectral sequence, and so for the cosimplicial spectrum $\mathcal E(\G^{\bullet +1}) = \E^{\otimes \bullet + 1}$ too.
\end{proof}


%% file: Sections/n=1.tex
In the previous parts we constructed pro\'etale models for the continuous action of $\G$ on Morava E-theory, its $\K$-local module \inftycat{} and its Picard spectrum. Respectively, these are $\mathcal E$ (constructed in \cref{sec:Morava E-theory}), $\Modk{\mathcal E}$ (in \cref{sec:modules}), and $\pic(\mathcal E)$ (in \cref{sec:picard sheaf}). We also described the resulting spectral sequences. In this section, we compute the Picard spectral sequence in the height one case and use this to give a new proof of the results of \cite{hopkins-mahowald-sadofsky} at all primes. As is common at height one, this splits into two cases: the case of odd primes, and the case $p=2$. In both cases, the strategy is first to compute the descent spectral sequence for $\mathcal E$ (which by \cref{lem:decalage of descent spectral sequence} is the $\K$-local $\E$-Adams spectral sequence) and then to use this to compute the Picard spectral sequence. However, the spectral sequences look somewhat different in the two cases. We start with some generalities, which are true uniformly in $h$ and $p$.

\subsection{Morava modules and the algebraic Picard group}
A productive strategy for studying $\Pic_h$ is to compare it to a certain algebraic variant $\Pic_h^{\mathrm{alg}}$ first defined in \cite[\S7]{hopkins-mahowald-sadofsky}; we recall its definition below. In \cite{pstragowski_picard}, Pstrągowski shows that the algebraic approximation is precise if $p \gg h$, a consequence of the fact that in this case the vanishing line in the Adams spectral sequence occurs at the starting page. While there is always a vanishing line, when $p-1 \mid h$ it appears only at a later page, since $\G$ no longer has finite cohomological dimension mod $p$. In this section we will show that this leaves a possibility for \emph{exotic} Picard elements, and more importantly explain how to identify these in the Picard spectral sequence (\cref{lem:algebraic Picard group}).

To do so, we recall that the \emph{completed E-homology} of a $\K$-local spectrum $X$ is
    \[ \E_*^\vee X \coloneqq \pi_* (\E \otimes X) = \pi_* L_{\K}(\E \wedge X). \]

This is naturally a $\pi_*\E$-module: indeed $\K$-localisation is symmetric monoidal, and therefore sends the $\E$-module $\E \wedge X$ to a module over $L_{\K} \E = \E$. As we discuss below, the abelian group $\E_*^\vee X$ has significantly more structure, and is a crucial tool in understanding the $\K$-local category. It was first studied by Hopkins, Mahowald and Sadofsky in \cite{hopkins-mahowald-sadofsky}, where it is denoted $\mathcal K_{h,*}(-)$; it is almost a homology theory, but fails to preserve infinite coproducts as a result of the failure of $\K$-localisation to be smashing. It is nevertheless an extremely effective invariant for Picard group computations, by virtue of the following theorem:

\begin{thm}[\cite{hopkins-mahowald-sadofsky}, Theorem~1.3] \label{thm:hms}
    $M \in \mathcal Sp_{\K}$ is invertible if and only if $\E_*^\vee M$ is a free $\pi_* \E$-module of rank one.
\end{thm}

In particular, \cref{thm:hms} implies that completed E-homology is not a useful invariant of invertible $\K$-local spectra when we only remember its structure as a $\pi_* \E$-module: since $\pi_0 \E$ is a Noetherian local ring, its Picard group is trivial. To get a more interesting invariant, we should remember the equivariant structure coming from the Morava action on $\E$. That is, if $X$ is a $\K$-local spectrum then $\G$ acts on $\E \otimes X$ by acting on the first factor, and therefore acts on $\E_*^\vee X$. This action makes $\E_*^\vee X$ into a \emph{twisted $\G$-$\pi_* \E$-module}, which means by definition that
\begin{equation*} \label{eqn:twisted module}
    g(a \cdot x) = ga \cdot gx
\end{equation*}
for $x \in \E_*^\vee X$, $a \in \pi_* \E$ and $g \in \G$. We will write $\Mod_{\pi_* \E}^\G$ for the category of twisted $\G$-$\pi_* \E$-modules.

\begin{remark}
For any twisted $\G$-$\pi_*\E$-module $M$, the $\G$-action is continuous for the $I_h$-adic topology: if $g x = y \in M$, then \cref{eqn:twisted module} implies that 
   \[ g(x + ax') = y + ga \cdot gx' \in y + I_h^k M \]
for $a \in I_h^k$ and $x' \in M$, since the action of $\G$ on $\pi_* \E$ fixes the $I_h$-adic filtration; that is, $g^{-1}(y + I_h^k)$ contains the open neighbourhood $x + I_h^{k}M$ (compare \cite[Lemma~5.2]{barthel-heard}).
\end{remark}


\begin{defn}
The \emph{algebraic Picard group} of $\mathcal Sp_{\K}$ is
    \[ \Pic_h^{\mathrm{alg}} \coloneqq \Pic \left(\Mod_{\pi_* \E}^{\G} \right). \]
The \emph{exotic Picard group} of $\mathcal Sp_{\K}$ is defined by the exact sequence of abelian groups
    \[ 1 \to \kappa_h \to \Pic_h \xrightarrow{\E_*^\vee} \Pic_h^{\mathrm{alg}}, \]
whose existence follows from \cref{thm:hms}. Restricting both Picard groups to their subgroups of elements concentrated in even degrees, one can equally obtain $\kappa_h$ as the kernel of the map
    \[ \Pic_h^0 \xrightarrow{\E^\vee_*} \Pic_h^{\mathrm{alg}, 0} \simeq \Pic(\Mod_{\pi_0 \E}^\G). \]
\end{defn}

One of the main theorems of \cite{hopkins-mahowald-sadofsky} is the computation $\kappa_1 \simeq \Z/2$ at the prime 2 (this is Theorem~3.3 therein). We want to show how the descent spectral sequence for $\pic (\mathcal E)$ recovers this computation, and we begin by identifying the algebraic Picard group in the spectral sequence. The aim of this subsection is therefore to prove the following result:

\begin{thm} \label{lem:algebraic Picard group}
    At arbitrary height $h$, the $1$-line of the descent spectral sequence for $\pic(\mathcal E)$ computes the image of $\Pic^0_h$ in $\Pic_h^{\mathrm{alg},0}$. The exotic Picard group $\kappa_h$ is computed by the subgroup in filtration $\geq 2$.
\end{thm}

Proving \cref{lem:algebraic Picard group} will require a short discussion of derived complete modules. Firstly, recall that $\pi_0 \E$ is a regular Noetherian local ring, with maximal ideal $I_h = (p, v_1, \dots, v_h)$. If $R$ is any such (classical) ring and $\mathfrak m$ its maximal ideal, the $\mathfrak m$-adic completion functor $(-)^\wedge_{\mathfrak{m}}$ has left-derived functors $L_i$, defined for example in \cite{greenlees-may_derivedcompletion}; this is in spite of $(-)^\wedge_{\mathfrak{m}}$ not being right-exact. For any $R$-module $M$ the completion map $M \to M^\wedge_{\mathfrak m}$ factors through the zero-th derived functor

    \[ M \xrightarrow{\eta_M} L_0 M \xrightarrow{\epsilon_M} M^\wedge_{\mathfrak m}, \]

and one says that $M$ is \emph{L-complete} or \emph{derived $\mathfrak m$-complete} if $\eta_M$ is an isomorphism. Hovey and Strickland prove the following facts about $L$-completion:
\renewcommand{\labelenumi}{($\arabic{enumi}$)}
\begin{enumerate}
    
    \item \label{item: pro-free modules are projectives} For any $R$, the full subcategory spanned by the $L$-complete modules is a thick abelian subcategory of $\Mod_R^{\heartsuit}$, with enough projective generators \cite[Theorem~A.6 and Corollary~A.12]{hovey-strickland}. This is denoted $\widehat{\mathcal M}$ in \textit{op.\ cit.}, but we will write $\Mod_R^{\heartsuit\mathrm{cpl}}$. The projective objects are precisely those $R$-modules which are \emph{pro-free}, that is $M = L_0 \bigoplus_S R$ for some (possibly infinite) set $S$.
    
    \item \label{item:L_0 is a localisation} The functor $L_0$ from modules to $L$-complete modules is a localisation. In particular, colimits in $L$-complete are computed as
        \[ L_0 \colim M , \]
    where $\colim M$ denotes the colimit at the level of modules. Thus $\Mod_{R}^{\heartsuit\mathrm{cpl}}$ is still generated under colimits by the $L$-completion of the unit, and in particular $\kappa$-presentable where $\kappa$ is chosen so that $L_0 R$ is $\kappa$-compact (note that we may have to take $\kappa > \omega$).

    \item \label{item: E^vee is L-complete}
    Any $\mathfrak{m}$-adically complete module is $L$-complete \cite[Theorem~A.6]{hovey-strickland}. In particular, $\pi_* \E$ is an $L$-complete module over itself.
    
    \item The category $\Mod_R^{\heartsuit\mathrm{cpl}}$ admits a unique symmetric monoidal product making $L_0$ a monoidal functor. This is given by the formula
        \[ M \widehat{\otimes}_{L_0 R} N = L_0 (M \otimes_R N) \]
    for $M$ and $N$ $L$-complete. For $\mathfrak{m}$-adically complete modules, one can also define the $\mathfrak m$-complete module
        \[ (M \otimes_R N)^\wedge_{\mathfrak{m}} \simeq (M \widehat{\otimes}_{L_0 R} N)^\wedge_{\mathfrak{m}}, \]
    but we will have no use for this.
    
    \item \label{item:derived completion of fg modules} Derived completion agrees with ordinary completion on finitely generated modules, and on projective modules: in other words, $\epsilon_M$ is an isomorphism in either of these cases. Moreover, for any $N$ the composition
        \[M \otimes_R L_0 N \to L_0 M \otimes_R L_0 N \to L_0(M \otimes_R N) = L_0 M \widehat{\otimes}_{L_0 R} L_0 N \] 
    is an isomorphism when $M$ is finitely generated \cite[Proposition~A.4]{hovey-strickland}. In particular, if $R$ is itself $L$-complete then finitely generated modules are complete, i.e. $M = L_0 M = M^\wedge_{\mathfrak m}$ and $L_i M = 0$ for $i > 0$.

\end{enumerate}
\renewcommand{\labelenumi}{(\textit{\roman{enumi}})}


If $A$ is an $L$-complete $R$-algebra (not necessarily Noetherian), we can define a category of $A$-modules which are $L$-complete with respect to $\mathfrak{m} \subset R$, $\Mod_A^{\heartsuit\mathrm{cpl}} \coloneqq \Mod_A( \Mod_R^{\heartsuit\mathrm{cpl}})$. 

We now specialise to the case $(R, \mathfrak m) = (\pi_0 \E, I_h)$, and work towards the proof of \cref{lem:algebraic Picard group}. It is shown in \cite[Proposition~8.4]{hovey-strickland} that for any $\K$-local spectrum $X$, the $\pi_0 \E$-module $\E_0^\vee X$ is $L$-complete, and that $\E_*^\vee X$ is finitely generated over $\pi_* \E$ if and only if $X$ is $\K$-locally dualisable \cite[Theorem~8.6]{hovey-strickland}. Our first task is to prove that the presheaf of $1$-categories
    \[ S \mapsto \Mod_{\pi_* \mathcal E(S)}\Lcomplete \]
is a stack on $\Free \G$.

\begin{warning}
We would like to proceed as in \cref{sec:modules}, but there is a small subtlety: namely, to deduce descent from the results of \cite{haine_descent_2022} (as in part (3) of the proof of \cref{lem:nu* Mod of free G-sets}), we would need to show that (the nerve of) $\Mod_{\pi_* \E}^{\heartsuit\mathrm{cpl}}$ is compactly generated; in fact, it would suffice to show it is compactly assembled. This is not clear: for example, Barthel and Frankland observe \cite[Appendix~A]{barthel-frankland} that the unit in $\Mod_{R}^{\heartsuit\mathrm{cpl}}$ (which \emph{is} a generator) cannot be compact. For our purposes it is enough to find \emph{any} (small) set of compact generators: by comparison, the $\K$-local category is compactly generated by the $\K$-localisation of any finite type-$h$ spectrum, even though its unit is not compact.

At this point it will be useful to pass to the \inftycat{} of complete modules over the discrete rings $\pi_0 \mathcal E(S)$, as defined in \cite[\S7.3]{sag} or \cite[\S3.4]{bs}. This \emph{can} be seen to be compactly generated by virtue of local (=Greenlees-May) duality; we will make use of this observation in the proof of \cref{lem:descent for L complete Mod pi_* mathcal E}. Given a discrete commutative ring $R$ complete with respect to a finitely generated maximal ideal $\mathfrak m$, we will view $R$ as a (connective) \Einfty-ring and write $\Mod_{R}^{\mathrm{cpl}} \subset \Mod_{R}$ for the sub-\inftycat{} of complete objects \cite[Definition~7.3.1.1]{sag}, which is a localisation of $\Mod_{R}$. There is a unique symmetric monoidal product on $\Mod_R\complete$ for which the localisation is a monoidal functor, and to avoid confusion we denote this by $\widehat{\otimes}_R^{\mathbb L}$. Moreover, the abelian category $\Mod_{R}^{\heartsuit\mathrm{cpl}}$ of $L$-complete discrete $R$-modules includes as the heart of $\Mod_R^{\mathrm{cpl}}$ for a t-structure constructed in \cite[Proposition~7.3.4.4]{sag}, and the $\infty$-categorical localisation functor $L: \Mod_R \to \Mod_R\complete$ agrees upon restriction with the (total) left derived functor of $L$-completion \cite[Corollary~7.3.7.5]{sag}; in particular, $L_0 \simeq \pi_0 L$.
\end{warning}

\begin{ex}
To give an example which makes the difference between $1$-categories and \inftycats{} apparent, one can consider the colimit along the multiplication-by-$p$ maps
    \[ \Z/p \to \Z/p^2 \to \cdots \]
over $\Z[p]$. The colimit in $\Mod_{\Z[p]}^{\heartsuit\mathrm{cpl}}$ is $L_0 (\Z/p^\infty) = 0$; on the other hand, in the $\infty$-categorical setting one has $L \varinjlim \Z/p = \Sigma L_1 (\Z/p^\infty) = \Sigma \Z[p]$. We do not know if $\Z[p]$ can be written as the filtered colimit of finite $p$-groups in $\Mod_{\Z[p]}\Lcomplete$ in some other way. 

In particular, this example shows that the t-structure on $\Mod_{R}\complete$ is generally \emph{not} compatible with colimits in the sense of \cite[Definition~1.2.2.12]{ha}.
\end{ex}

\begin{lem} \label{lem:L complete modules over pi_* mathcal E of free G sets is sheaves on proet}
For any profinite set $T$, we have an equivalence
\begin{align*} 
    \Mod_{\pi_0 \mathcal E(\G \times T)}\complete &\simeq \Sh[\Mod_{\pi_0 \E}\complete]{T},
    \numberthis \label{eqn:L-complete modules over pi_* mathcal E of a free G-set}
\end{align*}
\end{lem}

\begin{proof}
According to \cite[Corollary~2.5]{hovey_morava_2008}, if $i \mapsto X_i \in \mathcal Sp_{\K}$ is a filtered diagram such that $\E_*^\vee X_i$ is \emph{pro-free} for each $i \in I$, then the natural map
    \[ L_0 \varinjlim_i \E_0^\vee X_i \to \E_0^\vee \varinjlim_i X_i \]
is an isomorphism. In particular this applies to give the middle equivalence in
\begin{align} \label{eqn:colimit of pi_* mathcal E(T_i times G)}
    \pi_0 \mathcal E(\G \times T) \simeq \E_0^\vee \left(\varinjlim_i \bigoplus_{T_i} \bm 1_{\K} \right) \simeq L_0 \varinjlim_i \bigoplus_{T_i} \pi_0 \E = L_0 \varinjlim_i \pi_0 \mathcal E(\G \times T_i),
\end{align}
since $\E_0^\vee \left(\bigoplus_{T_i} \bm 1_{\K} \right) = \bigoplus_{T_i} \E_0^\vee \bm 1_{\K} = \bigoplus_{T_i} \pi_0 \E$ is certainly pro-free, each $T_i$ being finite. As a result, \cref{item:L_0 is a localisation} implies that $\pi_0 \mathcal E(\G \times T)$ is the colimit in $\Mod^{\heartsuit\mathrm{cpl}}_{\pi_0 \E}$ of the algebras $\pi_0 \mathcal E(\G \times T_i) = \bigoplus_{T_i} \pi_0 \E$. In fact this particular example \emph{is} also the limit in $\Mod_{\pi_0 \E}\complete$: indeed, for $s > 0$ we see that $L_s \varinjlim \pi_0 \mathcal E(\G \times T_i) = 0$ by \cite[Theorem~A.2(b)]{hovey-strickland}, since $\varinjlim \bigoplus_{T_i} \pi_0 \E$ is projective in $\Mod_{\pi_0 \E}^\heartsuit$.

This yields the first of the following equivalences of (presentably symmetric monoidal) \inftycats{}:
    \[ \Mod^{\mathrm{cpl}}_{\pi_0 \mathcal E(\G \times T)} \simeq \Mod^{\mathrm{cpl}}_{L \varinjlim \pi_0 \mathcal E(\G \times T_i)} \simeq \varinjlim \Mod^{\mathrm{cpl}}_{\pi_0 \mathcal E(\G \times T_i)} \simeq \Sh[\Mod^{\mathrm{cpl}}_{\pi_0 \E}]{T}. \]
The second equivalence follows from \cref{lem: colim of modules is modules over colim}, and the second can be proved identically to part (2) of the proof of \cref{lem:nu* Mod of free G-sets} (replacing $\mathcal E$ there by $\pi_0 \mathcal E$).
\end{proof}

\begin{prop} \label{lem:descent for L complete Mod pi_* mathcal E}
The presheaf
    \[ S \mapsto \Mod^{\heartsuit\mathrm{cpl}}_{\pi_0 \mathcal E(S)} \]
 it is a stack of $1$-categories; in other words, it is a sheaf of $1$-truncated symmetric monoidal \inftycats{} on $\Free \G$.
\end{prop}

\begin{proof}
We will proceed in a few steps: we will begin by showing that $S \mapsto \Mod_{\pi_0 \mathcal E(S)}^{\mathrm{cpl}}$ is a sheaf of \inftycats{}, and then deduce the desired result at the level of $1$-categories. 

\renewcommand{\labelenumi}{($\arabic{enumi}$)}
\begin{enumerate}[leftmargin = \parindent, itemindent = *]
    \item \label{item:descent for Mod pi_* mathcal e complete} Since any covering in $\Free \G$ is of the form 
    $ p \times \G \colon T' \times \G \to T \times \G $
for $p: T' \to T$ a covering of profinite sets, we can restrict attention to the presheaf 
\begin{equation} \label{eqn:presheaf of L complete modules over pi_* mathcal E}
    T \mapsto \Mod_{\pi_0 \mathcal E(T \times \G)}\complete
\end{equation}
on $\Profin \simeq {\proet \G}_{/\G}$. The local duality equivalence
    \[ \Mod_{\pi_0 \E}^{\mathrm{tors}} \simeq \Mod_{\pi_0 \E}^{\mathrm{cpl}} \]
of \cite[Theorem~3.7]{barthel-heard-valenzuela} or \cite[Proposition~7.3.1.3]{sag} implies that $\Mod_{\pi_0 \E}^{\mathrm{cpl}}$ is compactly generated, and hence dualisable in $\Pr^L$. Thus the presheaf
    \[ T \mapsto \Sh[\Mod_{\pi_0 \E}^{\mathrm{cpl}}]{T} \simeq \mathcal Sh(T) \otimes \Mod_{\pi_0 \E}^{\mathrm{cpl}} \]
is a sheaf on $\Profin$ by the main theorem of \cite{haine_descent_2022}, and so \cref{lem:L complete modules over pi_* mathcal E of free G sets is sheaves on proet} implies that \cref{eqn:presheaf of L complete modules over pi_* mathcal E} is a sheaf too.


\item \label{item:descent for Mod pi_* mathcal e heart complete} Next, we deduce descent for the presheaf of 1-categories $S \mapsto \Mod_{\pi_0 \mathcal E(S)}\Lcomplete$. Given a covering $p: T' \to T$ of profinite sets, we form the diagram
\[\begin{tikzcd}[ampersand replacement=\&,column sep=small]
	{\Mod_{\pi_0 \mathcal E(S)}^{\heartsuit\mathrm{cpl}}} \& \lim \& {\Mod_{\pi_0 \mathcal E({S'})}^{\heartsuit\mathrm{cpl}}} \& {\Mod_{\pi_0 \mathcal E({S'}^{\times_S 2} )}^{\heartsuit\mathrm{cpl}}} \& {\Mod_{\pi_0 \mathcal E({S'}^{\times_S 3} )}^{\heartsuit\mathrm{cpl}}}  \& \cdots \\
	{\Mod_{\pi_0 \mathcal E(S)}^{\mathrm{cpl}}} \& \lim \& {\Mod_{\pi_0 \mathcal E({S'})}^{\mathrm{cpl}}} \& {\Mod_{\pi_0 \mathcal E({S'}^{\times_S 2})}^{\mathrm{cpl}}} \& {\Mod_{\pi_0 \mathcal E({S'}^{\times_S 3} )}^{\mathrm{cpl}}} \& \cdots
	\arrow[hook, from=1-1, to=2-1]
	\arrow[hook, from=1-3, to=2-3]
	\arrow[hook, from=1-4, to=2-4]
	\arrow[hook, from=1-5, to=2-5]
	\arrow[hook, from=1-2, to=2-2]
	\arrow["\sim", from=2-1, to=2-2]
	\arrow[from=2-2, to=2-3]
	\arrow[shift right=1, from=2-3, to=2-4]
	\arrow[shift left=1, from=2-3, to=2-4]
	\arrow[shift left=2, from=2-4, to=2-5]
	\arrow[from=2-4, to=2-5]
	\arrow[shift right=2, from=2-4, to=2-5]
	\arrow[shift left=3, from=2-5, to=2-6]
	\arrow[shift right=3, from=2-5, to=2-6]
	\arrow[shift left=1, from=2-5, to=2-6]
	\arrow[shift right=1, from=2-5, to=2-6]
    \arrow[shift left=3, from=1-5, to=1-6]
	\arrow[shift right=3, from=1-5, to=1-6]
	\arrow[shift left=1, from=1-5, to=1-6]
	\arrow[shift right=1, from=1-5, to=1-6]
	\arrow[shift left=2, from=1-4, to=1-5]
	\arrow[shift right=2, from=1-4, to=1-5]
	\arrow[from=1-4, to=1-5]
	\arrow[shift left=1, from=1-3, to=1-4]
	\arrow[shift right=1, from=1-3, to=1-4]
	\arrow[from=1-2, to=1-3]
	\arrow[from=1-1, to=1-2, "\theta"]
\end{tikzcd}\]
where $S^{(i)} = T^{(i)} \times \G$ as usual. Note that the limit of the top row can be computed as the limit of the truncated diagram $\Delta_{\leq 3} \hookrightarrow \Delta \to \Cat_{\infty}$, since each term is a $1$-category. Moreover, the diagram shows that $\theta$ is fully faithful, and so to prove descent it remains to show that $\theta$ is essentially surjective. That is, given $M \in \Mod_{\pi_0 \mathcal E(S)}\complete$ with $M \widehat{\otimes}_{\pi_0 \mathcal E(S)}^{\mathbb L} \pi_0 \mathcal E(S')$ discrete, we must show that $M$ was discrete to begin with.

To this end, we claim first that $\pi_0 \mathcal E(S')$ is projective over $\pi_0 \mathcal E(S)$; since the graph of $p$ exhibits $\pi_0 \mathcal E(S') = \Cont(T', \pi_0 \E)$ as a retract of $\Cont(T' \times T, \pi_0 \E) = \Cont(T', \Cont(T, \pi_0 \E))$, it will suffice for this part to show that the latter is projective. But $\Cont(T, \pi_0 \E)$ is pro-discrete, which implies that
\begin{align*}
    \Cont(T', \Cont(T, \pi_0 \E)) &\cong \lim_I \varinjlim_i \bigoplus_{T'_i} \Cont(T, \pi_0 \E/I) \\
    &\cong \lim_I \varinjlim_i \bigoplus_{T'_i} \Cont(T, \pi_0 \E) \otimes_{\pi_0 \E} \pi_0 \E/I \\
    &\cong \lim_I \left[ \varinjlim_i \bigoplus_{T'_i} \Cont(T, \pi_0 \E) \right] \otimes_{\pi_0 \E} \pi_0 \E/I \\
    &\cong L_0 \varinjlim_i \bigoplus_{T'_i} \Cont(T, \pi_0 \E)
\end{align*}
is pro-free. To obtain the final isomorphism, we've used the fact that each term in the colimit is free over $\Cont(T, \pi_0 \E)$, so that the (uncompleted) colimit is projective. As a result, for any complete $\pi_0 \mathcal E(S)$-module spectrum $M$ we have
    \[ \pi_* \left(M \widehat{\otimes}_{\pi_0 \mathcal E(S)}^{\mathbb L} \pi_0 \mathcal E(S') \right) = (\pi_* M) \widehat{\otimes}_{\pi_0 \mathcal E(S)} \pi_0 \mathcal E(S'). \]
Since $\pi_0 \mathcal E(S')$ is faithful over $\pi_0 \mathcal E(S)$, we deduce that $M$ is discrete whenever its basechange is. \qedhere
\end{enumerate}
\end{proof}
\renewcommand{\labelenumi}{($\roman{enumi}$)}

\begin{remark}
For any graded ring $\bigoplus_{a \in A} R_a$ with $R_0$ a complete noetherian local ring, invertible objects in $\Mod_{R}\Lcomplete$ are locally free: this follows from \cref{lem:invertible sheaves of k modules are locally free} and Nakayama's lemma.
\end{remark}

 It is convenient at this stage to work with Picard spaces, which as usual we denote $\Picc$.

\begin{cor} \label{lem:descent for Pic of L complete Mod pi_* mathcal E}
The presheaf
    \[ S \mapsto \Picc \left(\Mod_{\pi_* \mathcal E(S)}\Lcomplete \right) \]
is a sheaf of groupoids on $\Free \G$.
\end{cor}

\begin{proof}
By \cref{lem:descent for L complete Mod pi_* mathcal E}, the assignment
\begin{align} \label{eqn:descent for Pic of L complete Mod pi_0 mathcal E}
    S \mapsto \Picc \left(\Mod_{\pi_0 \mathcal E(S)}\Lcomplete \right)
\end{align}
is a sheaf. Since invertible objects in $\Mod_{\pi_0 \mathcal E(S)}\Lcomplete$ are locally free, this extends to the graded case.
\end{proof}

We are now equipped to prove the promised result, identifying the algebraic elements in the Picard spectral sequence.

\begin{proof}[Proof. (\cref{lem:algebraic Picard group})]
In the 0-stem of the descent spectral sequence, the bottom two lines compute the image of the map $\pi_0 \Gamma \pic(\mathcal E) \to \pi_0 \Gamma \tau_{\leq 1} \pic(\mathcal E)$. We will argue by computing the target, identifying it with the algebraic Picard group.

First recall that by definition, $\pic(\mathcal E) = \pic(\Modk{\mathcal E})$. Observe that $\tau_{\leq 1} \Picc(\Modk{\mathcal E(S)}) = \Picc(h \Modk{\mathcal E(S)})$, and so $ \tau_{\leq 1} \Picc(\mathcal E)$ is the sheafification of
    \[ S \mapsto \Picc \left(h\Modk{\mathcal E(S)} \right). \]
Given $\mathcal E(S)$-modules $M$ and $M'$ with $M \in \Pic(\Modk{\mathcal E(S)})$, we saw in \cref{lem:invertible modules are locally free} that $M$ and so $\pi_* M$ is locally free, and therefore the latter is projective over $\pi_* \mathcal E(S)$. The universal coefficient spectral sequence \cite[Theorem~4.1]{ekmm} over $\mathcal E(S)$ therefore collapses. Thus
\begin{align*}
    [M, M']_{\mathcal E(S)} \simeq \Hom_{\pi_* \mathcal E (S)}(\pi_* M, \pi_* M'),
\end{align*}
and we see that the functor
\begin{align}   
    \pi_*: \Picc \left( h\Modk{\mathcal E(S)} \right) \to \Picc \left( \Mod_{\pi_* \mathcal E(S)}^{\heartsuit\mathrm{cpl}} \right)
\end{align}
is fully faithful. It is in fact an equivalence: any invertible $L$-complete module over $\pi_* \mathcal E(S)$ is locally free, and so projective, and in particular lifts to $\Modk{\mathcal E(S)}$ \cite[Theorem~3]{wolbert}.

By \cref{lem:descent for Pic of L complete Mod pi_* mathcal E}, no sheafification is therefore required when we restrict $\tau_{\leq 1} \Picc\left(\Modk{\mathcal E} \right)$ to $\Free \G$ and so we obtain
\begin{align*}
    \Gamma \tau_{\leq 1} \Picc (\mathcal E) 
    &\simeq \Tot \left[\Picc\left( h\Modk{\mathcal E(\G^{\bullet + 1} )}\right) \right] \\
    &\simeq \Tot \left[ \Picc \left( \Mod_{\Cont(\G^\bullet,\pi_* E)}^{\heartsuit\mathrm{cpl}}\right) \right] \\
    &\simeq \Tot \left[ \Picc\left( \Mod_{\operatorname{Cont}(\G^\bullet, \pi_* E)}\right) \right] \\
    &\simeq \Tot_3 \left[ \Picc\left( \Mod_{\operatorname{Cont}(\G^\bullet, \pi_* E)}\right) \right]
\end{align*}
This is the groupoid classifying twisted $\G$-$\pi_* \E$-modules with invertible underlying module; in particular, $\pi_0 \Gamma \tau_{\leq 1} \Picc(\mathcal E) \simeq \Pic_h^{\mathrm{alg}}$. On free $\G$-sets, the truncation map
    \[ \Picc \left(\Modk{\mathcal E(S)} \right) \to \tau_{\leq 1} \Picc \left( \Modk{\mathcal E(S)}^{\heartsuit\mathrm{cpl}} \right) \simeq \Picc \left(\Mod_{\pi_* \mathcal E(S)}^{\heartsuit\mathrm{cpl}} \right) \]
is just $M \mapsto \pi_* M$. On global sections, it therefore sends $M$ to the homotopy groups of the associated descent datum for the covering $\G \to *$; this is its Morava module $\E_*^\vee M = \pi_* \E \otimes M$, by definition.
\end{proof}

\begin{remark}
As a consequence, the map $\E_*^\vee: \Pic_h \to \Pic_h^{\mathrm{alg}}$ is:
\begin{enumerate}
    \item injective if and only if the zero stem in the $E_\infty$ term of the Picard spectral sequence is concentrated in filtration $\leq 1$;
    \item surjective if and only if there are no differentials in the Picard spectral sequence having source $(0,1)$ (the generator of the $\Z/2$ in bidegree $(0,0)$ certainly survives, and represents $\Sigma \bm 1_{\K}$).
\end{enumerate}
This refines the algebraicity results in \cite{pstragowski_picard}, although in practise it is hard to verify either assertion without assuming a horizontal vanishing line at $E_2$ in the ASS (which is what makes the results of \opcit\ go through).
\end{remark}

\begin{cor}[\cite{culver-zhang}, Proposition~1.25] \label{lem:n^2 = 2p-1}
If $p > 2, p-1 \nmid h$ and $h^2 \leq 4p-4$, there is an isomorphism
    \[ \kappa_h \cong H^{2p-1}(\G, \pi_{2p-2} \E). \]
\end{cor}

\begin{proof}
This follows from sparsity in the Adams spectral sequence. By \cite[Prop.~4.2.1]{heard_thesis}, the lowest-filtration contribution to the exotic Picard group comes from $E_\infty^{2p-1, 2p-1}$. If $p-1 \nmid h$, the vanishing line in the Adams spectral sequence occurs at the starting page, and if moreover $h^2 \leq 4p-4$ then the group $E_\infty^{2p-1, 2p-1}$ is the only possibly nonzero entry in this region, and hence fits in the exact sequence
    \[ 1 \to (\pi_0 \E)^\times \xrightarrow{d_{2p-1}} H^{2p-1}(\G, \pi_{2p-2} \E) \to E_\infty^{2p-1, 2p-1} \to 0. \]

In fact, the differential must vanish. Indeed, a weak form of chromatic vanishing proven in \cite[Lemma~1.33]{bobkova-goerss} shows that
    \[ H^0(\G, \pi_0 \E) \simeq H^0(\G, \W(\F[p^h])) \simeq H^0(\mathrm{Gal}(\F[p^h]/\F[p]), \W(\F[p^h])) \simeq \Z[p]. \]
This isomorphism is inverse to a component of the map on the $E_2$-pages of Adams spectral sequences associated to the diagram 
\[\begin{tikzcd}[ampersand replacement=\&]
	{\bm 1_p} \& {\bm 1_{\K}} \\
	{\mathbb{SW}_h} \& \E
	\arrow[from=1-1, to=1-2]
	\arrow[from=1-2, to=2-2]
	\arrow[from=1-1, to=2-1]
	\arrow[from=2-1, to=2-2]
\end{tikzcd}\]
Here $\bm 1_p$ denotes the $p$-complete sphere, and $\mathbb{SW}_h$ the spherical Witt vectors of $\F[p^h]$; the bottom map is defined under the universal property of $\mathbb{SW}_h$ \cite[Definition~5.2.1]{ellipticii} by the inclusion $\F[p^h] \hookrightarrow (\pi_0 \E)/p$. Since $(\pi_0 \mathbb S \W^{h\mathrm{Gal}})^\times = (\pi_0 \bm 1_p)^\times = \Z[p]^\times$,  the map on Picard spectral sequences induced by the above square implies that the group $E_2^{0,1} \simeq \Z[p]^\times$ in the Picard spectral sequence consists of permanent cycles.
\end{proof}

\begin{ex}
At the prime three this gives $\kappa_3 \cong H^5(\G, \pi_4 \E)$. In this case, the Morava stabiliser group has cohomological dimension nine.
\end{ex}
    
\begin{ex} \label{ex:n^2 = 2p-1}
In the boundary case $2p-1 = \operatorname{cd}_p(\G) = h^2$, we can use Poincar\'e duality to simplify the relevant cohomology group: this gives
\begin{equation}
    \kappa_h \cong H^{2p-1}(\G, \pi_{2p-2} \E) \cong H_0(\mathbb S, \pi_{2p-2} \E)^{\mathrm{Gal}}.
\end{equation}
Examples of such pairs $(h,p)$ are $(3,5)$, $(5,13)$, $(9,41)$ and $(11,61)$; in each case, this is the first prime for which \cite[Remark~2.6]{pstragowski_picard} leaves open the possibility of exotic Picard elements. The case $(h,p) = (3,5)$ case was considered by Culver and Zhang, using different methods; however, they show as above that Heard's spectral sequence combined with the conjectural vanishing
    \[ H_0(\mathbb S, \pi_{2p-2} \E) = 0, \]
would imply that $\kappa_h = 0$ \cite[Corollary~1.27]{culver-zhang}.

To our knowledge, it is not known if there are infinitely primes $p$ for which $2p-1$ is a perfect square; this is closely tied to Landau's (unsolved) fourth problem, which asks if there are infinitely many primes of the form $h^2 + 1$.
\end{ex}

\subsection{Picard groups at height one}
It is well-known that Morava $E$-theory at height one (and a fixed prime $p$) is $p$-completed complex K-theory $KU_p$, and as such its homotopy is given by
\begin{equation} \label{eqn:homotopy of KUp}
    \pi_* \E = \Z[p] [u^{\pm 1}],
\end{equation}

with $u \in \pi_{2} \E$ the Bott element. In this case, the Morava stabiliser group is isomorphic to the $p$-adic units $\Z[p]^\times$, acting on $KU_p$ by Adams operations
    \[ \psi^a: u \mapsto a u. \]

The $\K$-local $\E$-Adams spectral sequence therefore reads
\begin{equation} \label{eqn:ASS at height 1}
    E_{2,+}^{s,t} = H^{s}(\Z[p]^\times, \Z[p](t/2)) \implies \pi_{t-s} L_{\K} S,
\end{equation}
where $\Z[p](t/2)$ denotes the representation
    \[ \Z[p]^\times \xrightarrow{t/2} \Z[p]^\times \to \Z[p] \]
when $t$ is even, and zero when $t$ is odd. Note that these are never discrete $\Z[p]^\times$-modules, except at $t = 0$. Nevertheless, cohomology of continuous pro-$p$ modules is sensible for profinite groups $G$ of type $p$-$\mathrm{FP}_\infty$ \cite[\S4.2]{symonds-weigel}. The \emph{small} Morava stabiliser groups $\mathbb S$ are in general $p$-adic Lie groups, and so satisfy this assumption; this implies that $\G$ is type $p$-$\mathrm{FP}_\infty$, since $\mathbb S$ is a finite index subgroup. In this case cohomology is continuous, in the sense that its value on a pro-$p$ module is determined by its value on finite quotients:
    \[ H^*(G, \varprojlim M_i) = \varprojlim H^*(G, M_i). \]

Under the same assumption on $G$, there is also a Lyndon-Hochschild-Serre spectral sequence for any closed normal subgroup $N <_o G$ \cite[Theorem~4.2.6]{symonds-weigel}: that is,
\begin{equation} \label{eqn:Lyndon-Hochschild-Serre}
    E_2^{i,j} = H^i(G/N, H^{j}(N, M)) \implies H^{i+j}(G, M).
\end{equation}
See also \cite[Theorem~3.3]{jannsen} for a similar result for profinite coefficients that are not necessarily pro-$p$: one obtains a Lyndon-Hochschild-Serre spectral sequence by replacing the Galois covering $X' \to X$ therein by a map of sites $\et G \to \et N$.

The descent spectral sequence for $\pic(\mathcal E)$ at height one (and all primes) therefore has starting page
\begin{align} \label{eqn:Picard SS at height one}
    E_2^{s,t} =  \left\{
    \begin{array}{ll}
         H^{s}(\Z[p]^\times, \Z/2) & t = 0 \\
         H^{s}(\Z[p]^\times, \Z[p](0)^\times) & t = 1 \\
         H^{s}(\Z[p]^\times, \Z[p](\frac{t-1}{2})) & t \geq 2
    \end{array} \right.
\end{align}
The results of \cref{sec:picard sheaf} also tell us how to discern many differentials in the Picard spectral sequence from those in the Adams spectral sequence: in particular, we will make use of \cref{lem:additive differentials} and \cref{lem:nonlinear differential}. Our input is the well-known computation of the $\K$-local $\E$-Adams spectral sequence at height one. A convenient reference is \cite[\S4]{beaudry-goerss-henn_splitting}, but for completeness a different argument is presented in \cref{app:ASS n=1}.

\subsubsection{Odd primes}
\begin{figure}
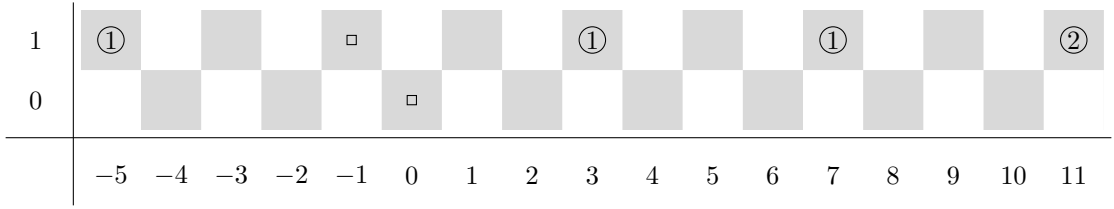

    \centering
    
    \printpage[name = HFPSSodd, grid = chess, page = 7, xrange = {-5}{11}]
    
    \caption[Descent spectral sequence at hight one, odd primes]{The $E_2$-page of the descent spectral sequence for $\E$ at odd primes (implicitly at $p=3$). Squares denote $\Z[p]$-summands, and circles are $p$-torsion summands (labelled by the degree of the torsion). This recovers the well-known computation of $\pi_* \bm{1}_{\K}$ at height 1 and odd primes.}
    \label{fig:HFPSS for n=1 odd primes}
\end{figure}

When $p > 2$, the Adams spectral sequence collapses immediately:

\begin{lem}[\cref{lem:ASS p odd app}]
The starting page of the descent spectral sequence for $\mathcal E$ is given by
\begin{equation}
    E_{2,+}^{s,t} = H^{s}(\Z[p]^\times, \pi_{t} \E) = \left\{
        \begin{array}{ll}
             \Z[p] & t = 0 \text{ and } s = 0, 1 \\
             \Z/p^{\nu_p(t) + 1} & t = 2(p-1)t' \neq 0 \text{ and } s = 1
        \end{array} \right.
\end{equation}
and zero otherwise. The result is displayed in \cref{fig:HFPSS for n=1 odd primes}.
\end{lem}

\begin{sseqdata}[Adams grading, scale = .8, name = HFPSSodd]
        \class[rectangle](0,0)
        \class[rectangle](-1,1)
        \class[minimum width = width("2") + 0.5em,"2"](11,1)
        \class[minimum width = width("2") + 0.5em,"1"](-5,1)
        \class[minimum width = width("2") + 0.5em,"1"](3,1)
        \class[minimum width = width("2") + 0.5em,"1"](7,1)
\end{sseqdata}

As a result of the vanishing line, the computation of $\Pic(\mathcal Sp_{\K})$ in this case depends \textit{only} on $H^*(\Z[p]^\times, \Pic(\E))$ and $H^*(\Z[p]^\times, (\pi_0 \E)^\times)$. This recovers the computation in \cite[Proposition~2.7]{hopkins-mahowald-sadofsky}.

\begin{prop} \label{lem:Pic_1 p odd}
The height one Picard group is algebraic at odd primes:
\begin{equation}
    \Pic_1 \cong \Pic_1^{\mathrm{alg}} \cong \Z[p] \times \Z/2(p-1)
\end{equation}
\end{prop}




\subsubsection{The case \texorpdfstring{$p = 2$}{p = 2}}
At the even prime, the Morava stabiliser group contains $2$-torsion, and therefore its cohomology with $2$-complete coefficients is periodic.

\begin{lem} \label{lem:ASS p odd}
The starting page of the descent spectral sequence for $\mathcal E$ is given by
\begin{equation}
    E_{2,+}^{s,t} = H^{s}(\Z[2]^\times, \pi_{t} \E) = \left\{
        \begin{array}{ll}
             \Z[2] & t = 0 \text{ and } s = 0,1 \\
             \Z/2 & t \equiv_4 2 \text{ and } s \geq 1 \\
             \Z/2^{\nu_2(t)+1} & 0 \neq t \equiv_4 0 \text{ and } s = 1 \\
             \Z/2 & t \equiv_4 0 \text{ and } s > 1
        \end{array} \right.
\end{equation}
and zero otherwise.
\end{lem}

This time we see that the spectral sequence can support many differentials. These can be computed by various methods, as for example in \cite{beaudry-goerss-henn_splitting}. We give another proof, more closely related to our methods, in \cref{app:ASS n=1}.

\begin{prop} \label{lem:differentials in ASS p=2}
The descent spectral sequence collapses at $E_4$ with a horizontal vanishing line. The differentials on the third page are displayed in \cref{fig:HFPSS for n=1 even primes}.
\end{prop}

As a result of the previous subsection, we can compute the groups of exotic and algebraic Picard elements at the prime 2. We will need one piece of the multiplicative structure: write $\eta$ for the generator in bidegree $(s,t) = (1,2)$, and $u^{-2} \eta^2$ for the generator in bidegree $(s,t) = (2,0)$.

\begin{lem}
In the descent spectral sequence for $\mathcal E$, the class
    \[ x \coloneqq u^{-2} \eta^2 \cdot \eta \in E_2^{3, 2} \]
is non-nilpotent. In particular, $x^j$ generates the group in bidegree $(s,t) = (3j,2j)$ of the descent spectral sequence for $\mathcal E$.
\end{lem}

\begin{proof}
The classes $u^{-2} \eta^2$ and $\eta$ are detected by elements of the same name in the HFPSS for the conjugation action on $KU_2$, under the map of spectral sequences induced by the square of Galois extensions
\[\begin{tikzcd}[ampersand replacement=\&]
	{\bm 1_{\K}} \& {KO_2} \\
	{KU_2} \& {KU_2}
	\arrow[from=1-1, to=2-1]
	\arrow[from=1-1, to=1-2]
	\arrow[from=1-2, to=2-2]
	\arrow[no head, equals, from=2-1, to=2-2]
\end{tikzcd}\]
To see this, one can trace through the computations of \cref{app:ASS n=1}: indeed, the proof of \cref{lem:differentials in ASS p=2 app} identifies the map of spectral sequences induced by the map of $C_2$-Galois extensions
\[\begin{tikzcd}[ampersand replacement=\&]
	{\bm 1_{\K}} \& {KO_2} \\
	{KU_2^{h(1+4\Z[2])}} \& {KU_2}
	\arrow[from=1-1, to=2-1]
	\arrow[from=1-1, to=1-2]
	\arrow[from=1-2, to=2-2]
	\arrow[from=2-1, to=2-2]
\end{tikzcd}\]
and \cref{lem:isomorphism between two spectral sequences n=1 p=2} identifies the descent spectral sequence for $\mathcal E$ with the HFPSS for $KU_2^{h(1+4\Z[2])}$, up to a filtration shift. But the starting page of the HFPSS for $KU_2$ is
    \[ H^*(C_2, \pi_* KU_2) = \Z[2][\eta, u^{\pm 2}]/2\eta, \]
and here in particular $(u^{-2} \eta^2 \cdot \eta)^j = u^{-2j} \eta^{3j} \neq 0$.
\end{proof}

\begin{prop} \label{lem:Pic_1 n=2}
At the prime $2$, the exotic Picard group $\kappa_1$ is $\Z/2$
.
\end{prop}

\begin{figure}
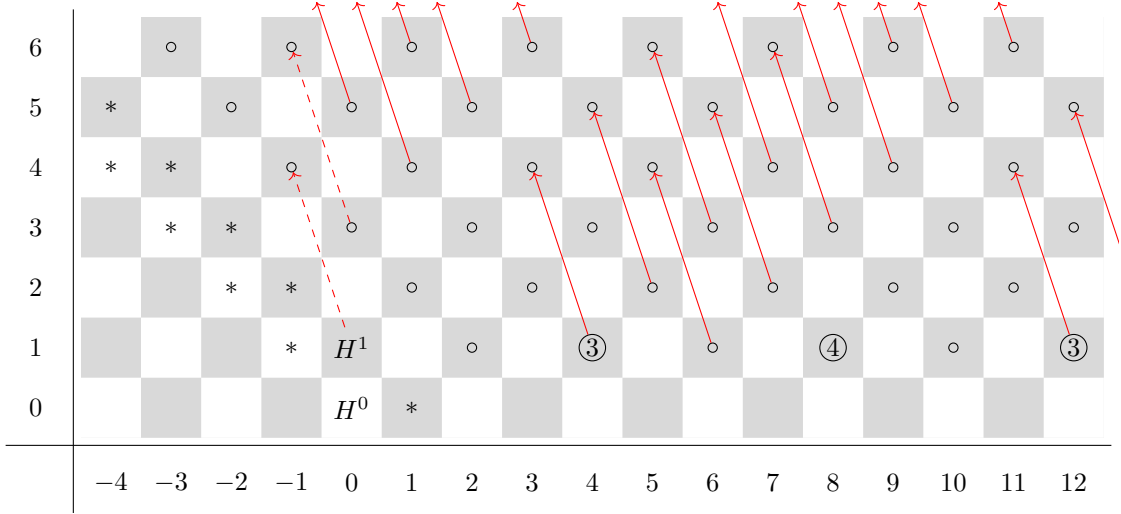

    \centering

    \printpage[x range = {-4}{12}, y range = {0}{6}, grid = chess, page = 2--3, name = PicardSSeven]
    
    \caption[$E_3$-page of Picard spectral sequence at height one, $p=2$]{The $E_3$-page of the Picard spectral sequence at $p=2$. $H^0$ denotes $H^0(\Z[2]^\times, \Pic(\E)) = \Z/2$, and $H^1 = H^1(\Z[2]^\times, (\pi_0 \E)^\times) = \Pic_1^{\mathrm{alg}, 0}$. The possible dashed differentials affects the Picard group; we have omitted some possible differentials with source $t-s \leq -1$.}
    \label{fig:Picard spectral sequence}
\end{figure}

\begin{proof}
    We will deduce this from \cref{lem:differentials in ASS p=2}, which implies that the descent spectral sequence for $\pic(\mathcal E)$ takes the form displayed in Figure \ref{fig:Picard spectral sequence}. According to \cref{lem:algebraic Picard group}, the only differential remaining for the computation of $\kappa_1$ is that on the class in bidegree $(s,t) = (3,0)$, which corresponds to the class $x \in E_{2,+}^{3,2}$ of the Adams spectral sequence. Applying the formula from Lemma \ref{lem:nonlinear differential},
        \[ d_{3}(x) = d_{3}^{ASS} (x) + x^2 = 2x^2 = 0. \]
    After this there is no space for further differentials on $x$, so $\kappa_1 = \Z/2$.
\end{proof}

\begin{sseqdata}[Adams grading, scale = .8, name = HFPSSeven, classes = {show name = {below = 0.1pt}}
]
    \foreach \x in {0,...,1}{                     \class[rectangle](-\x,\x) }
        
        
    \foreach \x in {-10,...,10} {
        \begin{scope}[xshift = 4*\x]
        \foreach \y in {0,...,20}{
            \class(1 - \y, \y + 1)
            }
        \end{scope}
    }
        
    \class[minimum width = width("3") + 0.5em, "3"](3,1)
    \class[minimum width = width("3") + 0.5em, "3"](11,1)
    \class[minimum width = width("3") + 0.5em, "3"](-5,1)
    \class[minimum width = width("3") + 0.5em, "3"](-13,1)
    \class[minimum width = width("3") + 0.5em, "4"](7,1)
    \class[minimum width = width("3") + 0.5em, "4"](-9,1)
    \class[minimum width = width("3") + 0.5em, "5"](15,1)
    
    \foreach \x in {-10,...,10} {
        \begin{scope}[xshift = 4*\x]
        \foreach \y in {0,...,20}{
            \class(-2 - \y, \y + 2) }
        \end{scope}}
            

    \classoptions[class labels = {below = 0.02em}, 
     "u^{-2} \eta^2"](-2,2)

    \classoptions[class labels = {below = 0.02em}, 
    "\eta"](1,1)

    \classoptions[class labels = {below = 0.02em}, 
    "x"](-1,3)
    

    \foreach \x in {-2,...,1} {
        \begin{scope}[xshift = 8*\x]
        \foreach \y in {0,...,10}{
            \d[red]3(5 + \y, 1 + \y)
            }
        \end{scope}
    }
    
    
    \foreach \x in {-1,...,1} {
        \begin{scope}[xshift = 8*\x]
            \d[red]3(3, 1) \replacesource[minimum width = width("3") + 0.5em, "2", label = none]
        \foreach \y in {0,...,10}{
            \d[red]3(4 + \y, 2 + \y)
            }
        \end{scope}
    }
    \begin{scope}
        \foreach \y in {0,...,3}{
            \d[red]3(-7 + \y, 7 + \y)
            }
        \end{scope}

    \foreach \x in {-1,...,1} {
        \begin{scope}[xshift = 8*\x]{
            \structline(3,1)(3,3)
            }
        \end{scope}
    }
\end{sseqdata}

\begin{figure}
    \centering

    \printpage[x range = {-5}{11}, y range = {0}{6}, grid = chess, page = 3, name = HFPSSeven]
    
    \caption[$E_3$-page of descent spectral sequence at height one, $p=2$]{The $E_3$-page of the descent spectral sequence for $\E$ at $p=2$.}
    \label{fig:HFPSS for n=1 even primes}
\end{figure}

\begin{figure}
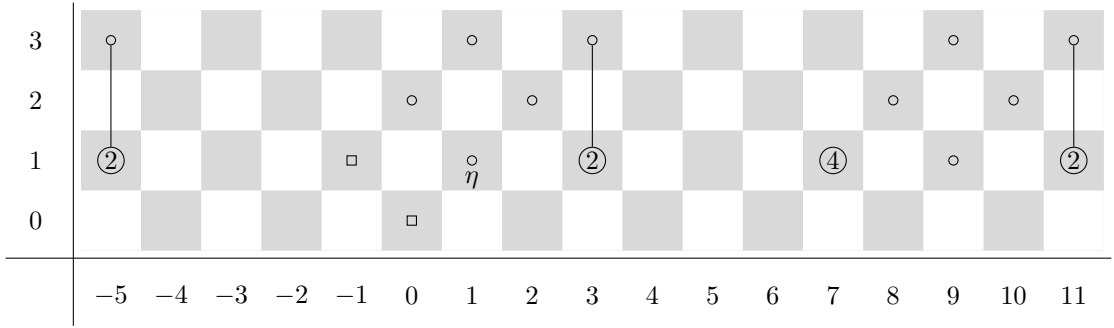

    \centering

    \printpage[x range = {-5}{11}, y range = {0}{3}, grid = chess, page = 4, name = HFPSSeven]
    
    \caption[$E_\infty$-page of descent spectral sequence at height one, $p=2$]{The $E_4 = E_\infty$-page of the descent spectral sequence for $\E$ at $p=2$. 
    }
    \label{fig:HFPSS for n=1 even primes E4}
\end{figure}

\begin{sseqdata}[Adams grading, scale = .8, name = PicardSSeven]
    \foreach \x in {2,...,15}{                     \class[draw = none, "\ast"](1-\x,\x) }
    
    \foreach \x in {1,...,15}{                     \class[draw = none, "\ast"](-\x,\x) }
    
    \class["H^0", draw = none](0,0)
    \class["H^1", draw = none](0,1)
    \class["\ast", draw = none](1,0)
    
        
    \foreach \x in {0,...,10} {
        \begin{scope}[xshift = 4*\x]
        \foreach \y in {0,...,20}{
            \class(2 - \y, \y + 1) }
        \end{scope}}
        
    \class[minimum width = width("3") + 0.5em, "3"](4,1)
    \class[minimum width = width("3") + 0.5em, "3"](12,1)
    \class[minimum width = width("3") + 0.5em, "4"](8,1)
            
    \foreach \x in {1,...,10} {
        \begin{scope}[xshift = 4*\x]
        \foreach \y in {0,...,20}{
            \class(-1 - \y, \y + 2) }
        \end{scope}}

    \foreach \x in {0,...,1} {
        \begin{scope}[xshift = 8*\x]
        \foreach \y in {0,...,10}{
            \d[red]3(6 + \y, 1 + \y)
            }
        \end{scope}
    }
    
        \begin{scope}
        \foreach \y in {0,...,8}{
            \d[red]3(1 + \y, 4 + \y)
            }
        \end{scope}
        \begin{scope}
        \foreach \y in {0,...,3}{
            \d[red]3(-3 + \y, 8 + \y)
            }
        \end{scope}
    
    \foreach \x in {0,...,1} {
        \begin{scope}[xshift = 8*\x]
        \foreach \y in {0,...,11}{
            \d[red]3(4 + \y, 1 + \y)
            }
        \end{scope}
    }
    \begin{scope}
        \foreach \y in {0,...,7}{
            \d[red]3(0 + \y, 5 + \y)
            }
        \end{scope}
    \begin{scope}
        \foreach \y in {0,...,3}{
            \d[red]3(-4 + \y, 9 + \y)
            }
        \end{scope}
        
    \d[red,dashed]3(0,3)
    \d[red,dashed]3(0,1)
    
\end{sseqdata}

%% file: Sections/Brauer.tex
In the previous sections, we considered Galois descent for the Picard spectrum of the $\K$-local category. Recall that the Picard spectrum deloops to the \emph{Brauer spectrum}, which classifies derived Azumaya algebras. In this section we consider Galois descent for the $\K$-local Brauer spectrum; see \cite{baker-richter-szymik, antieau-gepner, gl, antieau-meier-stojanoska} for related work. Unlike most of these sources, the unit in our context is \emph{not} compact, and this makes the corresponding descent statements slightly more delicate. We begin with the basic definitions:

\begin{defn}[\cite{hl}, Definition 2.2.1]
Suppose $\mathcal C$ is a stable homotopy theory, and $R \in \mathcal \CAlg(\mathcal C)$.
\begin{enumerate}
    \item Two $R$-algebras $A, B \in \Alg_R(\mathcal C)$ are \emph{Morita equivalent} if there is an $R$-linear equivalence $\catname{LMod}_A(\mathcal C) \simeq \catname{LMod}_B(\mathcal C)$. We will write $A \sim B$.

    \item An algebra $A \in \Alg_{R}(\mathcal C)$ is an \emph{Azumaya $R$-algebra} if there exists $B \in \Alg_{R}(\mathcal C)$ such that $A \otimes_R B \sim R$.
\end{enumerate}
We will mostly be interested in the case $\mathcal C = \mathcal Sp_\K$ and $R = \bm 1_\K$.
\end{defn}

\begin{remark}
Hopkins and Lurie show \cite[Corollary~2.2.3]{hl} that this definition is equivalent to a more familiar definition in terms of intrinsic properties of $A$. Using \cite[Prop.~2.9.6]{hl} together with faithfulness of the map $\bm 1_{\K} \to \E$, one sees that $A \in \mathcal Sp_\K$ is Azumaya if and only if
\begin{enumerate}
    \item $A$ is nonzero;
    \item $A$ is dualisable in $\Algk R$;
    \item The left/right-multiplication action $A \otimes_R A^{op} \to \operatorname{End}_R(A)$ is an equivalence.
\end{enumerate}
\end{remark}

We will define the space $\Azz(R) \subset \iota \Alg_{R, \K}$ to be the subgroupoid spanned by $\K$-local Azumaya algebras over $R$, and $\Az(R) \coloneqq \pi_0 \Azz(R)$.

\begin{defn}
The \emph{$\K$-local Brauer group} of $R$ is the set
    \[ \Br(R) \coloneqq \Az(R) / \sim, \]
equipped with the abelian group structure $[A] + [B] = [A \otimes_R B]$.
\end{defn}

\begin{remark}
This notation should not be confused with the group of Azumaya $R$-algebras in spectra, which is potentially different.
\end{remark}

In analogy with Picard groups we will write $\Br_h \coloneqq \Br(\bm 1_{\K})$. Our objective is to show how the Picard spectral sequence \cref{eqn:picard spectral sequence} can be used to compute this. First recall the main theorem of \textit{op.\ cit.}:

\begin{thm}[\cite{hl}, Theorem~1.0.11]
    There is an isomorphism
        \[ \Br(\E) \simeq \operatorname{BW}(\kappa) \times \Br'(\E), \]
    where $\operatorname{BW}(\kappa)$ is the \emph{Brauer-Wall} group classifying $\Z/2$-graded Azumaya algebras over $\kappa \coloneqq \pi_0 \E/I_h$, and $\operatorname{Br}'(\E)$ admits a filtration with associated graded $\bigoplus_{k \geq 2} I_h^k/I_h^{k+1}$.
\end{thm}

We will not discuss this result, and instead focus on the orthogonal problem of computing the group of $\K$-local Brauer algebras over the sphere which become Morita trivial over Morava E-theory; in the terminology of \cite{gl} this is the \emph{relative} Brauer group, and will be denoted by $\Br^0_h$. The full Brauer group can then (at least in theory) be obtained from the exact sequence
    \[ 1 \to \Br^0_h \to \Br_h \to \Br(\E)^{\G}. \]
While interesting, the problem of understanding the action of $\G$ on $\Br(\E)$ is somewhat separate, and again we do not attempt to tackle this.

\subsection{Brauer groups and descent} \label{sec:brauer descent}
In this subsection we show that the Picard spectral sequence gives an upper bound on the size of the relative Brauer group, as proven in \cite[Theorem~6.32]{gl} for finite Galois extensions of unlocalised ring spectra. To this end, we define a `cohomological' Brauer group; this might also be called a `Brauer-Grothendieck group' of the $\K$-local category, as opposed to the `Brauer-Azumaya' group discussed above. The ideas described in this section goes back to work of To\"en on Brauer groups in derived algebraic geometry \cite{toen}.

Given $R \in \CAlg(\mathcal Sp_{\K})$, the \inftycat{} $\Modk{R}$ is symmetric monoidal, and therefore defines an object of $\CAlg_{\mathcal Sp_{\K}}(\Pr^L)$. We will consider the symmetric monoidal \inftycat{}
    \[ \Catk R \coloneqq \Mod_{\Modk R}({\Pr}^L_\kappa), \]
where $\kappa$ is chosen to be large enough that $\Modk{R} \in \CAlg(\Pr^L_\kappa)$.

\begin{defn}
The \emph{cohomological $K$-local Brauer space} of $R$ is the Picard space
    \[ \Brr^{\mathrm{coh}}(R) \coloneqq \Picc(\Catk R). \]
Since $\Pr^L_\kappa$ is presentable \cite[Lemma~5.3.2.9(2)]{ha}, this is once again a small space. In analogy with the Picard case, we also write $\Brr_h^{0, \mathrm{coh}}$ for the full subspace of $\Brr_h^{\mathrm{coh}} \coloneqq \Brr^{\mathrm{coh}}(\bm 1_{\K})$ spanned by invertible $\mathcal Sp_{\K}$-modules $\mathcal C$ for which there is an $R$-linear equivalence $\mathcal C \otimes_{\mathcal Sp_{\K}} \Modk \E \simeq \Modk \E$.
\end{defn}

\begin{defn}
Passing to \inftycats{} of left modules defines a functor $\Alg(\Modk{R}) \to \Catk{R}$, with algebra maps acting by extension of scalars. This restricts by the Azumaya condition to
\begin{equation} \label{eqn:map from az to pic Cat}
    \Azz(R) \to \Brr^{\mathrm{coh}}(R)
\end{equation}
We define $\Brr(R) \subset \Brr^{\mathrm{coh}}(R)$ to be the full subgroupoid spanned by the essential image of $\Azz(R)$. We moreover define $\Brr_h \coloneqq \Brr(\bm 1_{\K})$ and $\Brr_h^0 \coloneqq \Brr_h \cap \Brr_h^{0, coh}$.
\end{defn}

\begin{warning}
When working with plain \Einfty-rings, one can take $\kappa = \omega$ in the definition of the cohomological Brauer group, and so the two groups agree by Schwede-Shipley theory. Indeed, a cohomological Brauer class is then an invertible compactly generated $R$-linear \inftycat{} $\mathcal C$, and in particular admits a \emph{finite set} $\{C_1, \dots, C_n\}$ of compact generators \cite[Lemma~3.9]{antieau-gepner}. Thus $\mathcal C \simeq \Mod_A$ for $A \coloneqq \End_{\mathcal C}(\bigoplus_i C_i)$ an Azumaya algebra. This argument fails in $\Pr^L_\kappa$ for $\kappa > \omega$; on the other hand, $\mathcal Sp_{\K} \not\in \CAlg(\Pr^L_\omega)$ since the unit is not compact. For \emph{relative} Brauer classes, we refer to \cite[§3]{relativebrauer} for an alternative solution.
\end{warning}

We now provide a descent formalism suitable for our context, based on the approach of \cite{gl}. If $\mathcal C$ is a presentably symmetric-monoidal \inftycat{} with $\kappa$-compact unit and $R \in \CAlg(\mathcal C)$, we will write $\Cat_R \coloneqq \Mod_{\Mod_R(\mathcal C)}\left(\Pr^L_\kappa \right)$. We will be interested in descent properties of the functor $\Cat_{(-)} \colon \CAlg(\mathcal C) \to \Prsm$: that is, if $R \to R'$ is a map of commutative algebras, we would like to know how close the functor $\theta$ below is to an equivalence:
\begin{equation}
    \Cat_R \xrightarrow{\theta} \lim \left[ \cosimp[\Cat_{R'}]{\Cat_{R' \otimes_R R'}} \right].
\end{equation}

\begin{lem} \label{lem:fully faithfulness of Fun between module categories}
If $R'$ is a descent $R$-algebra, then $\theta$ is fully faithful when restricted to the full subcategory spanned by \inftycats{} of left modules.
\end{lem}

\begin{proof}
Let $A, A' \in \Alg_R(\mathcal C)$. Writing $\LMod_A = \LMod_A (\mathcal C)$, we have equivalences of \inftycats{}
\begin{align*} \label{eqn:fully faithfulness of Fun between module categories}
    \Fun_R(\LMod_A, \LMod_{A'}) &\simeq \RMod_{A \otimes_{R} {A'}\op} \\
     &\simeq  \lim \RMod_{A \otimes_{R} {A'}\op \otimes_R {R'}^\bullet } \\
     &\simeq  \lim \RMod_{(A \otimes_{R} {R'}^\bullet) \otimes_{{R'}^\bullet} ({A'} \otimes_R {R'}^\bullet)\op } \\
     &\simeq \lim \Fun_{{R'}^\bullet}(\LMod_{A \otimes_{R} {R'}^\bullet}, \LMod_{A' \otimes_{R} {R'}^\bullet}). \numberthis
\end{align*}
Here we have twice appealed to \cite[Theorem~4.8.4.1]{ha} (applied to $A\op$ and ${A'}\op$), and to \cref{lem: basechange of descent} for the second equivalence. Passing to cores, we obtain
\begin{align*}
    \Map_{\Cat_R}(\LMod_A, \LMod_{A'}) &\simeq \lim \Map_{\Cat_{{R'}^\bullet}}(\LMod_{A \otimes_{R} {R'}^\bullet}, \LMod_{A' \otimes_{R} {R'}^\bullet}) \\
    &\simeq \Map_{\lim \Cat_{{R'}^\bullet}}(\LMod_{A \otimes_{R} {R'}^\bullet}, \LMod_{A' \otimes_{R} {R'}^\bullet}). \qedhere
\end{align*}
\end{proof}

\begin{cor} \label{lem:almost descent for mathcal E linear inftycats}
For any covering $S' \to S$ in $\proet \G$, the functor
        \[ \theta: \Catk{\mathcal E(S)} \to \lim \Catk{\mathcal E({S'}^{\times_S \bullet})} \]
    is fully faithful when restricted to left module \inftycats{}. \qed
\end{cor}

Using this, we can give a bound on the size of $\Br_h^0$.

\begin{prop} \label{lem:relative Brauer spectrum is a subspace of lim BPicc(R' bullet)}
The group $\Br_h^{0}$ is a subgroup of $\pi_0 \lim B\Picc(\E^{\bullet+1})$.
\end{prop}

\begin{proof}
If $S' \to S$ is a covering in $\proet \G$, we will write $R \to R'$ for the extension $\mathcal E({S}) \to \mathcal E({S'})$. Writing $\Brr(R \mid R')$ for the full subgroupoid of $\Brr(R)$ spanned by objects $\LMod_{A, \K}$ such that $\LMod_{A, \K} \otimes_{\Modk{R}} \Modk{R'} \simeq \Modk{R'}$, we will exhibit $\Brr(R \mid R')$ as a full subspace of $\lim B\Picc({R'}^\bullet)$. Taking the covering $\G \to *$ gives the desired result on $\pi_0$.

By definition, $\Brr(R)$ is the full subcategory of $\Brr^{\mathrm{coh}}(R)$ spanned by module categories. Using the inclusion of $B\Picc(R')$ as the component of the unit in $\Brr(R')$ we can form the diagram
\begin{align} \label{eqn:almost descent for relative Brauer}
\begin{tikzcd}[ampersand replacement=\&]
	{\Brr(R \mid R')} \\
	\& {\lim B\Picc({R'}^\bullet)} \& {B\Picc(R')} \& {B\Picc(R'\otimes_R R')} \& \cdots \\
	{\Brr (R)} \& {\lim \Brr({R'}^\bullet)} \& {\Brr(R')} \& {\Brr(R'\otimes_R R')} \& \cdots \\
	{\iota \Catk{R}} \& {\lim \iota \Catk{{R'}^\bullet}} \& {\iota \Catk{R'}} \& {\iota \Catk{R'\otimes_R R'}} \& \cdots
	\arrow[hook, from=3-1, to=4-1]
	\arrow[shift left=1, from=4-3, to=4-4]
	\arrow[shift right=1, from=4-3, to=4-4]
	\arrow[shift left=1, from=2-3, to=2-4]
	\arrow[shift right=1, from=2-3, to=2-4]
	\arrow[from=2-4, to=2-5]
	\arrow[shift left=2, from=2-4, to=2-5]
	\arrow[shift right=2, from=2-4, to=2-5]
	\arrow[shift left=2, from=4-4, to=4-5]
	\arrow[shift right=2, from=4-4, to=4-5]
	\arrow[from=4-4, to=4-5]
	\arrow[from=4-2, to=4-3]
	\arrow["{\theta}", from=4-1, to=4-2]
	\arrow[hook, from=1-1, to=3-1]
	\arrow[from=2-2, to=2-3]
	\arrow[dashed, from=1-1, to=2-2]
	\arrow[curve={height=-6pt}, from=1-1, to=2-3]
	\arrow[curve={height=-12pt}, from=1-1, to=2-4]
	\arrow["{(\star)}", hook, from=3-1, to=3-2]
	\arrow[from=3-2, to=3-3]
	\arrow[hook, from=2-2, to=3-2]
	\arrow[hook, from=3-2, to=4-2]
	\arrow[hook, from=2-3, to=3-3]
	\arrow[hook, from=3-3, to=4-3]
	\arrow[hook, from=2-4, to=3-4]
	\arrow[hook, from=3-4, to=4-4]
	\arrow[shift left=2, from=3-4, to=3-5]
	\arrow[shift right=2, from=3-4, to=3-5]
	\arrow[from=3-4, to=3-5]
	\arrow[shift left=1, from=3-3, to=3-4]
	\arrow[shift right=1, from=3-3, to=3-4]
\end{tikzcd}
\end{align}
Here hooked arrows denote fully faithful functors: this is essentially by definition in most cases, with the starred functor being fully faithful by virtue of \cref{lem:almost descent for mathcal E linear inftycats} (and the fact that passing to the core preserves limits). The dashed arrow, which clearly exists, is fully faithful by 2-out-of-3.
\end{proof}

In particular, the $(-1)$-stem in the descent spectral sequence for the Picard sheaf $\pic(\mathcal E)$ gives an upper bound on the size of the relative Brauer group. We will draw consequences from this in the next subsection.

As we now discuss, the cohomological Brauer space also admits a description in terms of the pro\'etale site. We will not use, but include a proof for completeness.

\begin{lem} \label{lem:BPicc is a sheaf on Free}
The restriction of the presheaf $B\Picc(\mathcal E)$ to $\Free \G$ is a hypercomplete sheaf of connective spectra.
\end{lem}

\begin{proof}
We will prove descent of the \v Cech nerve for a covering $T' \times \G \to T \times \G$. That is, we would like to show that the following is a limit diagram:
    \[ \augcosimp[B\Picc(\mathcal E(T))]{B\Picc(\mathcal E(T'))}{B\Picc(\mathcal E(T' \times_T T'))} \]
We will consider the Bousfield-Kan spectral sequence for the limit of the \v Cech complex, which reads
    \[ E_2^{s,t} = \pi^s \pi_t B\Picc(\mathcal E({T'}^{\times_T \bullet + 1})) \implies \pi_{t-s} \lim B\Picc(\mathcal E({T'}^{\times_T \bullet + 1})) \numberthis \label{eqn:BKSS for showing BPicc is a sheaf on Free} \]

By \cref{lem:pi_* pic mathcal E}, the homotopy presheaves of $B\Picc (\mathcal E)$ are
\begin{align*}
    \pi_t B\Picc(\mathcal E) = \Cont_\G (-, \pi_t B\Picc(\E)).
\end{align*}
We claim that $E_2^{s,s-1} = 0$ for $s > 0$, so that the map
    \[ B \Picc(\mathcal E(T)) \to \lim B\Picc(\mathcal E({T'}^{\times_T \bullet + 1})) \]
is an equivalence; note that the only possible difference is on $\pi_0$, since $\tau_{\geq 1} B\Picc(\mathcal E)$ is a sheaf of connected spaces. We are therefore interested in the cosimplicial abelian groups
    \[ \cosimp[\Cont(T', \pi_{s-1} B\Picc(\E))]{\Cont(T' \times_T T', \pi_{s-1} B\Picc(\E))}. \numberthis \label{eqn:cosimplicial preasheaf for showing BPicc is a sheaf on Free} \]
This is split when $T' \to T$ is a covering of finite sets, since there is a section at the level of coverings. As a result, by passing to a limit of finite coverings we see that \cref{eqn:cosimplicial preasheaf for showing BPicc is a sheaf on Free} is also split. In particular, the descent spectral sequence \cref{eqn:BKSS for showing BPicc is a sheaf on Free} collapses immediately to the $0$-line, which implies hyperdescent by \cite[Proposition~2.25]{cm}.
\end{proof}

\begin{defn}
Write $\Brr(\mathcal E \mid \E) \in \Sh[\mathcal S]{\proet \G}$ for the sheafification of $S \mapsto B\Pic(\mathcal E(S))$ on $\proet \G$. By \cref{lem:BPicc is a sheaf on Free}, no sheafification is required on the subsite $\Free \G$.
\end{defn}

\begin{cor}
For any closed subgroup $U \subset \G$, we have
    \[ \Brr^{\mathrm{coh}}(\E^{hU} \mid \E) \simeq \Gamma(\G/U, \Brr(\mathcal E \mid \E)). \]
In particular, $\Brr^{\mathrm{coh}}(\bm 1 \mid \E) \simeq \Gamma \Brr(\mathcal E \mid \E)$. \qed
\end{cor}

\subsection{Brauer groups at height one} \label{sec:Brauer n=1}

Using the descent results of the \cref{sec:brauer descent}, we now give bounds on the size of the relative Brauer groups $\Br_1^0$. As usual, the story differs depending on the parity of the prime. In fact, we give conjectural descriptions for some nontrivial Brauer elements; in upcoming work we will elaborate on these computations.

\subsubsection{Odd primes}
\begin{figure}
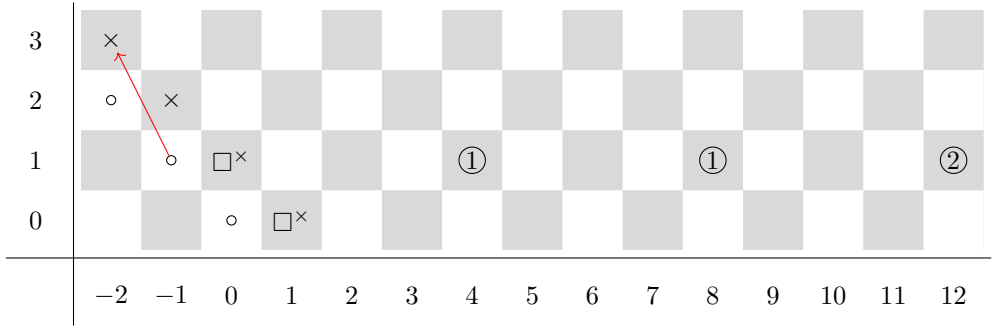

    \centering
    \printpage[grid = chess, xrange = {-2}{12}, yrange={0}{3}, page = 2, name = PicardSSOdd]
    \caption[Picard spectral sequence at height one, odd primes]{The height one Picard spectral sequence for odd primes (implicitly at $p = 3$). Classes are labelled as follows: $\circ = \Z/2$, $\square^\times = \Z[p]^\times$, $\times = \mu_{p-1}$, and circles denote $p$-power torsion (labelled by the torsion degree). Since $\Pic_1 \cong \Pic_1^{\mathrm{alg}} \cong \Z[p]^\times$, no differentials can hit the $(-1)$-stem. Differentials with source in stem $t-s \leq -2$ have been omitted.}
    \label{fig:Picard spectral sequence p odd}
\end{figure}

We first consider the case $p > 2$. The starting page of the Picard spectral sequence is recorded in \cref{lem:ASS p odd} (and computed in \cref{lem:ASS p odd app}). Using this, we obtain the following:

\begin{lem}
    At odd primes, the starting page of the Picard spectral sequence is given by
    \begin{align}
        E_2^{s,t} = H^s(\Z[p]^\times, \pi_t \pic(\E)) = \left\{ \begin{array}{ll}
             \Z/2 & t = 0 \text{ and } s \geq 0 \\
             \Z[p]^\times & t = 1 \text{ and } s = 0,1 \\
             \mu_{p-1} & t = 1 \text{ and } s \geq 2 \\
             \Z/p^{\nu_p(t) + 1} & t = 2(p-1)t'+1 \neq 1 \text{ and } s = 1
        \end{array} \right.
    \end{align} 
    This is displayed in \cref{fig:Picard spectral sequence p odd}. In particular, the spectral sequence collapses for degree reasons at the $E_3$-page.
\end{lem}

\begin{sseqdata}[Adams grading, scale = .8, name = PicardSSOdd ]
    \foreach \x in {0,...,5} {
        \class(-\x,\x)
        \class["\times", draw=none](-1-\x,2+\x)
    }
    \foreach \x in {0,...,1} {
        \class["\square^\times", draw=none](1-\x,\x)
    }
    \class[minimum width = width("2") + 0.5em,"1"](4,1)
    \class[minimum width = width("2") + 0.5em,"1"](8,1)
    \class[minimum width = width("2") + 0.5em,"2"](12,1)

    \d[red]2 (-1,1)
\end{sseqdata}

\begin{proof}
All that remains to compute is internal degrees $t = 0, 1$. We invoke the Lyndon-Hochschild-Serre spectral sequence, which collapses since the extension is split. In particular,
\begin{align*}
    H^*(\Z[p]^\times, \pi_0 \pic(E)) &\simeq H^*(\mu_{p-1}, H^*(\Z[p], \Z/2)) \\
    &\simeq H^*(\mu_{p-1}, \Z/2), \\
    H^*(\Z[p]^\times, \pi_1 \pic(E)) &\simeq H^*(\mu_{p-1}, H^*(\Z[p], \Z[p]^\times))\\
    &\simeq H^*(\mu_{p-1}, \Z[p]^\times \oplus \Sigma \Z[p]). \qedhere
\end{align*}
\end{proof}

Using the results of \cref{sec:brauer descent}, we obtain an upper bound on the relative Brauer group:

\begin{cor} \label{lem:Brauer n = 1 p odd}
At odd primes, $\Br_1^0$ is isomorphic to a subgroup of $\mu_{p-1}$.
\end{cor}

\begin{figure}
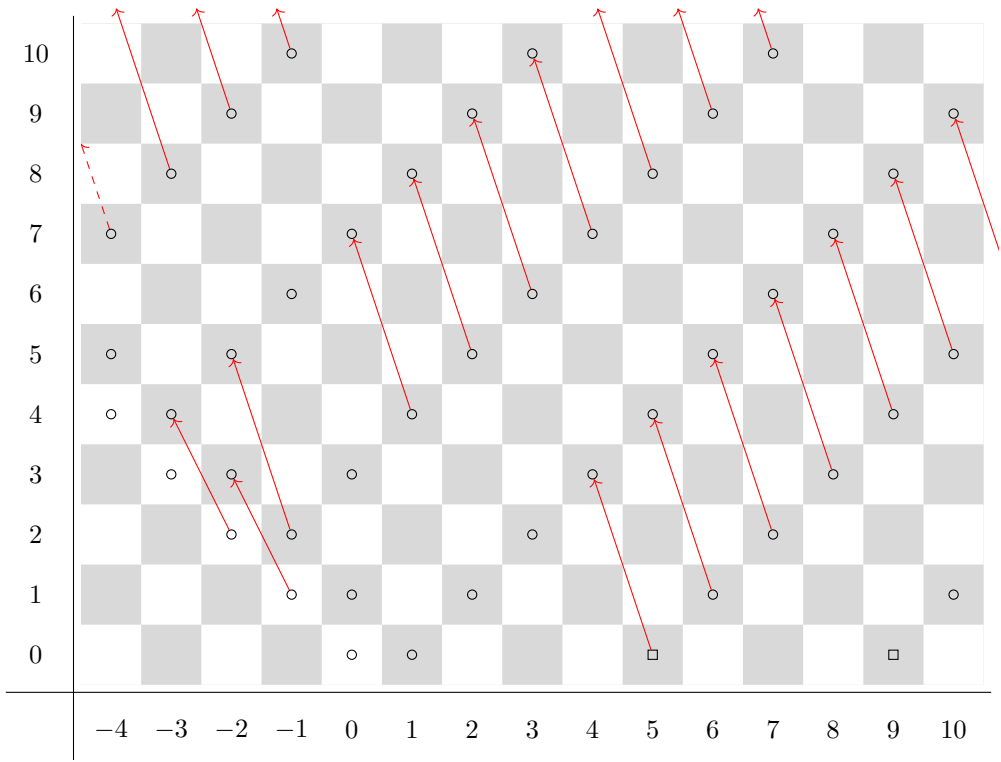

    \centering

    \printpage[x range = {-4}{10}, y range = {0}{10}, grid = chess, page = 2--3, name = BrauerKO]
    
    \caption[Picard spectral sequence for $KO$]{The $E_3$-page of the Picard spectral sequence for $KO$, as in \cite[Figure~7.2]{gl}.}
    \label{fig:Brauer spectral sequence for KO}
\end{figure}

\begin{proof}
The only possible differentials are $d_2$-differentials on classes in the $(-1)$-stem; note that there are no differentials \emph{into} the $(-1)$-stem, since every $E_2$-class in the $0$-stem is a permanent cycle. The generator in $E_2^{1,0}$ supports a $d_2$, since this is the case for the class in $E_2^{1,0}$ of the descent spectral sequence for the $C_2$-action on $KU$ \cite[Prop.~7.15]{gl} (this is displayed in \cref{fig:Brauer spectral sequence for KO}), and the span of Galois extensions
\[\begin{tikzcd}[ampersand replacement=\&]
	{\bm 1_{\K}} \& {KO_p} \& KO \\
	{KU_p} \& {KU_p} \& KU
	\arrow[from=1-2, to=2-2]
	\arrow[from=1-1, to=2-1]
	\arrow[from=1-1, to=1-2]
	\arrow[from=2-1, to=2-2]
	\arrow[from=1-3, to=1-2]
	\arrow[from=2-3, to=2-2]
	\arrow[from=1-3, to=2-3]
\end{tikzcd}\]
allows us to transport this differential. Note that the induced span on $E_2$-pages is
\[\begin{tikzcd}[ampersand replacement=\&,column sep=tiny]
	{\mu_{p-1}} \&\&\& {\mu_{p-1}} \&\&\& {\Z/2} \\
	\\
	\& {\Z/2} \&\&\& {\Z/2} \&\&\& {\Z/2}
	\arrow[red, "d_2", from=3-2, to=1-1]
	\arrow[red, "d_2", from=3-5, to=1-4]
	\arrow[red, "d_2", "{=}"', from=3-8, to=1-7]
	\arrow[shorten <=15pt, shorten >=15pt, hook', from=1-7, to=1-4]
	\arrow[shorten <=15pt, shorten >=15pt, no head, equals, from=1-1, to=1-4]
	\arrow[shorten <=15pt, shorten >=15pt, no head, equals, from=3-2, to=3-5]
	\arrow[shorten <=15pt, shorten >=15pt, no head, equals, from=3-8, to=3-5]
\end{tikzcd}\]
in bidegrees $(s,t) = (1,0)$ and $(3,1)$. Thus
    \[ \pi_0 \lim B\Picc(\Modk{\E^{\bullet + 1}}\perf) \cong \mu_{p-1}. \qedhere \]
\end{proof}

\begin{remark} \label{rem:cyclic algebras}
We now give a conjectural description of the possible nonzero elements of $\Br_1^0$. Recall the \emph{cyclic algebra} construction of \cite{baker-richter-szymik}: its input is $(i)$ a finite $G$-Galois extension $R \to S$, $(ii)$ an isomorphism $\chi: G \cong \Z/k$, and $(iii)$ a \emph{strict} unit $u \in \pi_0 \G_m(R) \coloneqq [\Z, \gll(R)]$. From this data, Baker, Richter and Szymik construct a twisting of the $G$-action on the matrix algebra $M_k(S)$, and hence an Azumaya $R$-algebra $A(S, \chi, u)$ such that $A(S, \chi, u) \otimes_R S \simeq M_k(S)$ by taking fixed points. This construction works equally well for a \emph{$\K$-local} cyclic Galois extension.

In \cite{carmeli_strictunits}, Carmeli shows that the subgroup $\mu_{p-1} \subset \Z[p]^\times = (\pi_0 \bm 1_p)^\times$ lifts to strict units, and hence we have $\mu_{p-1} \subset \G_m(\bm 1_{\K})$ too. Fixing $\chi: \mu_{p-1} \cong \Z/p-1$, any root of unity $\omega$ gives an element 
    \[ A(\E^{h(1 + p \Z[p])}, \chi, \omega) \in \Br(\bm 1_{\K} \mid \E^{h(1 + p \Z[p])}) \subset \Br(\bm 1_{\K} \mid \E), \]
and thus natural candidates for realisations of the classes in $E_\infty^{2,1}$. In fact, upcoming work of Shai Keidar shows that this is indeed the case: the cyclic algebra construction provides an isomorphism
    \[ A(\E^{h(1 + p \Z[p])}, \chi, -) : \mu_{p-1} \cong \Br(\bm 1_{\K} \mid \E). \]
We will give an alternative proof of this fact in upcoming work on the $p = 2$ case.
\end{remark}

\begin{remark} \label{rem:strict units of 1 K}
In fact, the machinery developed in \cref{sec:Morava E-theory,sec:descent theory} should be useful in computing the groups $\pi_0 \G_m(\bm 1_{\K})$ in general. Indeed, choose an algebraic closure $\overline{\F}_p$ and recall that the spectrum of strict units of Morava E-theory $\overline{\E} = E(\overline{\F}_p, \Gamma_n)$ based on $\overline{\F}_p$ is
    \[ \overline{\F}_p^\times \oplus \Sigma^{n+1} \Z[p] \]
by \cite[Theorem~8.17]{nullstellensatz}. In the form stated, the result in this form is due to Hopkins and Lurie---Burklund, Schlank and Yuan compute the strict \emph{Picard} spectrum of $\overline \E$. Since $\mathrm{Gal}(\overline{\F}_p/\F[p^h]) = n\widehat{\Z} \simeq \widehat{\Z} $ has finite cohomological dimension, it follows from \cref{lem:proetale descent for hypercomplete etale sheaves} that the extension $\E \to \overline \E$ is $n\widehat{\Z}$-Galois, and so the extension $\bm 1_{\K} \to \overline \E$ is $\overline{\G} \coloneqq \mathrm{Aut}(\Gamma_h, \overline{\F}_p) = \G \rtimes n\widehat{\Z}$-Galois. Taking strict units preserves limits, and so one obtains a descent spectral sequence
    \[ H^s(\proet{\overline{\G}}, \pi_t \G_m(\overline{\mathcal E})) \implies \pi_{t-s} \G_m (\bm 1_{\K}) \]
for the resulting pro\'etale sheaf $\overline{\mathcal E}$. To evaluate this spectral sequence, we must identify the sheaf cohomology on the starting page; this seems to be a surprisingly delicate problem. Assuming however that it is given by the expected group cohomology
    \[ H^*(\overline{\G}, \pi_* \G_m(\overline{\E})), \]
the spectral sequence collapses to yield $\G_m(\bm 1_{\K}) = \mu_{p-1} \oplus \Sigma^{n}\Z[p] \oplus \Sigma^{n+1} \Z[p]$.
\end{remark}


\subsubsection{The case \texorpdfstring{$p=2$}{p=2}}
We now proceed with the computation of the $(-1)$-stem for the even prime.

\begin{figure}
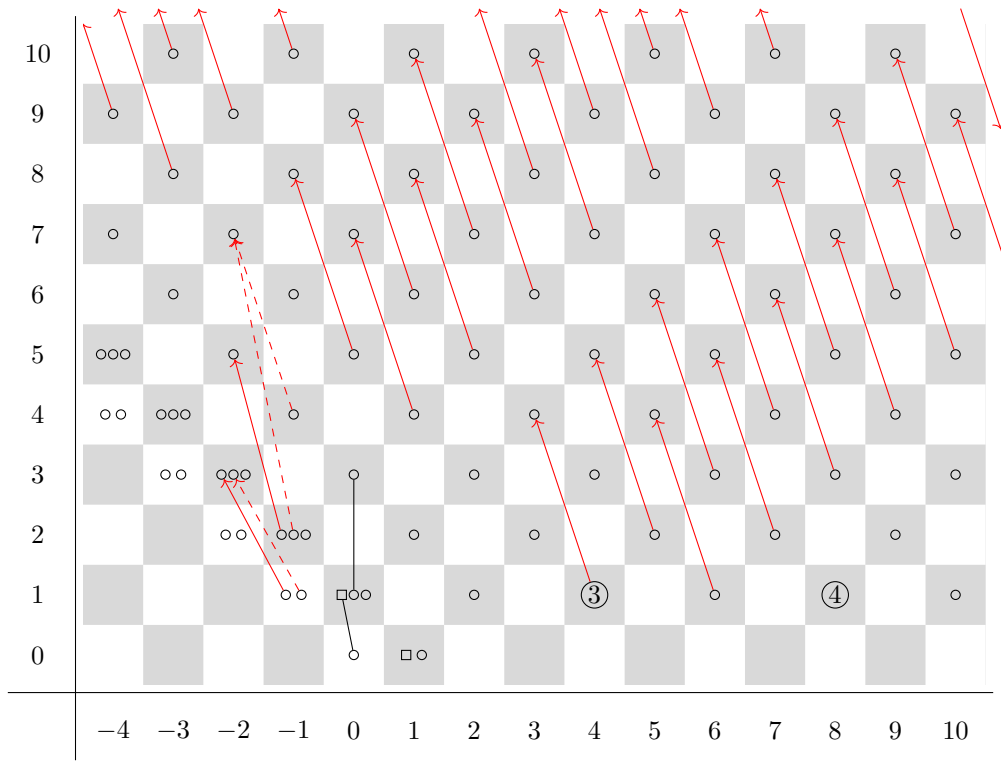

    \centering

    \printpage[x range = {-4}{10}, y range = {0}{10}, grid = chess, page = 2--5, name = Brauersseven]
    
    \caption[Picard spectral sequence at height one, $p=2$]{The $E_3$-page of the Picard spectral sequence at $p=2$. We know that all remaining classes in the $0$-stem survive, by comparing to the algebraic Picard group. Thus the only differentials that remain to compute are those \emph{out} of the $(-1)$-stem; some of these can be transported from the descent spectral sequence for $\Picc(KO)^{hC_2}$---see \cref{fig:Brauer spectral sequence for KO,fig:Brauer spectral sequence for KO2}.}
    
    \label{fig:Brauer spectral sequence p=2}
\end{figure}

\begin{lem}
We have
\begin{align*}
        H^s(\Z[2]^\times, \Pic(\E)) &= \left\{ \begin{array}{ll}
          \Z/2 & s = 0 \\
          (\Z/2)^2 & s \geq 1
    \end{array} \right. \\
    H^s(\Z[2]^\times, (\pi_0 \E)^\times) &= \left\{ \begin{array}{ll}
          \Z[2] \oplus \Z/2 & s = 0\\
          \Z[2] \oplus (\Z/2)^2 & s = 1 \\
          (\Z/2)^3 & s \geq 2
    \end{array} \right.
\end{align*}
The resulting spectral sequence is displayed in Figure \ref{fig:Brauer spectral sequence p=2}.
\end{lem}

\begin{proof}
We need to compute $H^*(\Z[2]^\times, \Z/2)$ and $H^*(\Z[2]^\times, \Z[2]^\times)$. Again we will use the LHSSS; for the first this reads
\begin{align*}
     E_2^{i,j} = H^i(C_2, H^j(\Z[2], \Z/2)) \cong \left[\Z/2[x, y]/x^2 \right]^{i,j} \implies H^{i+j}(\Z[2]^\times, \Z/2),
\end{align*}
The generators have $(i,j)$-bidegrees $\lvert x \rvert = (0,1)$ and $\lvert y \rvert = (1,0)$ respectively, and the spectral sequence is therefore determined by the differential $d_2(x) = \lambda y^2$, where $\lambda \in \Z/2$; since $H^1(\Z[2]^\times, \Z/2) = \Hom{\Z[2]^\times, \Z/2} = \Z/2 \oplus \Z/2$, we deduce that $\lambda = 0$. Finally, the potential extension $2x = y$ is ruled out by the fact that $\Z/2[y] = H^*(C_2, \Z/2)$ is a split summand.

For the second group, the splitting $\Z[2]^\times \cong \Z[2] \times \Z/2$ gives a summand isomorphic to $H^*(\Z[2]^\times, \Z/2)$. The complement is computed by the LHSSS
\begin{align*}
     E_2^{i,j} = H^i(C_2, H^j(\Z[2], \Z[2])) = \left[ \Z[2][w, z]/(w^2, 2z) \right]^{i,j} \implies H^{i+j}(\Z[2]^\times, \Z[2]),
\end{align*}
where now $\lvert w \rvert = (0,1)$ and $\lvert z \rvert = (2,0)$. In this case, the spectral sequence collapses immediately (again by computing $H^1$), with no space for extensions.
\end{proof}

\begin{cor} \label{lem:Brauer n = 1 p = 2}
At the prime two, $\lvert \Br_1^0 \rvert = 2^j$ for $j \leq 5$.
\end{cor}

\begin{proof}
Once again, what remains is to determine the possible differentials on classes in the $(-1)$-stem; these are displayed in \cref{fig:Brauer spectral sequence p=2}. Of these, two can be obtained by comparison to the Picard spectral sequence for the Galois extension $KO_2 \to KU_2$. There are three remaining undetermined differentials, which appear as the dashed arrows in \cref{fig:Brauer spectral sequence p=2}. Invoking \cref{lem:relative Brauer spectrum is a subspace of lim BPicc(R' bullet)} now gives the claimed upper bound.
\end{proof}

\begin{remark}
We comment on the realisation of classes in $\Br_1^0$ at the prime two. Firstly, observe that the surviving class in $E_\infty^{6,5}$ corresponds to the nontrivial Brauer class in $\Br(KO \mid KU)$, and in particular one might hope to show it is descended from a nonzero element of $\Br(KO_2)$. On the other hand, the conjectural computation of $\pi_0 \G_m(\bm 1_{\K})$ (\cref{rem:strict units of 1 K}) would imply that $\pi_0 \G_m(\bm 1_{\K}) = \Z/2 \langle 1 + \varepsilon \rangle$, where $\pi_0 \bm 1_{\K} = \Z[2][\varepsilon]/(2\varepsilon, \varepsilon^2)$ (see \cite[Corollary~5.5.5]{carmeli-yuan}); in particular, this unit becomes trivial in the $C_2$-Galois extension $\E^{h(1 + 4 \Z[2]}$, and is detected in positive filtration of the descent spectral sequence. Moreover, there is at least one further $E_\infty$-class in the Brauer spectral sequence, and thus the two candidates above are not enough to determine everything. Realisation of Brauer classes is therefore significantly more subtle, and we discuss this in upcoming work.

\end{remark}

\begin{sseqdata}[Adams grading, scale = .8, name = Brauersseven]
    \foreach \x in {2,...,15}{                         \class(1-\x,\x)
        \class(1-\x,\x)
        \class(1-\x,\x)}
    
    \foreach \x in {1,...,15}{                     \class(-\x,\x)
    \class(-\x,\x)}
    
    \class(0,0)
    \class[rectangle](0,1)
    \class(0,1)
    \class(0,1)
    \class[rectangle](1,0)
    \class(1,0)
    
        
    \foreach \x in {0,...,10} {
        \begin{scope}[xshift = 4*\x]
        \foreach \y in {0,...,20}{
            \class(2 - \y, \y + 1) }
        \end{scope}}
        
    \class[minimum width = width("3") + 0.5em, "3"](4,1)
    \class[minimum width = width("3") + 0.5em, "3"](12,1)
    \class[minimum width = width("3") + 0.5em, "4"](8,1)
            
    \foreach \x in {1,...,10} {
        \begin{scope}[xshift = 4*\x]
        \foreach \y in {0,...,20}{
            \class(-1 - \y, \y + 2) }
        \end{scope}}

    \foreach \x in {0,...,1} {
        \begin{scope}[xshift = 8*\x]
        \foreach \y in {0,...,10}{
            \d[red]3(6 + \y, 1 + \y)
            }
        \end{scope}
    }
    
        \begin{scope}
        \foreach \y in {0,...,8}{
            \d[red]3(1 + \y, 4 + \y)
            }
        \end{scope}
        \begin{scope}
        \foreach \y in {0,...,3}{
            \d[red]3(-3 + \y, 8 + \y)
            }
        \end{scope}
    \d[red,dashed]3(-1,4) \replacetarget
    \d[red,dashed]2(-1,1,2,2)
    \d[red]2(-1,1)
    \d[red]3(-1,2)
    \d[red,dashed]5(-1,2,2)
    
    \foreach \x in {0,...,1} {
        \begin{scope}[xshift = 8*\x]
        \foreach \y in {0,...,11}{
            \d[red]3(4 + \y, 1 + \y)
            }
        \end{scope}
    }
    \begin{scope}
        \foreach \y in {0,...,7}{
            \d[red]3(0 + \y, 5 + \y)
            }
        \end{scope}
    \begin{scope}
        \foreach \y in {0,...,3}{
            \d[red]3(-4 + \y, 9 + \y)
            }
        \end{scope}

    \structline(0,0)(0,1,1)
    \structline(0,1,2)(0,3)
\end{sseqdata}

\begin{sseqdata}[Adams grading, scale = .8, name = BrauerKO]
    \foreach \x in {2,...,15}{                         \class(1-\x,\x)}
    
    \foreach \x in {1,...,15}{                     \class(-\x,\x)}
    
    \class(0,0)
    \class(0,1)
    \class(1,0)
    
        
    \foreach \x in {0,...,10} {
        \begin{scope}[xshift = 4*\x]
        \foreach \y in {0,...,20}{
            \class(2 - 2*\y, 2*\y + 1) }
        \end{scope}}
            
    \foreach \x in {1,...,10} {
    \class[rectangle](4*\x+1,0)
        \begin{scope}[xshift = 4*\x]
        \foreach \y in {1,...,20}{
            \class(1 - 2*\y, 2*\y) }
        \end{scope}}

    \foreach \x in {0,...,1} {
        \begin{scope}[xshift = 8*\x]
        \foreach \y in {0,...,10}{
            \d[red]3(5 + \y, \y)
            }
        \end{scope}
    }
    
        \begin{scope}
        \foreach \y in {0,...,8}{
            \d[red]3(1 + \y, 4 + \y)
            }
        \end{scope}
        \begin{scope}
        \foreach \y in {0,...,3}{
            \d[red]3(-3 + \y, 8 + \y)
            }
        \end{scope}
        
    \d[red]3(-1,2)
    \d[red,dashed]3(-4,7)
    \d[red]2(-1,1)
    \d[red]2(-2,2)
    

\end{sseqdata}

\begin{sseqdata}[Adams grading, scale = .8, name = BrauerKO2]
    \foreach \x in {1,...,15}{                         \class(1-\x,\x)}
    
    \foreach \x in {1,...,7}{                         \class(1-2*\x, 2*\x)}
    
    \foreach \x in {1,...,15}{                     \class(-\x,\x)}
    
    \class[rectangle](1,0)
    \class(0,0)
    \class(1,0)
    
        
    \foreach \x in {0,...,10} {
        \begin{scope}[xshift = 4*\x]
        \foreach \y in {0,...,20}{
            \class(2 - 2*\y, 2*\y + 1) }
        \end{scope}}
            
    \foreach \x in {1,...,10} {
        \begin{scope}[xshift = 4*\x]
        \class[rectangle](1,0)
        \foreach \y in {1,...,20}{
            \class(1 - 2*\y, 2*\y) }
        \end{scope}}

    \foreach \x in {0,...,1} {
        \begin{scope}[xshift = 8*\x]
        \foreach \y in {0,...,10}{
            \d[red]3(5 + \y, \y)
            }
        \end{scope}
    }
    
        \begin{scope}
        \foreach \y in {0,...,8}{
            \d[red]3(1 + \y, 4 + \y)
            }
        \end{scope}
        \begin{scope}
        \foreach \y in {0,...,3}{
            \d[red]3(-3 + \y, 8 + \y)
            }
        \end{scope}
        
    \d[red]3(-1,2,1)
    \d[red,dashed]3(-4,7)
    \d[red]2(-1,1)
    \d[red]2(-2,2)
    

\end{sseqdata}

\begin{figure}
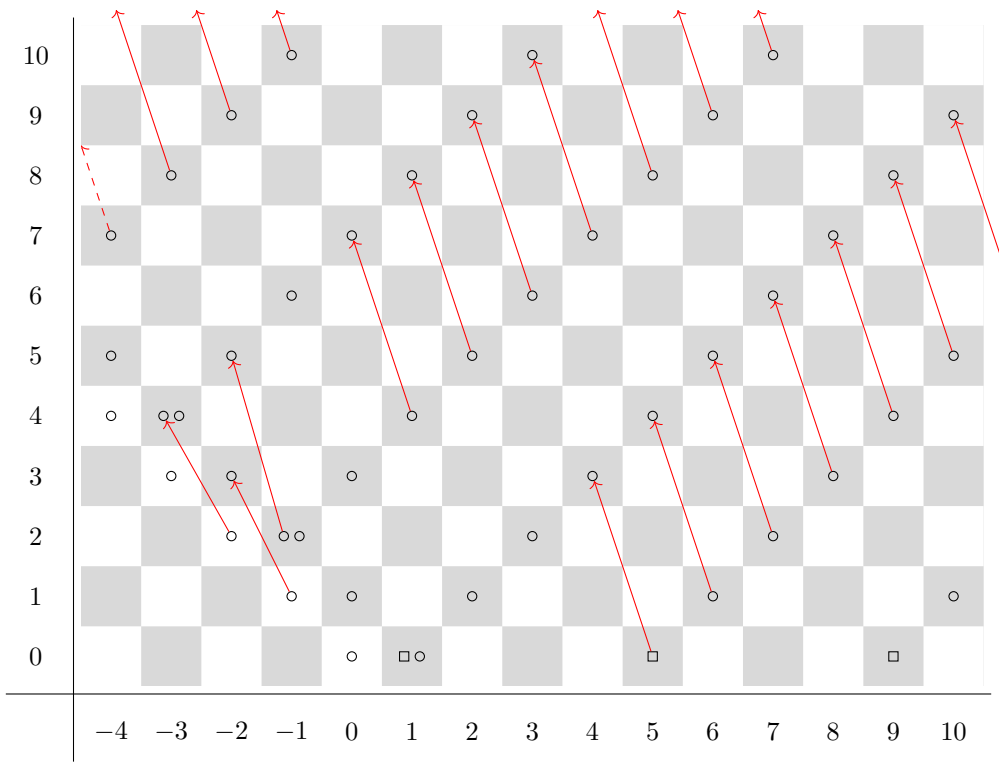

    \centering

    \printpage[x range = {-4}{10}, y range = {0}{10}, grid = chess, page = 2--3, name = BrauerKO2]
    
    \caption[Picard spectral sequence for $KO_2$]{The Picard spectral sequence for $KO_2$.}
    \label{fig:Brauer spectral sequence for KO2}
\end{figure}


%% file: Sections/zDecalage.tex
We make the derivation of the descent spectral sequence a little more explicit, and relate it to the spectral sequence obtained from the covering $\G \twoheadrightarrow *$. For clarity most of this section is written in a general context, but the main result is \cref{lem:decalage for sheaves of spectra appendix}, which will be used to relate the descent spectral sequence to the $\K$-local $\E$-Adams spectral sequence through d\'ecalage.

Let $\mathcal C$ be a site, and write $
\mathcal A \coloneqq \Sh[\mathcal Sp]{\mathcal C}$. Given any object $\mathcal F \in \mathcal A$, there are two natural filtrations one can consider:

\begin{enumerate}
    \item The usual t-structure on $\mathcal Sp$ induces one on $\mathcal A$, defined by the property that $\mathcal F' \in \mathcal A$ is $t$-truncated if and only if $\mathcal F'(X) \in \mathcal Sp_{\leq t}$ for every $X \in \mathcal C$. One can therefore form the Postnikov tower
        \[ \mathcal F \simeq \lim( \cdots \to \tau_{\leq t} \mathcal F \to \cdots) \]
    in $\mathcal A$, and obtain a tower of spectra on taking global sections.
    
    \item Suppose $U \twoheadrightarrow *$ is a covering of the terminal object. Since $A$ is a sheaf, we can recover $\Gamma \mathcal F$ as the limit of its \v Cech complex for the covering, and this has an associated tower. Explicitly, write $\operatorname{Tot}_q = \lim_{\bm \Delta_{\leq q}}$ so that
        \[ \Gamma \mathcal F \simeq \operatorname{Tot} \mathcal F(U^\bullet) \simeq \lim (\cdots \to \operatorname{Tot}_0 \mathcal F(U^\bullet) ). \]
\end{enumerate}
    

Any tower of spectra $X = \lim(\cdots \to X_t \to \cdots)$ gives rise to a (conditionally convergent) spectral sequence, as in \cref{lem:constructing the descent spectral sequence}. Respectively, in the cases above these read
\begin{equation} \label{eqn:descent spectral sequence for a sheaf}
        E^{s,t}_2 = \pi_{t-s} \Gamma \tau_t \mathcal F \implies \pi_{t-s} \Gamma \mathcal F,
\end{equation}
and
\begin{equation} \label{eqn:Cech spectral sequence for a sheaf}
    \check E_2^{p,q} = \pi_{q-p} f_q \mathcal F(U^\bullet) \implies \pi_{q-p} \Gamma \mathcal F,
\end{equation}
where $f_q$ denotes the fibre of the natural transformation $\operatorname{Tot}_q \to \operatorname{Tot}_{q-1}$. In each case, the $d_r$ differential has bidegree $(r+1, r)$ in the displayed grading.

For our purposes, the first spectral sequence, whose $E_2$-page and differentials are both defined at the level of truncations, is useful for interpreting the descent spectral sequence: for example, we use this in \cref{lem:algebraic Picard group}. On the other hand, the second spectral sequence is more easily evaluated once we know the value of a pro\'etale sheaf on the free $\G$-sets. It will therefore be important to be able to compare the two, and this comparison is made using the \textit{d\'ecalage} technique originally due to Deligne in \cite{hodgeii}. The following material is well known (see for example \cite{levine}), but we include an exposition for convenience and to fix indexing conventions. The d\'ecalage construction of \cite{hedenlund_thesis}, which `turns the page' of a spectral sequence by functorially associating to a filtered spectrum its \emph{decal\'ee} filtration, is closely related but not immediately equivalent.

Recall that any tower dualises to a filtration (this will be recounted below); the proof is cleanest in the slightly more general context of \emph{bifiltered} spectra, and so we will make the connection between \cref{eqn:descent spectral sequence for a sheaf} and \cref{eqn:Cech spectral sequence for a sheaf} explicit after proving a slightly more general result.

We suppose therefore that $X$ is a spectrum equipped with a (complete and decreasing) bifiltration. That is, we have a diagram of spectra $X^{t,q}: \Z[]\op \times \Z[]\op \to \mathcal Sp$ such that $X = \colim_{p,q} X^{t,q}$. We will write $ X^{-\infty,q} \coloneqq \colim_t X^{t, q}$ for any fixed $q$, and likewise $ X^{t,-\infty} \coloneqq \colim_q X^{t, q}$ for fixed $t$. Finally, write $X^{t/t', q/q'} \coloneqq \cofib (X^{t',q'} \to X^{t,q})$.

\begin{prop} \label{lem:decalage of a bifiltered spectrum appendix}
    Let $X$ be a bifiltered spectrum, and suppose that for all $t$ and $q$ we have $\pi_s X^{t/t+1, q/q+1} = 0$ unless $s = t - q$. Then there is an isomorphism
        \begin{equation} \label{eqn:decalage on first page}
            {}^1 E_2^{s,t} \simeq {}^2 E_3^{2s-t,s},
        \end{equation} 
    where the left-hand side is the spectral sequence for the filtration $X = \colim_t X^{t,-\infty}$ and the right-hand for $X = \colim_q X^{-\infty,q}$.
    
    This isomorphism is compatible with differentials, and so extends to isomorphisms
        \[ {}^1 E_r^{s,t} \simeq {}^2 E_{r+1}^{2s-t,s} \]
    for $1 \leq r \leq \infty$.
\end{prop}

\begin{proof}
We begin with the isomorphism \cref{eqn:decalage on first page}: it is obtained by considering the following trigraded spectral sequences, which converge to the respective $E_2$-pages.
\begin{align}
    E_2^{s,t,q} = \pi_{q-s} X^{t/t+1, q+1/q} &\implies  \pi_{q-s} X^{t/t+1, -\infty} = {}^1 E_2^{s+t-q,t}, \label{eqn:decalage trigraded sseq for first filtration}\\
    E_2^{s,t,q} = \pi_{t-s} X^{t+1/t, q+1/q} &\implies \pi_{t-s} X^{-\infty, q+1/q} = {}^2 E_2^{s + q - t,q}. \label{eqn:decalage trigraded sseq for second filtration}
\end{align}
The $d_r$ differentials have $(s,t,q)$-tridegrees $(r+1, 0, r)$ and $(r+1, r, 0)$ respectively. Both spectral sequences therefore take a very simple form, because we have assumed each object $X^{t+1/t, q+1/q}$ is Eilenberg-Mac Lane. They are displayed in \cref{fig:decalage trigraded sseqs}.

\begin{sseqdata}
    [name = firsttrigraded,
    Adams grading, scale = .5,
    classes = {draw = none}, 
    x label = {$q-s$},
    y label = {$s$},
    no ticks
    ]

    \foreach \x in {-1,...,10}{
        \class["*"](\x, 1 - 2*\x)
        }
    \foreach \x in {0,...,10}{
        \d[red]2(\x, 1 - 2*\x)}
        
\end{sseqdata}

\begin{sseqdata}[name = secondtrigraded,
    Adams grading, scale = .5,
    classes = {draw = none}, 
    x label = {$t-s$},
    y label = {$s$},
    no ticks
    ]
    \foreach \x in {0,...,10}{
        \class["*"](\x, 0)
        }
\end{sseqdata}

\begin{figure}
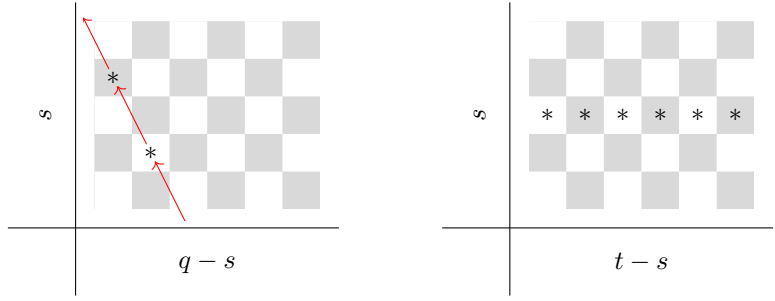

    \centering
    \begin{subfigure}{.3\textwidth}
    \centering
    \printpage[x range = {0}{5}, y range = {-2}{2}, grid = chess, page = 2, name = firsttrigraded]
    \caption{Spectral sequence \cref{eqn:decalage trigraded sseq for first filtration} at fixed $t$. The $y$-intercept is $s = t$.}
    \label{fig:decalage trigraded sseq for first filtration}
    \end{subfigure} \qquad %
    \begin{subfigure}{.3\textwidth}
    \centering
    \printpage[x range = {0}{5}, y range = {-2}{2}, grid = chess, page = 3, name = secondtrigraded]
    \caption{Spectral sequence \cref{eqn:decalage trigraded sseq for second filtration} at fixed $q$. The $y$-intercept is $s = q$.}
    \end{subfigure}
    \caption[Trigraded spectral sequences]{The trigraded spectral sequences computing ${}^1 E_2^{*,t}$ and ${}^2 E_2^{*,q}$ respectively.}
    \label{fig:decalage trigraded sseqs}
\end{figure}

In particular, the first collapses after the $E_2$-page, so that ${}^1E_2^{*,t}$ is the cohomology of the complex
\begin{equation}
    \cdots \to \pi_{t-q} X^{t/t+1, q/q+1} \to \pi_{t-q-1} X^{t/t+1, q+1/q+2} \to \cdots
    \label{eqn:complex from decalage trigraded spectral sequence}
\end{equation}
whose differentials are induced by the composites
    \[ X^{t+1/t, q/q+1} \to \Sigma X^{t+1/t, q+1} \to \Sigma X^{t+1/t, q+1/q+2} .\]
More precisely, ${}^1 E_2^{s,t} = \pi_{t-s} X^{t/t+1, -\infty} \cong H^{s} \left(\pi_{t-*} X^{t/t+1,*/*+1} \right)$.

The second spectral sequence is even simpler, collapsing immediately to give
    \[ {}^2 E_2^{2s-t,s} =\pi_{t-s} X^{-\infty, s/s+1} \cong \pi_{t-s} X^{t/t+1, s/s+1}. \]
We can further identify the first differential on ${}^2E_2^{2s-t,s}$: it is induced by the map $X^{-\infty, 2s-t/2s-t+1} \to \Sigma X^{-\infty, 2s-t+1/2s-t+2}$, and so the identification ${}^1 E_2^{s, t} \cong {}^2 E_3^{2s-t, s}$ follows from the commutative diagram below, in which the top row is part of the complex \cref{eqn:complex from decalage trigraded spectral sequence}.

\[\begin{tikzcd}[column sep = small]
	{\pi_{t-s} X^{t/t+1, s/s+1}} & {\pi_{t-s} \Sigma X^{t/t+1, s+1}} & {\pi_{t-s} \Sigma X^{t/t+1, s+1/s+2}} \\
	{\pi_{t-s} X^{-\infty, s/s+1}} & {\pi_{t-s} \Sigma X^{-\infty, s+1}} & {\pi_{t-s} \Sigma X^{-\infty, s+1/s+2}}
	\arrow[from=2-1, to=2-2]
	\arrow[from=2-2, to=2-3]
	\arrow[from=1-1, to=1-2]
	\arrow[from=1-2, to=1-3]
	\arrow[from=1-1, to=2-1, "\sim"']
	\arrow[from=1-2, to=2-2]
	\arrow[from=1-3, to=2-3, "\sim"]
\end{tikzcd}\]

We next argue that this extends to an isomorphism of spectral sequences ${}^1 E_r^{s,t} \cong {}^2 E_{r+1}^{2s-t,s}$. To do so we will give a map of exact couples as below:
\begin{align*}
    \begin{array}{ccc}
         \left( 
\begin{tikzcd}[ampersand replacement = \&, column sep = tiny]
	{{}^1 D_2^{s,t}} \&\& {{}^1 D_2^{s,t}} \\
	\& {{}^1 E_2^{s,t}}
	\arrow[from=1-1, to=1-3]
	\arrow[from=1-3, to=2-2]
	\arrow[from=2-2, to=1-1]
\end{tikzcd} \right) & \to & \left( 
\begin{tikzcd}[ampersand replacement = \&, column sep = tiny]
	{{}^2 D_3^{2s-t,s}} \&\& {{}^2 D_3^{2s-t,s}} \\
	\& {{}^2 E_3^{2s-t,s}}
	\arrow[from=1-1, to=1-3]
	\arrow[from=1-3, to=2-2]
	\arrow[from=2-2, to=1-1]
\end{tikzcd} \right)
    \end{array}
\end{align*}
By definition of the derived couple on the right-hand side, this amounts to giving maps
    \[ \pi_{t-s} X^{t, -\infty} \to \im(\pi_{t-s} X^{-\infty, s} \to \pi_{t-s} X^{-\infty, s-1}) \]
for all $s$ and $t$, subject to appropriate naturality conditions. To this end we claim that the first map in each of the spans below induces an isomorphism on $\pi_{t-s}$:
    \[ X^{t, -\infty} \leftarrow X^{t, s-1} \to X^{-\infty, s-1}. \]

Indeed, writing $C \coloneqq \cofib(X^{t, s-1} \to X^{t, -\infty})$, we have a filtration $C = \colim C^{t'}$, where $t' \geq t$ and $C^{t'} = \cofib(X^{t', q-1} \to X^{t', -\infty})$. The resulting spectral sequence reads
\begin{equation} \label{eqn:decalage spectral sequence for connectivity argument}
    \pi_{t'-s'} C^{t'/t'+1} \implies \pi_{t'-s'} C \qquad (t' \geq t), 
\end{equation} 
and its $E_2$-page is in turn computed by a trigraded spectral sequence
\begin{equation} \label{eqn:decalage truncated trigraded spectral sequence}
    \pi_{q'-s'} X^{t'/t'+1, q'/q'+1} X \implies \pi_{q'-s'} C^{t'/t'+1} \qquad (q' \leq s-2).
\end{equation}
Spectral sequence \cref{eqn:decalage truncated trigraded spectral sequence} can be thought of as a truncation of \cref{fig:decalage trigraded sseq for first filtration}; it is in turn displayed in \cref{fig:decalage truncated trigraded spectral sequence}. the form of its $E_2$-page implies that $C^{t'/t'+1}$ is $(t' - s +2)$-connected, so that \cref{eqn:decalage spectral sequence for connectivity argument} takes the form displayed in \cref{fig:decalage spectral sequence for connectivity argument}.

\begin{sseqdata}
    [name = truncatedtrigraded,
    Adams grading, scale = .5,
    classes = {draw = none}, 
    x label = {$q'-s'$},
    y label = {$s'$},
    no ticks
    ]

    \foreach \x in {1,...,10}{
        \class["*"](2 + \x, 2 - 2*\x)
        }
    \foreach \x in {2,...,10}{
        \d[red]2(2 + \x, 2 - 2*\x)
        }
        
\end{sseqdata}

\begin{sseqdata}[name = decalage connectivity,
    Adams grading, scale = .5,
    classes = {draw = none}, 
    x label = {$t'-s'$},
    y label = {$s'$},
    no ticks
    ]

    \foreach \y in {0,...,10}{
            \foreach \x in {1,...,10}{
                \class["*"](1 + \x + \y, - \y)
        }
    }
    
\end{sseqdata}

\begin{figure}
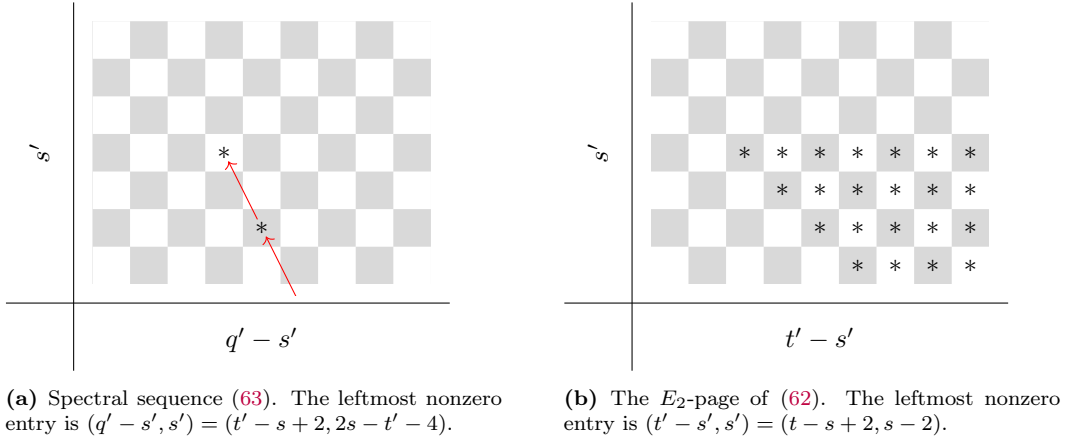

    \centering
    \begin{subfigure}{.4\textwidth}
    \printpage[x range = {0}{8}, y range = {-3}{3}, grid = chess, page = 2, name = truncatedtrigraded]
    \caption{Spectral sequence \cref{eqn:decalage truncated trigraded spectral sequence}. The leftmost nonzero entry is $(q'-s', s') = (t' - s + 2 , 2s - t' - 4)$.}
    \label{fig:decalage truncated trigraded spectral sequence}
    \end{subfigure} \qquad%
    \begin{subfigure}{.4\textwidth}
    \printpage[x range = {0}{8}, y range = {-3}{3}, grid = chess, page = 3, name = decalage connectivity]
    \caption{The $E_2$-page of \cref{eqn:decalage spectral sequence for connectivity argument}. The leftmost nonzero entry is $(t'-s', s') = (t-s+2, s-2)$.}
    \label{fig:decalage spectral sequence for connectivity argument}
    \end{subfigure}
    \caption[Truncated spectral sequences]{Truncated spectral sequences computing $\pi_{*} C$.}
\end{figure}

In particular, $C$ is $t-s+1$-connected and so the map $X^{t, s - 1} \to X^{t, -\infty}$ is $(t-s+1)$-connected. Applying an identical analysis to $X^{t, s} \to X^{t, -\infty}$ shows that this map has $(t-s)$-connected cofibre, and so induces a surjection on $\pi_{t-s}$; the diagram below therefore produces the requisite map ${}^1D_2^{s,t} \to {}^2D_3^{2s-t,s}$.
\[\begin{tikzcd}[column sep = small]
	{\pi_{t-s} X^{t,s}} & {\pi_{t-s} X^{t,s-1}} & {\pi_{t-s} X^{t, -\infty}} \\
	{\pi_{t-s} X^{-\infty, s}} & {\pi_{t-s} X^{-\infty, s - 1}}
	\arrow[two heads, from=1-1, to=1-2]
	\arrow["\sim", from=1-2, to=1-3]
	\arrow[from=1-2, to=2-2]
	\arrow[from=1-1, to=2-1]
	\arrow[from=2-1, to=2-2]
	\arrow[from=1-3, to=2-2, "\exists"]
\end{tikzcd}\]

We will not show that this indeed defines a map of exact couples, since this result is well-known: for example, see \cite[\S6]{levine}.
\end{proof}

We now relate this to the context of Postnikov towers and cobar complexes. If $X \simeq \lim_t( \cdots \to X_t \to \cdots)$ is a convergent tower of spectra with $\colim X_t = 0$, we can form a \emph{dual} filtered spectrum
        \[ X^t \coloneqq \fib (X \to X_{t-1}). \]

Then $\colim X^t \xrightarrow{\sim} X$, and we get another spectral sequence
    \[ \pi_{t-s} f^t X \implies \pi_{t-s} X, \]
where once again $X^{t/t+1} \coloneqq \cofib(X^{t+1} \to X^t)$. Observe that $X^{t/t+1} \simeq X_{t/t-1} \coloneqq \fib(X_t \to X_{t-1})$, by the octahedral axiom:
\[\begin{tikzcd}[column sep = tiny, row sep = tiny]
	& {X^{t+1}} \\
	&&& {X^t} \\
	{} && X &&& {X_{t/t-1}} \\
	&&& {X_t} \\
	& {X_{t-1}}
	\arrow[from=1-2, to=2-4]
	\arrow[from=2-4, to=3-3]
	\arrow[from=1-2, to=3-3]
	\arrow[from=3-3, to=4-4]
	\arrow[from=3-3, to=5-2]
	\arrow[from=4-4, to=5-2]
	\arrow[from=3-6, to=4-4]
	\arrow[from=2-4, to=3-6]
\end{tikzcd}\]

\begin{lem}
    The spectral sequences for $(X_t)$ and $(X^t)$ agree.
\end{lem}

\begin{proof}
The observation above is that the $E_2$-pages agree. To show that the entire spectral sequences match it is enough to show that the differentials $d_2$ do, in other words that the outer diagram below commutes.
\[\begin{tikzcd}[column sep = tiny]
	{X_{t/t-1} = X^{t/t+1}} & {\Sigma X^{t+1}} \\
	{X_t} & {\Sigma X_{t+1/t} = \Sigma X^{t+1/t+2}}
	\arrow[from=1-1, to=2-1]
	\arrow[from=2-1, to=2-2]
	\arrow[from=1-2, to=2-2]
	\arrow[from=1-1, to=1-2]
	\arrow[dashed, from=2-1, to=1-2]
\end{tikzcd}\]

The dashed arrow is given by applying the $3 \times 3$-lemma \cite{neeman_triangulated} to obtain the diagrams below; note that each triangle above only anticommutes, and so the outer square is commutative.

\begin{align*}
    \begin{tikzcd}[ampersand replacement = \&, column sep = tiny]
	{\Sigma^{-1}X^{t+1}} \& {\Sigma^{-1}X^{t}} \& {\Sigma^{-1} X^{t/t+1}} \& {X^{t+1}} \\
	{\Sigma^{-1}X} \& {\Sigma^{-1}X} \& 0 \& X \\
	{\Sigma^{-1} X_t} \& {\Sigma^{-1}X_{t-1}} \& {X_{t/t-1}} \& {X_t} \\
	{X^{t+1}} \& {X^t} \& {X^{t/t+1}} \& {\Sigma X^{t+1}}
	\arrow[from=1-1, to=2-1]
	\arrow[from=2-1, to=3-1]
	\arrow[from=3-1, to=4-1]
	\arrow[from=4-1, to=4-2]
	\arrow[from=3-1, to=3-2]
	\arrow[equals, from=2-1, to=2-2]
	\arrow[from=1-1, to=1-2]
	\arrow[from=1-2, to=1-3]
	\arrow[from=1-3, to=1-4]
	\arrow[from=1-4, to=2-4]
	\arrow[from=2-4, to=3-4]
	\arrow[from=3-4, to=4-4, dotted]
	\arrow[from=4-3, to=4-4]
	\arrow[from=4-2, to=4-3]
	\arrow[from=3-2, to=3-3]
	\arrow[from=3-3, to=3-4]
	\arrow[from=2-2, to=2-3]
	\arrow[from=2-3, to=2-4]
	\arrow[from=1-3, to=2-3]
	\arrow[from=2-3, to=3-3]
	\arrow[from=3-3, to=4-3]
	\arrow[from=1-2, to=2-2]
	\arrow[from=2-2, to=3-2]
	\arrow[from=3-2, to=4-2]
	\arrow["{(-1)}"{description}, draw=none, from=3-3, to=4-4]
\end{tikzcd} & \hskip 20pt &
\begin{tikzcd}[ampersand replacement = \&, column sep = tiny]
	{X_{t-1}} \& {\Sigma X_{t/t-1}} \& {\Sigma X_t} \& {\Sigma X_{t-1}} \\
	{\Sigma X^t} \& {\Sigma X^{t/t+1}} \& {\Sigma^2 X^{t+1}} \& {\Sigma^2 X^t} \\
	{\Sigma X} \& 0 \& {\Sigma^2 X} \& {\Sigma^2X} \\
	{\Sigma X_{t-1}} \& {\Sigma^2 X_{t/t-1}} \& {\Sigma^2 X_t} \& {\Sigma^2 X_{t-1}}
	\arrow[from=1-1, to=2-1, dotted]
	\arrow[from=2-1, to=3-1]
	\arrow[from=3-1, to=4-1]
	\arrow[from=4-1, to=4-2]
	\arrow[from=3-1, to=3-2]
	\arrow[from=2-1, to=2-2]
	\arrow[from=1-1, to=1-2]
	\arrow[from=1-2, to=1-3]
	\arrow[from=1-3, to=1-4]
	\arrow[from=1-4, to=2-4]
	\arrow[from=2-4, to=3-4]
	\arrow[from=3-4, to=4-4]
	\arrow[from=4-3, to=4-4]
	\arrow[from=4-2, to=4-3]
	\arrow[from=3-2, to=3-3]
	\arrow[equals, from=3-3, to=3-4]
	\arrow[from=2-2, to=2-3]
	\arrow[from=2-3, to=2-4]
	\arrow[from=1-3, to=2-3]
	\arrow[from=2-3, to=3-3]
	\arrow[from=3-3, to=4-3]
	\arrow[from=1-2, to=2-2]
	\arrow[from=2-2, to=3-2]
	\arrow[from=3-2, to=4-2]
	\arrow["{(-1)}"{description}, draw=none, from=1-1, to=2-2]
\end{tikzcd}
\end{align*} \qedhere
\end{proof}

By dualising the Postnikov and $\Tot$-towers, it will therefore suffice to verify that the induced bifiltration satisfies the assumptions of \cref{lem:decalage of a bifiltered spectrum appendix}.

\begin{lem} \label{lem:decalage for sheaves of spectra appendix}
    Let $\mathcal F$ be a sheaf of spectra on a site $\mathcal C$, and let $X \twoheadrightarrow *$ be a covering of the terminal object. Suppose that for every $t$ and every $q > 0$ we have $\Gamma(X^q, \tau_t \mathcal F) = \tau_t \Gamma(X^q, \mathcal F).$ Then there is an isomorphism between the descent and Bousfield-Kan spectral sequences, up to reindexing: for all $r$,
        \[ E_r^{s,t} \cong \check E_{r+1}^{2s-t,s}. \]
\end{lem}

\begin{proof}
We form the bifiltrations $(\Gamma \mathcal F)_{t,q} = \Tot_q \Gamma(U^\bullet, \tau_{\leq t} \mathcal F)$. Then
	\[ (\Gamma \mathcal F)_{t,-\infty} = \lim_q \Tot_q \Gamma(U^\bullet, \tau_{\leq t} \mathcal F) = \Tot \Gamma(U^\bullet, \tau_{\leq t} \mathcal F) = \Gamma \tau_{\leq t} \mathcal F, \]
while
	\[ (\Gamma \mathcal F)_{-\infty, q} = \lim_t \Tot_q \Gamma(U^\bullet, \tau_{\leq t} \mathcal F) = \Tot_q \Gamma(U^\bullet, \lim_t \tau_{\leq t} \mathcal F) = \Tot_q \Gamma(U^\bullet, \mathcal F). \]
Applying \cref{lem:decalage of a bifiltered spectrum appendix} to the dual filtration, we need to verify that
    \[ \Tot^{q/q+1} \Gamma(X^\bullet, \tau_t \mathcal F) \]
is Eilenberg-Mac Lane. But recall that for any cosimplicial spectrum $B^\bullet$ we have
    \[ \fib(\Tot_q B^\bullet \to \Tot_{q-1} B^\bullet) \simeq \Omega^q N^q B^\bullet, \]
where $N^q$ denotes the fibre of the map from $X^q$ to the $q$-th matching spectrum; $N^q B^\bullet$ is a pointed space with
    \[ \pi_j N^q B^\bullet = \pi_j B^q \cap \ker s^0 \cap \ker s^{q-1}. \]
In the case of pointed spaces, this fact is \cite[Prop.~X.6.3]{bousfield-kan}; the proof, which appears also as \cite[Lemma~VIII.1.8]{goerss-jardine}, works equally well for a cosimplicial spectrum\footnote{Note that the inductive argument there applies in the `cosimplicial' direction, i.e., in the notation of \loccit one shows for any fixed $t$ that $N^{n,k} \pi_t X = \ker(\pi_t X \to M^{n,k} \pi_t X)$ for $k, n \leq s$.}. By abuse of notation, we also denote this group by $N^q \pi_j B^\bullet$. Thus
\begin{align*}
    \pi_j \Tot^{q/q+1} \Gamma(X^\bullet, \tau_t \mathcal F) &\simeq \pi_j \Omega^q N^q \Gamma(X^\bullet, \tau_t \mathcal F) \\
    &\simeq  N^q \pi_{j + q} \Gamma(X^\bullet, \tau_t \mathcal F) \\
    &\simeq N^q \pi_{j + q} \tau_t \Gamma(X^\bullet, \mathcal F) \\
    & \subset \pi_{j+q} \tau_t \Gamma(X^\bullet, \mathcal F) = \left\{ \begin{array}{ll}
         \pi_t \Gamma(X^\bullet, \mathcal F) & j = t - q \\
         0 & \text{otherwise.}
    \end{array} \right. \qedhere
\end{align*}
\end{proof}

\begin{remark}
On the starting pages, one has
\begin{align*}
    E_2^{s,t} &= \pi_{t-s} \Gamma \tau_t \mathcal F \simeq H^{s}(\mathcal C, \pi_t \mathcal F)
\end{align*}
and
\begin{align*}
    \check E_3^{2s-t, s} &= H^{s} (\pi_{t-*} \Omega^* N^* \Gamma(X^\bullet, \mathcal F) ) = H^{s} ( N^* \pi_{t} \Gamma(X^\bullet, \mathcal F) ) \simeq \check H^s(X \twoheadrightarrow *, \pi_t \mathcal F).
\end{align*}
In particular, note that our assumption implies that the \v Cech-to-derived functor spectral sequence collapses.
\end{remark}

%% file: Sections/zASSodd.tex
At odd primes, the multiplicative lift gives a splitting
    \[ \Z[p]^\times \simeq  (1+ p \Z[p]) \times \mu_{p-1} \simeq  \Z[p] \times \mu_{p-1}, \]
with the second isomorphism given by the $p$-adic logarithm. A pair $(a, b) \in \Z[p] \times \mu_{p-1}$ therefore acts on $\pi_t \E$ as $(b\operatorname{exp}_p(a))^{t}$.

\begin{lem} \label{lem:ASS p odd app}
The starting page of the $\K$-local $\E$-Adams spectral sequence is
\begin{equation}
    E_2^{s,t} = H^{s}(\Z[p]^\times, \pi_{t} \E) = \left\{
        \begin{array}{ll}
             \Z[p] & t = 0 \text{ and } s = 0, 1 \\
             \Z/p^{\nu_p(t') + 1} & t = 2(p-1)t' \neq 0 \text{ and } s = 1
        \end{array} \right.
\end{equation}
and zero otherwise. In particular, it collapses immediately to $E_\infty$.
\end{lem}

\begin{proof}
We will use the Lyndon-Hochschild-Serre spectral sequence \cite[Theorem~4.2.6]{symonds-weigel} for the inclusion $\Z[p] \simeq 1 + p \Z[p] \hookrightarrow \Z[p]^\times$, which reads
    \[ H^i(\mu_{p-1}, H^j(\Z[p], \pi_t KU_p)) \implies  H^{i+j}(\Z[p]^\times, \pi_t KU_p). \]
Since everything is $(p)$-local, taking $\mu_{p-1}$ fixed-points is exact. The spectral sequence therefore collapses, and what remains is to compute $\Z[p]$-cohomology.

By \cite[\S3.2]{symonds-weigel}, the trivial pro-$p$ module $\Z[p]$ admits a projective resolution
    \[ 0 \to \Z[p][[\Z[p]]] \xrightarrow{\zeta - 1} \Z[p][[\Z[p]]] \to \Z[p] \to 0 \]
in the (abelian) category $\mathfrak C_p(\Z[p])$ of pro-$p$ continuous $\Z[p]^\times$-modules; here we write $\zeta$ for a topological generator, and $\Z[p][[G]] \coloneqq \lim_{U <_o G} \Z[p][G/U]$ for the \emph{completed} group algebra of a profinite group $G$. In particular,
\begin{align*}
    H^j(\Z[p], \pi_{2t} KU_p) = H^j \left(\Z[p] \xrightarrow{\exp_p(\zeta)^t - 1} \Z[p] \right) \simeq \left\{ \begin{array}{ll}
         \Z[p] & t = 0 \text{ and } j = 0,1 \\
         \Z[p]/p^{\nu_p(t)+ 1} & t \neq 0 \text{ and } j = 1
    \end{array} \right.
\end{align*}
 where for the final isomorphism we have used the isomorphism $\exp_p: p^{j-1} \Z[p]/p^j \Z[p] \simeq 1 + p^j \Z[p]^\times/1 + p^{j+1} \Z[p]^\times$ to obtain $\nu_p(\exp_p(\zeta)^t-1) = \nu_p(t \zeta) + 1 = \nu_p(t) + 1$. The $\mu_{p-1}$-action has fixed points precisely when it is trivial, i.e. when $p-1 \mid t$, which gives the stated form.
\end{proof}

%% file: Sections/zASSp=2.tex
What changes at even primes? In this case the multiplicative lift is instead defined on $(\Z/4)^\times$, and provides a splitting
    \[ \Z[2]^\times \simeq (1+ 4 \Z[2]) \times C_2 \simeq \Z[2] \times C_2. \]

This implies the following more complicated form for the starting page, since $\mathrm{cd}_2(C_2) = \infty$.

\begin{lem} \label{lem:E2 page of ASS n=1 p=2 app}
The starting page of the descent spectral sequence for the action of $\G$ on $\E$ is given by
\begin{equation} \label{eqn:E2 page of ASS n=1 p=1 app}
    E_2^{s,t} = H^{s}(\Z[2]^\times, \pi_{t} \E) = \left\{
        \begin{array}{ll}
             \Z[2] & t = 0 \text{ and } s = 0,1 \\
             \Z/2 & t \equiv_4 2 \text{ and } s \text{ odd} \\
             \Z/2 & t \equiv_4 0 \text{ and } s > 1 \text{ odd} \\
             \Z/2^{\nu_2(t)+2} & 0 \neq t \equiv_4 0 \text{ and } s = 1
        \end{array} \right.
\end{equation}
and zero otherwise. The result is displayed in \cref{fig:HFPSS for n=1 even primes appendix}, which is reproduced for convenience of the reader.
\end{lem}

\begin{proof}
We will again use the Lyndon-Hochschild-Serre spectral sequence for the inclusion $\Z[2] \simeq 1 + 4 \Z[2] \hookrightarrow \Z[2]^\times$, which reads
    \[ H^i(C_2, H^j(\Z[2], \pi_t KU_2)) \implies  H^{i+j}(\Z[2]^\times, \pi_t KU_2). \]
Since $C_2$ is $2$-torsion, we will have higher $C_2$-cohomology contributions. Nevertheless, the computation of $\Z[2]$-cohomology is identical to the odd-prime case, except that now one has $\nu_2(\exp_2(\zeta)^t - 1) = \nu_2(t) + 2$. Thus
\begin{align*}
    H^j(\Z[2], \pi_t KU_2) = H^j \left(\Z[2] \xrightarrow{\exp_p(\zeta)^t - 1} \Z[2] \right) \simeq \left\{ \begin{array}{ll}
         \Z[2] & t = 0 \text{ and } j = 0,1 \\
         \Z[2]/2^{\nu_2(t) + 2} & t \neq 0 \text{ and } j = 1
    \end{array} \right.
\end{align*}
For $t \neq 0$, the $E_2$-page of the Lyndon-Hochschild-Serre is therefore concentrated in degrees $j = 1$, and so collapses. This yields
    \[ H^s(\Z[2]^\times, \pi_t KU_2) \simeq H^{s-1}(C_2, \Z/2^{\nu_2(t) + 2}), \]
which accounts for most of the groups in \cref{eqn:E2 page of ASS n=1 p=1 app}. For $t = 0$, it instead takes the form
    \[ E_2^{*,*} = \Z[2][x^{(0,1)},y^{(2,0)}]/(2x, 2y, x^2) \implies H^{i+j}(\Z[2]^\times, \pi_0 KU_2). \]
In particular, the entire spectral sequence is determined by the differential $d_2(x) = 0$, which implies that all other differentials vanish by multiplicativity. To deduce this differential, note that the edge map $H^1(\Z[2]^\times, \pi_0 KU_2) \to H^0(C_2, H^1(\Z[2], \pi_0 KU_2))$ can be interpreted as the restriction map
    \[ \Hom{\Z[2]^\times, \Z[2]} \to \Hom{\Z[2], \Z[2]}, \]
along $\exp_2: \Z[2] \hookrightarrow \Z[2]^\times$. This is an isomorphism since $\Hom{\Z/2, \Z[2]} = 0$, so $d_2$ must act trivially on bidegree $(0,1)$.
\end{proof}

Our next task is to compute the differentials. In the rest of the appendix we will prove the following:

\begin{prop} \label{lem:differentials in ASS p=2 app}
The differentials on the third page are as displayed in \cref{fig:HFPSS for n=1 even primes appendix}. The spectral sequence collapses at $E_4$ with a horizontal vanishing line.
\end{prop}

\begin{sseqdata}[Adams grading, scale = .5, name = HFPSSevenApp
]
    \foreach \x in {0,...,1}{                     \class[rectangle](-\x,\x) }
        
        
    \foreach \x in {-10,...,10} {
        \begin{scope}[xshift = 4*\x]
        \foreach \y in {0,...,20}{
            \class(1 - \y, \y + 1)
            }
        \end{scope}
    }
        
    \class[minimum width = width("3") + 0.5em, "3"](3,1)
    \class[minimum width = width("3") + 0.5em, "3"](11,1)
    \class[minimum width = width("3") + 0.5em, "3"](-5,1)
    \class[minimum width = width("3") + 0.5em, "3"](-13,1)
    \class[minimum width = width("3") + 0.5em, "4"](7,1)
    \class[minimum width = width("3") + 0.5em, "4"](-9,1)
    \class[minimum width = width("3") + 0.5em, "5"](15,1)
    
    \foreach \x in {-10,...,10} {
        \begin{scope}[xshift = 4*\x]
        \foreach \y in {0,...,20}{
            \class(-2 - \y, \y + 2) }
        \end{scope}}
    
    \foreach \x in {-1,...,1} {
        \begin{scope}[xshift = 8*\x]
        \foreach \y in {0,...,10}{
            \d[red]3(5 + \y, 1 + \y)
            }
        \end{scope}
    }
    
    \foreach \x in {-1,...,1} {
    \begin{scope}[xshift = 8*\x]
        \foreach \y in {0,...,10}{
            \d[red]3(3 + \y, 1 + \y)
            }
        \end{scope}
    }
\end{sseqdata}

\begin{figure}
    \centering

    \printpage[x range = {-5}{15}, y range = {0}{5}, grid = chess, page = 3, name = HFPSSevenApp]
    
    \caption[$E_3$-page of descent spectral sequence at height one, $p = 2$]{The $E_3$-page of the descent spectral sequence for $\E$ at $p=2$.}
    \label{fig:HFPSS for n=1 even primes appendix}
\end{figure}

\begin{sseqdata}[Adams grading, scale = .5, name = HFPSSKO
]
        
        
    \foreach \x in {-10,...,10} {
        \begin{scope}[xshift = 4*\x]
        \class[rectangle](0,0)
        \foreach \y in {0,...,20}{
            \class(1 + \y, 1+ \y)
            \structline(\y, \y)(1+\y,1+\y)
            }
        
        \end{scope}
    }

    \classoptions[class labels = {left = 0.1em}, "u^2"](4,0)
    \classoptions[class labels = {below = 0.1em}, "\eta"](1,1)

    \foreach \x in {-3,...,3} {
        \begin{scope}[xshift = 8*\x]
        \foreach \y in {0,...,10}{
            \d[red]3(4 + \y, \y)
            }
        \end{scope}
    }

\end{sseqdata}

\begin{figure}
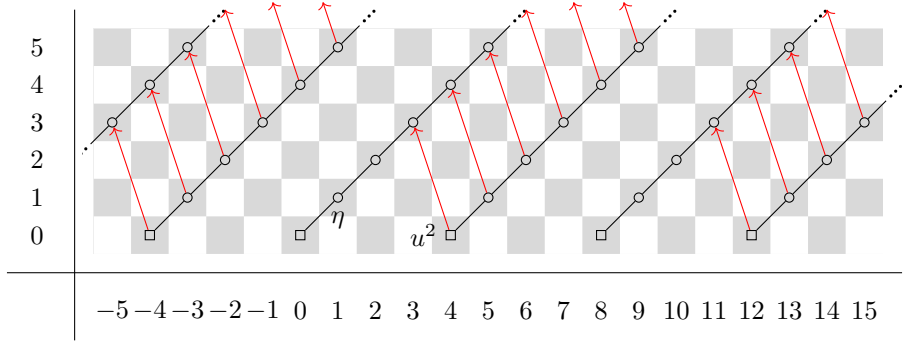

    \centering

    \printpage[x range = {-5}{15}, y range = {0}{5}, grid = chess, page = 3, name = HFPSSKO]
    
    \caption[HFPSS for $KO$]{The HFPSS for the $C_2$-Galois extension $KO \to KU$. The class $\eta$ represents $\eta$ in the Hurewicz image in $\pi_* KO$, and towers of slope one are related by $\eta$-multiplications. In particular, $\eta^4$ cannot survive to $E_\infty$ since $\eta^4 = 0 \in \pi_* \mathbb S$. The only option is $d_3(u^2 \eta) = \eta^4$, which implies the rest by multiplicativity.}
    \label{fig:HFPSS for KO to KU}
\end{figure}

We will compute these differentials by comparing to the HFPSS for the conjugation action on $KU_2$ (\cref{fig:HFPSS for KO to KU}), which reads
    \[ E_2^{s,t} = H^s(C_2, \pi_* KU_2) = \Z[2][\eta, u^{\pm 2}]/2\eta \implies \pi_{t-s} KO_2. \]
The following result folklore:

\begin{lem} \label{lem:fibre sequence for L_K S}
The $\K$-local sphere fits in a fibre sequence
    \[ \bm 1_{\K} \to KO_2 \xrightarrow{\psi^5 - 1} KO_2. \]
\end{lem}

\begin{proof}
We first consider the map $\psi^5 - 1: KU_2 \to KU_2$. Certainly $\psi^5$ acts trivially on $KU_2^{1+ 4\Z[2]} \coloneqq \Gamma(\Z[2]^\times/1 + 4\Z[2], \mathcal E)$,
so there is a map
\begin{align} \label{eqn:fibre of psi - 1}
    KU_2^{h(1 + 4\Z[2])} \to \operatorname{fib}(KU_2 \xrightarrow{\psi^5 - 1} KU_2)    
\end{align}
induced by the inclusion of fixed points $KU_2^{1 + 4
\Z[2]} \to KU_2$. As observed in \cref{lem:E2 page of ASS n=1 p=2 app}, the HFPSS computing $\pi_* KU_2^{h(1 + 4\Z[2])}$ collapses at $E_2$ with horizontal vanishing above filtration one, and one therefore sees that \cref{eqn:fibre of psi - 1} must be an equivalence by computing homotopy groups of $\fib(\psi^5 -1)$ using the exact sequence.

Taking fixed points for the $C_2$ action now yields the result:
\begin{align*}
    \bm 1_{\K} &\simeq (KU_2^{h(1 + 4\Z[2])})^{C_2} \\
    & \simeq \operatorname{fib}(KU_2 \xrightarrow{\psi^5 - 1} KU_2)^{C_2} \\
    & \simeq \operatorname{fib}(KU_2^{C_2} \xrightarrow{\psi^5 - 1} KU_2^{C_2}) \\
    & \simeq \operatorname{fib}(KO_2 \xrightarrow{\psi^5 - 1} KO_2). \qedhere
\end{align*}
\end{proof} 

\begin{proof}[Proof (\cref{lem:differentials in ASS p=2 app})]
The previous lemma gives the diagram
\begin{align} \label{eqn:fibre sequence KO to LS to KO}
\begin{tikzcd}[ampersand replacement = \&]
	{\Sigma^{-1}KO_2} \& {\bm 1_{\K}} \& {KO_2} \\
	{\Sigma^{-1}KU_2} \& {KU_2^{h(1 + 4\Z[2])}} \& {KU_2}
	\arrow[from=1-1, to=2-1]
	\arrow[from=2-1, to=2-2]
	\arrow[from=1-1, to=1-2]
	\arrow[from=1-2, to=2-2]
	\arrow[from=2-2, to=2-3]
	\arrow[from=1-3, to=2-3]
	\arrow[from=1-2, to=1-3]
\end{tikzcd}
\end{align}
in which the top row is obtained as $C_2$-fixed points of the bottom. The HFPSS for the middle map,
\begin{equation} \label{eqn:HFPSS for L_K S to KU^1+4Z2}
    E_2^{s, t} = H^s(C_2, \pi_t KU_2^{1+4\Z[2]}) \implies \pi_{t-s} \bm 1_{\K},
\end{equation}
is very closely to the descent spectral sequence; it is displayed in \cref{fig:mu_2 spectral sequence for KU^1+4 Z_2}. In fact, in \cref{lem:isomorphism between two spectral sequences n=1 p=2} we will show that the two spectral sequences are isomorphic (including differentials), up to a certain filtration shift. To infer the differentials in \cref{fig:HFPSS for n=1 even primes appendix}, it is therefore enough to compute the differentials in \cref{fig:mu_2 spectral sequence for KU^1+4 Z_2}.

Next observe that \cref{lem:fibre sequence for L_K S} implies existence of exact sequences
\begin{align*}
    0 \to \pi_{2t} KU_2^{h(1 + 4\Z[2])} \to \Z[2] & \xrightarrow{5^t - 1} \Z[2] \to \pi_{2t - 1} KU_2^{h(1 + 4\Z[2])} \to 0,
\end{align*}
in which the leftmost term vanishes for $t \neq 0$ and the middle map is null for $t = 0$. Taking $C_2$-cohomology yields isomorphisms
\begin{align}
    H^*(C_2, \pi_0 KU_2^{1 + 4 \Z[2]}) & \xrightarrow{\sim} H^*(C_2, \pi_0 KU_2) \label{eqn:H^*(mu_2 pi_0 KU^1+4Z2)} \\
    H^*(C_2, \pi_0 KU_2) &\xrightarrow{\sim} H^*(C_2, \pi_{-1} KU_2^{1 + 4 \Z[2]}) \label{eqn:H^*(mu_2 pi_-1 KU^1+4Z2)},
\end{align}
and an exact sequence
\begin{equation*}
    0 \to H^{s-1}(C_2, \pi_{2t-1} KU^{h(1 + 4\Z[2])}) \to H^s(C_2, \pi_{2t} KU_2) \to H^s(C_2, \pi_{2t} KU_2) \to H^{s}(C_2, \pi_{2t-1} KU^{h(1 + 4\Z[2])}) \to 0
\end{equation*}
 for $s \geq 1$ (and $t \neq 0$). The middle terms are either both zero or both $\Z/2$, and in the latter case we obtain the following further isomorphisms:
\begin{align}
    H^{s-1}(C_2, \pi_{2t-1} KU_2^{1+4\Z[2]}) & \xrightarrow{\sim} H^s(C_2, \pi_{2t} KU_2) \simeq \Z/2, \label{eqn:H^s-1(mu_2 pi_2t-1 KU^1+4Z2)}\\
    \Z/2 \simeq H^s(C_2, \pi_{2t} KU_2) & \xrightarrow{\sim} H^s(C_2, \pi_{2t-1} \label{eqn:H^s(mu_2 pi_2t-1 KU^1+4Z2)}KU_2^{1+4\Z[2]}).
\end{align}
For $s = 0$ and $t \neq 0$ even, we instead have an exact sequence
\begin{equation}
    0 \to \Z[2] \xrightarrow{5^t-1} \Z[2] \to H^0(C_2, \pi_{2t-1} KU_2^{h(1 + 4\Z[2])}) \to 0. \label{eqn:H^0(mu_2 pi_2t-1 KU^1+4Z2)}
\end{equation} 
Equations \cref{eqn:H^*(mu_2 pi_-1 KU^1+4Z2),eqn:H^*(mu_2 pi_0 KU^1+4Z2),eqn:H^0(mu_2 pi_2t-1 KU^1+4Z2),eqn:H^s(mu_2 pi_2t-1 KU^1+4Z2),eqn:H^s-1(mu_2 pi_2t-1 KU^1+4Z2)} compute the effect of the maps in \cref{eqn:fibre sequence KO to LS to KO} on spectral sequences. Namely:
\begin{enumerate}
    \item The map $\Sigma^{-1} KU_2 \to KU_2^{1+4\Z[2]}$ induces a (filtration preserving) surjection from $E_2(KU)$ onto the unfilled classes in \cref{fig:mu_2 spectral sequence for KU^1+4 Z_2}.
    \item The map $KU_2^{1+4\Z[2]} \to KU_2$ is injective on the solid classes in \cref{fig:mu_2 spectral sequence for KU^1+4 Z_2}. It induces a filtration-preserving isomorphism on the subalgebras in internal degree $t = 0$, and away from this increases filtration by one.
\end{enumerate}
The differentials in \cref{fig:mu_2 spectral sequence for KU^1+4 Z_2} follow almost immediately: on unfilled classes, they are images of differentials in the HFPSS for $KU_2$, and on most solid classes they are detected by differentials in the HFPSS for $KU_2$. We are left to determine a small number of differentials on classes with internal degree $t$ close to zero.

\begin{enumerate}
    \item The exact sequence induced by \cref{lem:fibre sequence for L_K S} shows that the map $\bm 1_{\K} \to KO_2$ induces an isomorphism on $\pi_0$. As a result, the unit in the HFPSS for $KU^{h(1 + 4\Z[2])}$ is a permanent cycle.
    
    \item Write $u^{-2} \eta^2$ for the generator in bidegree $(s,t) = (2,0)$, which maps to a class of the same name in the HFPSS for $KU_2$. This cannot survive in the HFPSS for $KU_2^{h(1 + 4\Z[2])}$: one sees that $\pi_{-2} \bm 1_{\K} = 0$ using the fibre sequence of \cref{lem:fibre sequence for L_K S} and the fact that $\pi_{-2}KO_2 = \pi_{-1} KO_2 = 0$. The only possibility is a nonzero $d_2$ since all other possible differentials occur at $E_4$ or later, when there are no possible targets left (by virtue of the known $d_3$-differentials). Likewise, this implies
        \[ d_2((u^{-2} \eta^2)^{2j+1}) \equiv_2 (u^{-2} \eta^2)^{2j} d_2(u^{-2} \eta^2) \neq 0 \]
    by the Leibniz rule.

    \item Write $z$ for the generator in bidegree $(s, t-s) = (0, -3)$, which is detected in filtration one by $u^{-2} \eta$. As above one computes that $\pi_{-3} \bm 1_{\K} = 0$, and so this class must die; the only option is a nonzero $d_4$ on $z$. In fact, when $j$ is even the exact sequence
        \[ \pi_{4j-2} KO_2 \to \pi_{4j-3} \bm 1_{\K} \to \pi_{4j-3} KO_2 \]
    implies that $\pi_{4j-3} \bm 1_{\K} = 0$, so that $(u^{-2} \eta^2)^{2j}z$ supports a $d_4$ by the same argument. When $j$ is odd the sequence only gives a bound $|\pi_{4j-3} \bm 1_{\K}| \leq 4$, but this is already populated by elements in filtration $s \leq 2$ (which survive by comparison to the HFPSS for $KU_2$). Thus $(u^{-2} \eta^2)^{2j} z$ supports a $d_4$ in these cases too.
\end{enumerate}
After this, the spectral sequence collapses by sparsity.
\end{proof}

\begin{sseqdata}[Adams grading, scale = .5, name = mu2SS
]
    \class[rectangle, fill](0,0)
    \class[rectangle](-1,0)
        
        
    \foreach \x in {0,...,10} {
        \begin{scope}[xshift = 4*\x]
        \foreach \y in {0,...,20}{
            \class[fill](1 + \y, \y)
            }
        \end{scope}
    }

    \foreach \x in {-10,...,-1} {
        \begin{scope}[xshift = 4*\x]
        \foreach \y in {0,...,20}{
            \class[fill](1 + \y - 2*\x, \y - 2*\x)
            }
        \end{scope}
    }

    \class[fill](-3,0)

    \foreach \y in {0,...,2}{
        \class[fill](-7 + \y, \y)
        }
        
    \class[minimum width = width("3") + 0.5em, "3"](3,0)
    \class[minimum width = width("3") + 0.5em, "3"](11,0)
    \class[minimum width = width("3") + 0.5em, "3"](-5,0)
    \class[minimum width = width("3") + 0.5em, "3"](-13,0)
    \class[minimum width = width("3") + 0.5em, "4"](7,0)
    \class[minimum width = width("3") + 0.5em, "4"](-9,0)
    \class[minimum width = width("3") + 0.5em, "5"](15,0)
    
    \foreach \x in {-10,...,10} {
        \begin{scope}[xshift = 4*\x]
        \foreach \y in {0,...,20}{
            \class(\y, \y + 1) }
        \end{scope}}
        
    \foreach \y in {1,...,10}{
        \class[fill](-2*\y, 2*\y) }

    \classoptions[class labels = {below = 0.1em}, "u^{-2} \eta^2"](-2,2)
    \classoptions[class labels = {below = 0.1em}, "z"](-3,0)


    \foreach \x in {0,...,0} {
        \d[red]2(-2 - 4*\x, 2 + 4*\x)
        \d[red]4(-3 - 4*\x, 4*\x)
    }
    
    \foreach \x in {0,...,1} {
        \begin{scope}[xshift = 8*\x]
        \foreach \y in {0,...,11}{
            \d[red]3(3 + \y, \y)
            \d[red]3(5 + \y, \y)
            }
        \end{scope}
    }

    \foreach \x in {-2,...,-1} {
        \begin{scope}[xshift = 8*\x]
        \foreach \y in {0,...,11}{
            \d[red]3(3 + \y - 4*\x, \y - 4*\x - 2)
            \d[red]3(3 + \y, \y)
            }
        \end{scope}
    }

\end{sseqdata}

\begin{figure}
    \centering

    \printpage[x range = {-5}{15}, y range = {0}{5}, grid = chess, page = 2--4, name = mu2SS]
    
    \caption[HFPSS for $(KU^{h(1 + 4 \Z[2]})^{hC_2}$]{The HFPSS for $C_2$ acting on $KU^{h(1 + 4 \Z[2])}$. Solid classes are detected in the HFPSS for $KU_2$ by the map $KU_2^{h(1 + 4\Z[2])} \to KU_2$, with a filtration shift of one in internal degrees $t \neq 0$. Unfilled classes are in the image of the HFPSS for $KU_2$ under $\Sigma^{-1} KU_2 \to KU_2^{1+4\Z[2]}$. The only differentials not immediately determined by this are the $d_2$ and $d_4$ differentials on classes in internal degree $t = 0$ and $-3$ respectively, which are treated at the end of \cref{lem:differentials in ASS p=2 app}.}
    \label{fig:mu_2 spectral sequence for KU^1+4 Z_2}
\end{figure}

To conclude, we must show that the two spectral sequences
        \begin{align}
            H^*(\G, \pi_* KU_2) &\implies \pi_* \bm 1_{\K} \label{eqn:first spectral sequence for L_K S}\\
            H^*(C_2, \pi_* KU_2^{1+4\Z[2]}) &\implies \pi_* \bm 1_{\K} \hspace{50pt} \label{eqn:second spectral sequence for L_K S}
        \end{align}
are isomorphic, up to a shift in filtration.

\begin{lem} \label{lem:isomorphism between two spectral sequences n=1 p=2}
For each $s \geq 0$ and $t$ there are isomorphisms
\begin{align*}
    \phantom{\qquad (t \neq 1)} H^{s+1}(\G, \pi_t KU_2) &\simeq H^{s}(C_2, \pi_{t-1} KU^{1+4\Z[2]}_2) \qquad (t \neq 1), \\
    H^s(\G, \pi_0 KU_2) &\simeq H^s(C_2, \pi_0 KU_2).
\end{align*}
These are compatible with differentials, and together yield a (filtration-shifting) isomorphism of spectral sequences between \cref{eqn:first spectral sequence for L_K S,eqn:second spectral sequence for L_K S}.
\end{lem}

\begin{remark}
In other words, when passing from the HFPSS for $KU_2^{1+4\Z[2]}$ \cref{eqn:second spectral sequence for L_K S} to the descent spectral sequence \cref{eqn:first spectral sequence for L_K S}, the filtration shift is precisely by one away from internal degree $t = 0$, and zero otherwise.
\end{remark}

\begin{proof}
Note first that the $E_1$-pages are abstractly isomorphic:
\begin{equation} \label{eqn:isomorphism between E_1 pages of two spectral sequences for n=1 p=2}
    H^{i+j}(\G, \pi_{2t} KU_2) \simeq H^i (C_2, H^j(H, \pi_{2t} KU_2)) \simeq H^i(C_2, \pi_{2t-j} KU_2^{1+ 4\Z[2]}),
\end{equation}
where $j = 1$ unless $t = 0$, in which case we also have $j = 0$. The first isomorphism is given by the Lyndon-Hochschild-Serre spectral sequence, which collapses with each degree of the abutment concentrated in a single filtration; the second isomorphism comes from the same fact about the homotopy fixed point spectral sequence for the action of $1 + 4\Z[2]$ on $KU_2$. Note that in both cases the abutment is sometimes in positive filtration, and this will account for the shift.
    
Next recall that the descent spectral sequence comes from global sections of the Postnikov tower for $\mathcal E \in \Shh[\mathcal Sp]{\proet \G}$; on the other hand, the $C_2$-fixed points spectral sequence is induced by global sections of the sheaf $j_* \mathcal E \in \Shh[\mathcal Sp]{\proet{(C_2)}}$, where $j: \G \to C_2$ is the quotient by the open subgroup $1 + 4\Z[2]$. Since each $j_* \tau_{\leq t} \mathcal E(T) = \tau_{\leq t} \mathcal E (\res[C_2]{\G} T)$ is $t$-truncated, we obtain a map of towers in $\Shh[\mathcal Sp]{\proet{(C_2)}}$,
    
\[\begin{tikzcd}
	{j_* \mathcal E} & \cdots & {\tau_{\leq t} j_*\mathcal E} & {\tau_{\leq t-1} j_*\mathcal E} & \cdots \\
	{j_* \mathcal E} & \cdots & {j_* \tau_{\leq t} \mathcal E} & {j_* \tau_{\leq t-1} \mathcal E} & \cdots
	\arrow[from=1-1, to=2-1, equals]
	\arrow[from=1-1, to=1-2]
	\arrow[from=1-2, to=1-3]
	\arrow[from=1-3, to=1-4]
	\arrow[from=1-4, to=1-5]
	\arrow[from=2-1, to=2-2]
	\arrow[from=2-2, to=2-3]
	\arrow[from=2-3, to=2-4]
	\arrow[from=2-4, to=2-5]
	\arrow[from=1-3, to=2-3]
	\arrow[from=1-4, to=2-4]
\end{tikzcd}\]
    
Taking global sections of the top row yields the HFPSS \cref{eqn:second spectral sequence for L_K S}, while the bottom yields the descent spectral sequence \cref{eqn:first spectral sequence for L_K S} (bearing in mind that $\Gamma(C_2/C_2, j_*(-)) \simeq \Gamma(\G/\G, -)$).
    
To proceed, we compute the sections of these towers. Note that any cover of $C_2 \in \proet{(C_2)}$ must split, and so sheafification on $\proet{(C_2)}$ preserves any multiplicative presheaf when restricted to the generating sub-site $\Free{C_2}$, i.e. any $\mathcal F$ satisfying $\mathcal F(S \sqcup S') \simeq \mathcal F(S) \times \mathcal F(S')$. Thus
    \[\tau_{t} j_* \mathcal E: C_2 \mapsto \tau_t \mathcal E(\Z[2]^\times/1 + 4\Z[2]) = \tau_t KU_2^{1+4\Z[2]} \]
This sheaf is therefore zero unless $t$ is odd or zero. On the other hand,
    \[ \pi_s \Gamma j_* \tau_t \mathcal E: C_2 \mapsto \pi_s \Gamma (C_2, j_* \tau_t \mathcal E) = \pi_s \Gamma(\Z[2]^\times/1 + 4\Z[2], \tau_t \mathcal E) = H^{t-s}(1+4 \Z[2], \pi_t KU_2). \]
This is zero unless $t$ is even and $s = 1$, or $t = s =0$. In particular, for $t \neq 0$ we have
    \[ \tau_{2t-1} j_* \mathcal E \simeq \tau_{2t-1} j_* \tau_{2t} \mathcal E \simeq j_* \tau_{2t} \mathcal E, \]
while $j_* \tau_0 \mathcal E$ has homotopy concentrated in degrees $\{-1, 0\}$. Note that this also implies that $\Gamma j_* \tau_{\leq 2t} \mathcal E$ is ($2t-1$)-truncated for $t \neq 0$, since $\pi_{2t} \Gamma j_* \tau_{\leq 2t} \mathcal E = \pi_{2t} \Gamma j_* \tau_{2t} \mathcal E = 0$.


The isomorphisms (\ref{eqn:isomorphism between E_1 pages of two spectral sequences for n=1 p=2}) can be interpreted as arising from the two trigraded spectral sequences
    \begin{equation}  \label{eqn:trigraded spectral sequences for n=1 p=2}
        \begin{aligned}
        H^i (C_2, H^j(1 + 4\Z[2], \pi_{2t}KU_2)) &\implies H^{i+j}(\G, \pi_{2t} KU_2), \\
        H^i (C_2, H^j(1 + 4\Z[2], \pi_{2t}KU_2)) &\implies H^i(C_2, \pi_{2t-j} KU_2^{h(1 + 4\Z[2])})
        \end{aligned}
    \end{equation}
coming from the bifiltration $\Gamma \mathcal E = \Gamma \tau_{\leq j} j_* \tau_{\leq 2t} \mathcal E$ (c.f. \cref{lem:decalage of a bifiltered spectrum appendix}). At a fixed $t$, the first is associated to the filtration $\Gamma \tau_{2t} \mathcal E = \Gamma j_* \tau_{2t} \mathcal E = \lim_j \Gamma \tau_{\leq j} j_* \tau_{2t} \mathcal E$, or equivalently to the 
resolution of $H^*(1 + 4 \Z[2], \pi_{2t} KU_2)$ by acyclic $C_2$-modules, and is therefore the LHSSS. On the other hand, the second is $C_2$-cohomology applied pointwise to the HFPSS for the $1 + 4\Z[2]$-action (at a fixed $j$); in particular, both collapse at $E_1$.

For $t \neq 0$, the towers therefore look as in the following diagram, in which we have identified in both rows those consecutive layers with zero fibre, i.e. we run the tower `at double speed'.
    
\[\begin{tikzcd}[column sep = small, row sep = small]
	\cdots & {\tau _{\leq 2t-1} j_*\mathcal E} && {\tau _{\leq 2t-3} j_*\mathcal E} & \cdots \\
	&& {\tau _{2t-1} j_*\mathcal E} \\
	\cdots & {j_* \tau_{\leq 2t} \mathcal E} && {j_* \tau_{\leq 2t-2} \mathcal E} & \cdots \\
	&& {j_* \tau_{2t} \mathcal E}
	\arrow[from=1-1, to=1-2]
	\arrow[from=1-2, to=1-4]
	\arrow[from=1-4, to=1-5]
	\arrow[from=3-1, to=3-2]
	\arrow[from=3-2, to=3-4]
	\arrow[from=3-4, to=3-5]
	\arrow[from=1-2, to=3-2]
	\arrow[from=1-4, to=3-4]
	\arrow[from=2-3, to=4-3, crossing over]
	\arrow[from=3-4, to=4-3]
	\arrow[from=4-3, to=3-2]
	\arrow[from=2-3, to=1-2]
	\arrow[from=1-4, to=2-3]
\end{tikzcd}\]

A large but routine diagram verifies that the map $\tau_{2t-1} j_* \mathcal E \to j_* \tau_{2t} \mathcal E$ agrees on homotopy with \cref{eqn:isomorphism between E_1 pages of two spectral sequences for n=1 p=2}.

Near zero, we have instead the diagram 
    
\[\begin{tikzcd}[column sep = small, row sep = small]
	{\tau _{\leq 1} j_*\mathcal E} && {\tau _{\leq 0 } j_*\mathcal E} && {\tau _{\leq -1} j_*\mathcal E} && {\tau _{\leq -2} j_*\mathcal E} \\
	& {\tau_{1} j_* \mathcal E} && {\tau_{0} j_* \mathcal E} && {\tau_{-1} j_* \mathcal E} \\
	{j_* \tau_{\leq 1} \mathcal E} && {j_* \tau_{\leq 0} \mathcal E} && {j_* \tau_{\leq -1} \mathcal E} && {j_* \tau_{\leq -2}\mathcal E} \\
	& 0 && {j_* \tau_0 \mathcal E} && 0
	\arrow[from=1-3, to=3-3]
	\arrow[from=1-3, to=1-5]
	\arrow[from=3-3, to=3-5]
	\arrow[from=1-5, to=3-5]
	\arrow[from=3-5, to=3-7, equals]
	\arrow[from=1-5, to=1-7]
	\arrow[from=2-4, to=4-4, crossing over]
	\arrow[from=3-5, to=4-4]
	\arrow[from=2-4, to=1-3]
	\arrow[from=2-6, to=4-6, crossing over]
	\arrow[from=2-6, to=1-5]
	\arrow[from=4-6, to=3-5]
	\arrow[from=3-7, to=4-6]
	\arrow[from=1-7, to=3-7]
	\arrow[from=1-7, to=2-6]
	\arrow[from=1-1, to=1-3]
	\arrow[from=3-1, to=3-3, equals]
	\arrow[from=4-4, to=3-3]
	\arrow[from=3-3, to=4-2]
	\arrow[from=4-2, to=3-1]
	\arrow[from=1-1, to=3-1]
	\arrow[from=2-2, to=4-2, crossing over]
	\arrow[from=1-3, to=2-2]
	\arrow[from=2-2, to=1-1]
	\arrow[dashed, from=2-6, to=4-4, crossing over, bend left = 40pt, "\beta"]
	\arrow[from=1-5, to=2-4, crossing over]
\end{tikzcd}\]

To see that the dashed arrow exists, note that
    \[ d_1: \tau_{-1} j_* \mathcal E \to \Sigma \tau_0 j_* \mathcal E \]
is zero on homotopy: in the proof \cref{lem:differentials in ASS p=2 app} we computed the only nontrivial $d_2$ differentials, which have source in internal degree $t = 0$. As a map between Eilenberg-Mac Lane objects, it is in fact null, so we can lift as below:
\[\begin{tikzcd}[ampersand replacement=\&, column sep = small]
	{\tau_{-1}j_* \mathcal E} \&\& {\Sigma \tau_0 j_* \mathcal E} \\
	\& {\tau_{\leq -1}j_* \mathcal E} \\
	{\tau_{\leq 0}j_* \mathcal E}
	\arrow["{\exists \alpha}"', dashed, from=1-1, to=3-1]
	\arrow[from=3-1, to=2-2]
	\arrow[from=1-1, to=2-2]
	\arrow[from=2-2, to=1-3]
	\arrow["{d_1}"', "0", from=1-1, to=1-3]
\end{tikzcd} \qquad 
\begin{tikzcd}[ampersand replacement=\&, column sep = small]
	{\tau_{-1}j_* \mathcal E} \&\&\& 0 \\
	\& {\tau_{\leq 0}j_* \mathcal E} \\
	{j_* \tau_{0} \mathcal E} \&\& {j_* \tau_{\leq 0} \mathcal E} \& {j_* \tau_{\leq -1} \mathcal E}
	\arrow["\alpha", from=1-1, to=2-2]
	\arrow[from=2-2, to=3-3]
	\arrow[from=3-1, to=3-3]
	\arrow["{\exists \beta}"', dashed, from=1-1, to=3-1]
	\arrow[from=3-3, to=3-4]
	\arrow[from=1-1, to=1-4]
	\arrow[from=1-4, to=3-4]
\end{tikzcd}\]

We must show that evident map $\tau_0 j_* \mathcal E \to j_* \tau_0 \mathcal E$ and the map $\beta$ agree with the respective compositions of edge maps. To deduce this, it is enough to contemplate the diagrams below, in which the dashed arrow are equivalences and the dotted arrows admit right inverses.
\begin{align*}
\begin{tikzcd}[ampersand replacement=\&,column sep=tiny]
	{\tau_{0}j_*\mathcal E} \&\&\& {\tau_{\leq 0} j_*\mathcal E} \\
	{\tau_{0}j_* \tau_{\geq 0}\mathcal E} \& {\tau_{0} j_* \tau_{\leq 0}\mathcal E} \\
	{\tau_{0}j_* \tau_0 \mathcal E} \&\& {\tau_{\leq 0}j_* \tau_{\leq 0}\mathcal E} \\
	\& {\tau_{\leq 0}j_* \tau_{p}\mathcal E} \\
	{j_* \tau_0 \mathcal E} \&\&\& {j_* \tau_{\leq 0} \mathcal E}
	\arrow[curve={height=30pt}, squiggly, from=1-1, to=3-1]
	\arrow[squiggly, from=3-1, to=5-1]
	\arrow[from=1-1, to=2-2]
	\arrow[from=3-1, to=2-2]
	\arrow[from=3-1, to=4-2]
	\arrow[dashed, from=5-1, to=4-2]
	\arrow[from=1-4, to=5-4]
	\arrow[from=5-1, to=5-4]
	\arrow[dashed, from=5-4, to=3-3]
	\arrow[from=1-4, to=3-3]
	\arrow[from=2-2, to=3-3]
	\arrow[from=4-2, to=3-3]
	\arrow[from=1-1, to=1-4]
	\arrow[dashed, from=2-1, to=1-1]
	\arrow[from=2-1, to=3-1]
\end{tikzcd}
    \qquad
\begin{tikzcd}[ampersand replacement=\&, column sep=tiny]
	{\tau_{-1}j_*\mathcal E} \&\&\& {\tau_{\leq 0} j_*\mathcal E} \\
	{\tau_{-1}j_* \tau_{\geq 0} \mathcal E} \& {\tau_{-1} j_* \tau_{\leq 0} \mathcal E} \\
	{\tau_{-1}j_* \tau_0 \mathcal E} \&\& {\tau_{\leq -1} j_* \tau_{\leq 0} \mathcal E} \\
	\& {\tau_{\leq -1}j_* \tau_0 \mathcal  E} \\
	{j_* \tau_0 \mathcal E} \&\&\& {j_* \tau_{\leq 0} \mathcal E}
	\arrow[curve={height=40pt}, squiggly, from=1-1, to=3-1]
	\arrow[squiggly, from=3-1, to=5-1]
	\arrow[from=1-4, to=5-4]
	\arrow[from=5-1, to=5-4]
	\arrow[from=1-1, to=1-4]
	\arrow[dashed, from=2-1, to=1-1]
	\arrow[from=2-1, to=3-1]
	\arrow[from=3-1, to=4-2]
	\arrow[dotted, from=5-1, to=4-2]
	\arrow[from=3-1, to=2-2]
	\arrow[from=1-1, to=2-2]
	\arrow[from=2-2, to=3-3]
	\arrow[from=4-2, to=3-3]
	\arrow[dotted, from=5-4, to=3-3]
	\arrow[from=1-4, to=3-3]
\end{tikzcd}
\end{align*}
The squiggly arrows are the edges maps from the trigraded spectral sequences, which are isomorphisms thanks to the collapse of the two trigraded spectral sequences.
\end{proof}